\documentclass[aap]{imsart}

\usepackage[utf8]{inputenc} 
\usepackage[T1]{fontenc}    
\usepackage{fancyref}
\usepackage{hyperref}       
\usepackage{url}            
\usepackage{booktabs}       

\usepackage{amsfonts}       
\usepackage{amsmath}
\usepackage{amsthm}
\usepackage{amssymb}
\usepackage{appendix}
\usepackage{bbm}
\usepackage{bm} 
\usepackage{commath}
\usepackage{comment}
\usepackage{enumitem}
\usepackage{float}
\usepackage{graphicx}
\usepackage{lipsum}
\usepackage{mathrsfs}
\usepackage{mathtools}
\usepackage{microtype}      
\usepackage{nicefrac}       
\usepackage{thmtools}
\usepackage{xcolor}
\usepackage{mathrsfs}

\graphicspath{ {./images/} }

\usepackage{listings}
\usepackage{nameref}
\usepackage[section]{placeins}
\usepackage{soul}
\usepackage[normalem]{ulem}
\usepackage{xcolor}

\theoremstyle{plain}
\newtheorem{theorem}{Theorem}
\newtheorem{proposition}{Proposition}
\newtheorem{lemma}{Lemma}
\newtheorem{corollary}{Corollary}

\theoremstyle{definition}
\newtheorem*{definition}{Definition}

\newtheorem{assumption}{Assumptions}

\newcommand{\rd}{\mathrm{d}}
\newcommand{\N}{\mathbb{N}}

\newcommand{\R}{\mathbb{R}}

\newcommand{\T}{\mathbb{T}}
\newcommand{\E}{\mathbb{E}}

\newcommand{\defeq}{\vcentcolon=}
\newcommand{\eqdef}{=\vcentcolon}

\newcommand{\wsto}{\overset{\ast}{\rightharpoonup}}

\startlocaldefs

\endlocaldefs

\begin{document}
\sloppy

\begin{frontmatter}
\title{Dense networks of integrate-and-fire neurons: Spatially-extended mean-field limit of the empirical measure}
\runtitle{Dense networks of integrate-and-fire neurons}

\begin{aug}
\author[A]{\fnms{Pierre-Emmanuel}~\snm{Jabin}\ead[label=e1]{pejabin@psu.edu}\orcid{0000-0001-5998-0347}},
\author[B]{\fnms{Valentin}~\snm{Schmutz}\ead[label=e2]{v.schmutz@ucl.ac.uk}\orcid{0000-0002-0935-6121}}
\and
\author[C]{\fnms{Datong}~\snm{Zhou}\ead[label=e3]{datong.zhou@sorbonne-universite.fr}}
\address[A]{Department of Mathematics, Pennsylvania State University\printead[presep={,\ }]{e1}}

\address[B]{UCL Queen Square Institute of Neurology, University College London\printead[presep={,\ }]{e2}}

\address[C]{Laboratoire de Probabilités, Statistique et Modélisation, Sorbonne Université\printead[presep={,\ }]{e3}}
\end{aug}

\begin{abstract}
The dynamics of spatially-structured networks of $N$ interacting stochastic neurons can be described by deterministic population equations in the mean-field limit. While this is known, a general question has remained unanswered: does synaptic weight scaling suffice, by itself, to guarantee the convergence of network dynamics to a deterministic population equation, even when networks are not assumed to be homogeneous or spatially structured? In this work, we consider networks of stochastic integrate-and-fire neurons with arbitrary synaptic weights satisfying a $\mathcal{O}(1/N)$ scaling condition. Borrowing results from the theory of dense graph limits, or graphons, we prove that, as $N\to\infty$, and up to the extraction of a subsequence, the empirical measure of the neurons' membrane potentials converges to the solution of a spatially-extended mean-field partial differential equation (PDE). Our proof requires analytical techniques that go beyond standard propagation of chaos methods. In particular, we introduce a weak metric that depends on the dense graph limit kernel and we show how the weak convergence of the initial data can be obtained by propagating the regularity of the limit kernel along the dual-backward equation associated with the spatially-extended mean-field PDE. Overall, this result invites us to reinterpret spatially-extended population equations as universal mean-field limits of networks of neurons with $\mathcal{O}(1/N)$ synaptic weight scaling. 
\end{abstract}



\end{frontmatter}

\tableofcontents


\section{Introduction} \label{sec:introduction}

\subsection{General motivation}
Spatially-structured population equations were among the first mathematical models proposed to describe the macroscopic dynamics of large networks of neurons in the brain \cite{Beu56, Gri63}. Fifty years ago, Wilson and Cowan derived a famous integro-differential equation approximating the collective dynamics of spatially-structured networks of integrate-and-fire-type neurons \cite{WilCow73}. Motivated by early neuronal recordings in sensory cortices \cite{Mou57,HubWie62}, their model was based on the idea that the cortex is functionally organized as a sheet of radially redundant neurons called cortical columns \cite{Mou97}. While this view is biologically simplistic from today's standpoint, spatially-structured population equations---or neural field models---have established themselves as classical models in neuroscience \cite{Nun74, Ama77, BenBar95, Erm98, Bre11, Bre17}.

Integro-differential equations such as that of Wilson and Cowan \cite{WilCow73} are not exact mean-field limits of spatially-structured networks of integrate-and-fire-type neurons because their derivation relied on a form of coarse-graining in time. The first exact mean-field limit equation that does not neglect the fine time structure of integrate-and-fire-type neurons appears to have been derived by Gerstner \cite{Ger95} and it takes the form of a spatially-structured partial differential equation (PDE). It is therefore known in theoretical neuroscience that spatially-structured networks can be described, as the number of neurons $N$ tends to infinity and if synaptic weights scale as $\mathcal{O}(1/N)$, by a spatially-structured PDE. A much less trivial question is following: In the absence of any prescribed spatial structure in the synaptic weights, does $\mathcal{O}(1/N)$ weight scaling suffice, by itself, to give rise to spatially-structured population dynamics in the mean-field limit? If the answer is positive, the spatial structure in the mean-field limit has to be interpreted as an emergent property that does not need to be related to the location of the neurons in putative cortical columns. Recent applications of the theory of dense graph limits \cite{LovSze06,Lov12}, or graphons, to the study of mean-field limits of non-exchangeable systems \cite{KalMed18,JabPoy24} suggest that the answer to the question is indeed positive.

The goal of this work is to present a transparent proof for the convergence of non-exchangeable systems of interacting integrate-and-fire-type neurons to a spatially-extended PDE in the mean-field limit. In order to do so, we will consider the mean-field limit of dense networks only (instead of the networks with intermediate sparsity considered previously in \cite{JabZho23}). While the spatial extension of the mean-field PDE we obtain resembles that of classical models such as neural field models, its interpretation is radically different: The spatial extension does not come from any prescribed spatial structure, but it reflects the organization of the emergent redundancies of individual neuron trajectories induced by $\mathcal{O}(1/N)$ synaptic weight scaling in the mean-field limit.

\subsection{Model: Non-exchangeable systems of interacting neurons}\label{sec:non-exchangeable}
We consider networks of biological neurons where neurons are modeled as integrate-and-fire neurons with escape noise \cite{Ger00}, a network model which has already been studied in several mathematical works, see, e.g.,~\cite{DemGal15, FouLoe16, CorTan20, CorTan21, SchLoe23}. 

For any network size $N$, each neuron $i \in \{1, \dots, N\}$ has a variable $\bm{X}^{i;N}(t) \in \R$ that represents its membrane potential. The membrane potential of a neuron $\bm{X}^{i;N}(t)$ determines its instantaneous probability of emitting a spike, that is, its conditional intensity $f(\bm{X}^{i;N}(t-))$, through the non-negative and non-decreasing intensity function $f:\R \to \R_+$. Each time neuron $i$ emits a spike, its membrane potential is reset to $0$ and the membrane potentials of all other neurons $j \neq i$ make a jump of $\frac{1}{N}w^N_{j,i}$, respectively, where $w^N_{j,i}$ denotes the (rescaled) synaptic weight from neuron $i$ to neuron $j$ (by convention, we put $w^N_{i,i}=0$ for all $i$). Between spikes, the membrane potential $\bm{X}^{i;N}(t)$ drifts according to the velocity field $b : \R \to\R$. 

More formally, given the initial conditions $\{\bm{X}^{i;N}(0)=\bm{X}^{i;N}_0\}_{i=1}^N$, the stochastic dynamics of the network model is characterized by the following system of stochastic differential equations (SDEs): For $i \in \{1, \dots, N\}$,
\begin{equation} \label{eqn:system}
\begin{aligned}
\rd \bm{X}^{i;N}(t) & = b(\bm{X}^{i;N}(t)) \rd t + \frac{1}{N}\sum_{j=1}^N w^N_{i,j}\rd \bm{Z}^{j;N}(t) - \bm{X}^{i;N}(t-)\rd \bm{Z}^{i;N}(t),
\\
\bm{Z}^{i;N}(t) & \defeq \int_{[0,t]\times \R_+}\mathbbm{1}_{\{z\leq f(\bm{X}^{i;N}(s-))\}}\bm{\Pi}^{i;N}(\rd s, \rd z),
\\
\bm{X}^{i;N}(0) & = \bm{X}^{i;N}_0,
\end{aligned}
\end{equation}
where the $\{\bm{\Pi}^{i;N}(\rd s, \rd z)\}_{i=1}^N$ are independent Poisson random measures on $\R_+ \times \R_+$ with intensity $\rd s \rd z$\footnote{See \cite{Kin92, Loe17} for introductions on Poisson random measures.} and the initial data $\{\bm{X}^{i;N}_0\}_{i=1}^N$ are $N$ real-valued random variables.

We consider the pathwise unique $(\mathscr{F}_t)_{t\in\R_+}$-adapted càdlàg strong solution to \eqref{eqn:system} in the filtered probability space $\left(\Omega, \mathscr{F}, (\mathscr{F}_t)_{t\in\R_+}, \mathbb{P}\right)$, where the Poisson random measures $\bm{\Pi}^{i;N}$ are independent of $\mathscr{F}_0 = \sigma(\{\bm{X}^{i;N}_0\}_{i=1}^N)$ and where $\mathscr{F}_t = \mathscr{F}_0 \cup \sigma\left(\{\bm{\Pi}^{i}\}_{i \in \mathbb{N}^*, B \in \mathscr{B}(\R_+)}\right)$ ($\mathscr{B}(\R_+)$ denotes the Borel $\sigma$-algebra on $\R_+$). Furthermore, we work with the following assumptions on the intensity function $f$ and the velocity field $b$.

\begin{assumption}\label{assumption:b_f} 
The intensity function $f :\R \to \R_+$ is bounded and continuously differentiable; $f\in C^1_b(\R,\R_+)$. The velocity field $b : \R \to \R$ is bounded and continuously differentiable; $b\in C^1_b(\R)$.
\end{assumption}

Under Assumptions~\ref{assumption:b_f}, the system of SDEs~\eqref{eqn:system} defines a Piecewise-Deterministic Markov Process with bounded variation: The membrane potentials follow the deterministic dynamics $\frac{\rd}{\rd t} \bm{X}^{i;N}(t) = b(\bm{X}^{i;N}(t))$ between spikes and the probability of spike emission is bounded by $\|f\|_\infty < +\infty$. 


This work focuses on the behavior of the system~\eqref{eqn:system} as $N \to \infty$, when synaptic weights scale as $\mathcal{O}(1/N)$. For clarity, in \eqref{eqn:system}, the $\mathcal{O}(1/N)$ scaling of synaptic weights is explicitly indicated by the factor $1/N$. Consistently, the rescaled synaptic weights $w^N_{i,j}$ will be assumed to be of order $\mathcal{O}(1)$ as $N \to \infty$. The goal of this work is to describe the mean-field limits of the system \eqref{eqn:system} when no particular assumptions are made on the rescaled synaptic weights $w^N_{i,j}$. As explained in \cite{JabPoy24}, because of the $\mathcal{O}(1/N)$ scaling of interactions between neurons, the stochastic multi-neuron system \eqref{eqn:system} is expected to effectively behave as a deterministic system based on the neurons' law (or density function) as $N\to\infty$; this intuition can be seen as a generalization of standard \textit{propagation of chaos} arguments, which we recall in Section~\ref{sec:exchangeable}. The general non-exchangeable system~\eqref{eqn:system} can be reduced to an exchangeable system by assuming that the initial data $\{\bm{X}^{i;N}_0\}_{i=1}^N$ are \textit{i.i.d.} and rescaled synaptic weights $w^N_{i,j}$ are also \textit{i.i.d.}; in this case, the mean-field limit of the system can be described by a mean-field PDE with no spatial extension, as proved in \cite{DemGal15,FouLoe16}.

To give an intuition for the type of mean-field limit we can expect in the case of general non-exchangeable systems, let us consider the example of \textit{growing uniform attachment graphs} (see \cite{Lov12}) as an example of sequences of random synaptic weights giving rise to non-exchangeable systems. In this example, the synaptic weights $w^{N}_{i,j}$ are generated via an iterative procedure. The procedure starts at $N=1$ with a single neuron. At the $N+1$-th iteration step, a neuron is added to the network of the previous step and is randomly and independently connected to each of the $N$ previous neurons with probability $1/(N+1)$. Two neurons are said to be connected if $w_{i,j} = w_{j,i} = 1$ and not connected if $w_{i,j} = w_{j,i} = 0$ (synaptic weights are symmetric). This construction leads to a sequence of non-exchangeable systems since neurons added at the beginning of the procedure tend to have more connections than neurons added at the end. Although the algorithm for generating growing uniform attachment graphs does not rely on an explicit spatial structure, the order in which neurons are added actually defines an implicit spatial structure which makes these graphs equal in law to $W$-random graphs. Indeed, as explained in \cite[Example 11.39]{Lov12}, for any network size $N$, if the indices $i$ of neurons correspond to the iteration step at which they were added to the network, then each neuron $i$ can be assigned to a location $\xi_i = (i-1)/N$ on the interval $[0,1]$. Then, one can verify that the connections between pairs of distinct neurons are independent $\{0,1\}$-valued random variables with expectation $\E[w^{N}_{i,j}]=1 - \max(\xi_i,\xi_j)$, for all $i\neq j$.

The theory of graphons \cite{Lov12}, which will be briefly reviewed in Sec.~\ref{sec:graphon}, tells us that growing uniform attachment graphs converge, as $N\to\infty$, to the limit graphon $(\xi,\zeta) \mapsto 1 - \max(\xi,\zeta) \in L^\infty([0,1]^2)$, in the topology induced by the \textit{cut norm} (see \cite[Proposition 11.40]{Lov12}). This suggests that, in this example, the dynamics of the system~\eqref{eqn:system} could converge, in the mean-field limit, to the solution to a deterministic equation with a spatial extension given by the interval $[0,1]$. In this example, neurons can be ordered on the interval $[0,1]$ in a way that reveals the implicit spatial structure of the model. However, one can imagine sequences of dense graphs with no apparent spatial structure. In the latter case, the theory of graphons tells us that \textit{any} sequence of dense graphs converges, up to the extraction of a subsequence, to a limit graphon $w\in L^\infty([0,1]^2)$, in the topology of the \textit{cut distance}, which is the infimum of the cut norm over all possible orderings of the nodes (neurons) on the interval $[0,1]$ (see Sec.~\ref{sec:graphon} for precise statements). This highly nontrivial result from the theory of graphons suggests that, in some appropriate topology, the dynamics of \textit{any} sequence of dense networks~\eqref{eqn:system} has to converge---up to neuron re-orderings and up to the extraction of a subsequence---to the solution to a deterministic equation with a spatial extension on the interval $[0,1]$ and a limit kernel $w\in L^\infty([0,1]^2)$. Therefore, to fully exploit the theory of graphons and obtain such a general result for non-exchangeable systems, we will consider mean-field limits up to neuron re-orderings and up to the extraction of a subsequence.

From a technical point of view, our task will be to identify a topology in which the convergence of the synaptic weights to the limit kernel $w$, as $N\to\infty$, can be well propagated to the dynamics of the system as a whole. 

\subsection{Main result: Mean-field limits of non-exchangeable systems}

As mentioned in the previous section, to leverage the theory of graphons and obtain a result for general non-exchangeable systems, we need to study mean-field limits up to neuron re-orderings. To allow neurons to be arbitrarily re-ordered, we assign to each neuron $i\in\{1, \dots, N\}$ a measurable set of Lebesgue measure $1/N$ on the interval $[0,1]$ such that the $N$ sets form a partition of $[0,1]$. More specifically, we will use the following notion of \textit{almost everywhere partition}.

\begin{definition}
A set of measurable sets $\{E^{i;N}\}_{i = 1}^N \subset [0,1]^N$ is an \textit{almost everywhere partition} of $\{1,\dots,N\}$ if, for all $1 \leq i \leq N$, $E^{i;N}$ has Lebesgue measure $1/N$ and the identity $\mathbbm{1}_{[0,1]} = \sum_{i = 1}^N \mathbbm{1}_{E^{i;N}}$ holds almost everywhere on $[0,1]$.
\newline

Conveniently, almost everywhere partitions will allow us to deal with neuron re-orderings in an implicit manner. Given an almost everywhere partition $\{E^{i;N}\}_{i = 1}^N$, we define the \textit{extended empirical measure} of $\{\bm{X}^{i;N}(t)\}_{i=1}^N$ as
\begin{equation} \label{eqn:extended_empirical_measure}
\begin{aligned}
\bm{\mu}^N(t,\xi,\rd x) \defeq \sum_{i=1}^N\delta_{\bm{X}^{i;N}(t)}(\rd x)\mathbbm{1}_{E^{i;N}}(\xi), \quad \forall t\geq 0,\;\forall \xi \in [0,1],
\end{aligned}
\end{equation}
\end{definition}
where $\delta_X$ denotes the Dirac measure located in $X$. The auxiliary variable $\xi\in[0,1]$ in this definition will play the same role as the auxiliary variable $\xi$ in the theory of graphons (see Sec.~\ref{sec:graphon}). 

Informally, our main result, Theorem~\ref{thm:main}, states that, up to neuron re-orderings (which are taken care of by the almost everywhere partitions), the extended empirical measures $\bm{\mu}^N$ of the non-exchangeable system~\eqref{eqn:system} converges in the mean-field limit $N\to\infty$, in a weak-* sense and up to the extraction of a subsequence, to the solution $\mu$ to the spatially-extended mean-field PDE
\begin{subequations} \label{eqn:Vlasov}
\begin{align} 
&\partial_t  \mu(t,\xi,\rd x) + \partial_x \Big[\big(b(x) + h(t,\xi)\big) \mu(t,\xi, \rd x)\Big] + f(x) \mu(t,\xi,\rd x) - r(t,\xi)\delta_0(\rd x) = 0,
\label{eqn:Vlasov_field} \\ 
& r(t,\xi) = \int_{\R}f(x) \mu(t,\xi,\rd x), \quad h(t,\xi) = \int_{[0,1]} w(\xi,\zeta) r(t,\zeta) \rd \zeta, \label{eqn:Vlasov_r_h} \\ 
&\mu(0,\xi,\rd x) = \mu_0(\xi, \rd x), \notag
\end{align}
\end{subequations}
where $w \in L^\infty([0,1]^2)$ is a limit kernel. The PDE~\eqref{eqn:Vlasov} can be seen a spatially-extended version of the mean-field PDE rigorously derived by De Masi, Galves, Löcherbach, and Presutti~\cite{DemGal15} in the case of exchangeable systems of integrate-and-fire neurons with escape noise. In~\eqref{eqn:Vlasov_r_h}, $r(t,\xi)$ represents the firing rate of neurons at ``location'' $\xi$, and $h(t,\xi)$ represents the postsynaptic input received by neurons at location $\xi$. Importantly, the result does not require the networks~\eqref{eqn:system} to have a prescribed spatial structure as it applies to general non-exchangeable systems with $\mathcal{O}(1/N)$ synaptic weight scaling. This generality implies that, in a the mean-field limit, the dynamics of large networks of neurons with $\mathcal{O}(1/N)$ synaptic weight scaling can always be described by a deterministic PDE with spatial extension over the interval $[0,1]$.

\subsubsection*{General notations for function spaces}
The function space $C_b(\R)$ denotes the space of continuous bounded functions in $\R$, equipped with the uniform norm, and $C^k_b(\R)$ denotes the subspace of $C_b(\R)$ with continuous and bounded derivatives up to the $k$-th order. The function space $C([0,t])$ denotes the space of continuous functions on the finite time interval $[0,t]$. Lebesgue and Sobolev spaces are denoted $L^p$ and $W^{s,p}$ respectively.

We denote by $\mathcal{M}(\R)$ the space of signed Borel measures with bounded \emph{total variation norm} on $\R$, $\mathcal{M}_+(\R)$ the subspace of non-negative measures, and $\mathcal{P}(\R)$ the subspace probability measures. 
As for the choice of the topology on spaces of measures, we will use the standard notion of \emph{weak-*} topology and use $\wsto$ to denote the weak-* convergence.

\begin{theorem} \label{thm:main}
Grant Assumptions~\ref{assumption:b_f}. Let $\{\bm{X}^{i;N}\}_{i = 1}^N$, $N\to\infty$, be a sequence of solutions to \eqref{eqn:system} with synaptic weight matrices $\{w^N_{i,j}\}_{i,j = 1}^N$ satisfying the uniform boundedness condition 
\begin{equation}\label{eq:boundedness_condition}
    \sup_{N} \max_{1 \leq i,j \leq N} |w^N_{i,j}| < \infty.
\end{equation}
Moreover, assume that the initial data $\{\bm{X}^{i;N}_0\}_{i = 1}^N$ are independent random variables. Finally, assume the moment bound for the initial data
\begin{equation}\label{eq:th1_moment_bound}
\begin{aligned}
\sup_{N} \sup_{1 \leq i \leq N} \E\big[|\bm{X}^{i;N}_0|^2\big] < \infty.
\end{aligned}
\end{equation}

Then, there exist a limit kernel $w \in L^\infty([0,1]^2)$, a solution $\mu \in L^\infty([0,t_*] \times [0,1]; \mathcal{M}(\R))$ to the PDE~\eqref{eqn:Vlasov} in the sense of characteristics for all $t_* > 0$, and almost everywhere partitions $\{E^{i;N}\}_{i = 1}^N$ such that the sequence $\{\bm{\mu}^N\}$ as defined in \eqref{eqn:extended_empirical_measure} converges, up to the extraction of a subsequence, to $\mu$ in the following sense: For any $t \in[0,t^*]$,
\begin{equation}\label{eq:th1_weak-star}
\begin{aligned}
\bm{\mu}^N(t,\cdot,\cdot) \wsto \mu(t,\cdot,\cdot)\ \mbox{in}\ \mathcal{M} ([0,1] \times \R)\quad \mbox{a.s.}, \quad \text{as } N \to \infty. 
\end{aligned}
\end{equation} 
\end{theorem}
The exact definition of the weak-* convergence above is presented in Section~\ref{sec:metric}, where we also construct the topology of weak-* convergence through quantitative norms. (Using such a constructed weak norm, the convergence in \eqref{eq:th1_weak-star} can actually be re-stated as a uniform-in-time convergence on the finite interval $[0,t_*]$; see Theorem~\ref{thm:main_metric} in Sec.~\ref{sec:metric}.)

Note that the weak-* convergence $\bm{\mu}^N \wsto \mu$ also includes the weak-* convergence of the extended empirical measure at time $0$, which is nontrivial because the sequences of partitions $\{E^{i;N}\}_{i = 1}^N$ we use. Intuitively, the sequence of partitions $\{E^{i;N}\}_{i = 1}^N$ re-order the ``locations'' of the neurons on $[0,1]$ for each network size $N$, as $N\to\infty$. Therefore, we need to work with a notion of convergence that is weak enough in the $\xi$-dimension so that the extended empirical measure at time $t=0$ converges despite the re-orderings caused by the sequences of partitions $\{E^{i;N}\}_{i = 1}^N$. 

Importantly, the convergence of the extended empirical measure $\bm{\mu}^N$ in the mean-field limit, established in Theorem~\ref{thm:main}, keeps track of individual neuron trajectories $\bm{X}^{i;N}(t)$, as stated in the following corollary.

\begin{corollary} \label{cor:trajectory}

Grant all the assumptions and results of Theorem~\ref{thm:main}. For all $N\geq 1$, let $\{\widetilde{\bm{X}}^{i;N}(t)\}_{i=1}^N$ be the auxiliary processes where the postsynaptic input terms $\frac{1}{N}\sum_{j=1}^N w^N_{i,j}\rd \bm{Z}^{j;N}(t)$ are replaced by the respective mean-field inputs, namely
\begin{equation}\label{eqn:auxiliary}
\begin{aligned}
\rd \widetilde{\bm{X}}^{i;N}(t) & = b(\widetilde{\bm{X}}^{i;N}(t)) \rd t + \int_{[0,1]^2} \mathbbm{1}_{E^{i;N}}(\xi) w(\xi,\zeta) \int_{\R} f(y)\mu(t,\zeta,y) \rd y \rd\zeta \rd\xi 
\\
& \hspace{7cm} - \widetilde{\bm{X}}^{i;N}(t-)\rd \widetilde{\bm{Z}}^{i;N}(t),
\\
\widetilde{\bm{Z}}^{i;N}(t) & \defeq \int_{[0,t]\times \R_+}\mathbbm{1}_{\{z\leq f(\widetilde{\bm{X}}^{i;N}(s-))\}}\bm{\Pi}^{i;N}(\rd s, \rd z),
\\
\widetilde{\bm{X}}^{i;N}(0) & = \bm{X}^{i;N}_0.
\end{aligned}
\end{equation}
Then, for all $t_* \geq 0$,
\begin{equation*}
\begin{aligned}
& \E\bigg[\frac{1}{N} \sum_{i=1}^N  \sup_{t \in [0,t_*]} \min\Big( \big| \bm{X}^{i;N}(t) - \widetilde{\bm{X}}^{i;N}(t) \big|, 1 \Big) \bigg] \to 0, \ \text{ as } \ N \to \infty.
\end{aligned}
\end{equation*}
\end{corollary}

The fine tracking of individual neuron trajectories in the mean-field limit is one of the main advantages of the proof method presented in this work, compared to method of macroscopic ``observables'' used in \cite{JabZho23} (see Sec.~\ref{sec:previous}). Note that the auxiliary processes $\widetilde{\bm{X}}^{i;N}$ in 
Corollary~\ref{cor:trajectory} do not form a closed system since \eqref{eqn:auxiliary} involves the limit empirical measure $\mu$.

\subsection{Previous works}\label{sec:previous}
Mean-field models have played a major role in shaping our understanding of emergent behaviors in large networks of interacting neurons. We refer the interested reader to \cite{GerKis14,Tal10,Tal11,HelDah20} for recent and comprehensive textbooks covering different types of mean-field theories in neuroscience.
Although mean-field theory in neuroscience is a vast topic, in the following, we restrict our survey of previous works to exact mean-field models for networks of integrate-and-fire-type neurons, with an emphasis on the mathematical literature most relevant to the present work. 
For recent results on the replica-mean-field limit of networks of neurons, a distinct notion of limit which will not be covered here, we refer the reader to~\cite{Dav23,AviDav25} and references therein. 

Exact mean-field models for networks of integrate-and-fire-type neurons can be broadly classified into two classes, depending on how the intrinsic noise of neurons is modeled: Neuronal noise can be modeled with ``escape noise'', i.e., a soft, probabilistic firing threshold, or it can be modeled with diffusive noise in the membrane potential dynamics. Mean-field models with escape noise were introduced in theoretical neuroscience by \cite{GerHem92, Ger95, Ger00} (see also \cite[Chaps.~9~\&~14]{GerKis14}, \cite{GalLoe16} and \cite[Chap.~1.2]{Sch22thesis} for introductions). The convergence of homogeneous networks, i.e., exchangeable systems, in the mean-field limit to a nonlinear population equation was proved in \cite{DemGal15, FouLoe16}, and the dynamics of the limit equation has been analyzed in \cite{RobTou16,CorTan20,CorTan21,Cor24}. Integrate-and-fire neurons with escape noise are, from a modeling point of view, very similar to age-dependent Hawkes processes (see \cite[p.~10]{Sch22thesis}). The convergence of exchangeable systems of interacting age-dependent Hawkes processes to a mean-field limit was proved in \cite{Che17} and this convergence proof has since been generalized to the case of multiple interacting populations of neurons \cite{RaaDit20} and to the case of neurons with additional short-term memory variables \cite{Sch22}.
The dynamics of the limit mean-field equation, known as the Time Elapsed Neuron Network Model, has been analyzed in \cite{PakPer10, PakPer13, MisQui18, MisWen18,CanYol19} (see also \cite{PakPer14, TorPer22, FonSch22} for works on the dynamics of generalized versions of the original limit equation). Concurrently, mean-field models with diffusive noise were introduced in theoretical neuroscience by Brunel and collaborators \cite{BruHak99,Bru00,OstBru09}, the convergence of exchangeable systems to a mean-field limit was proved in \cite{DelIng15b,IngTal15}, and the dynamics of the limit mean-field equation---known as the Nonlinear Noisy Leaky Integrate and Fire (NNLIF) model---has been analyzed in \cite{CacCar11, CarGon13, CacPer14, CarPer15, DelIng15a}. 

The mean-field limit of spatially-structured networks of integrate-and-fire neurons, which are examples of non-exchangeable systems, has received much less attention in mathematical neuroscience (but see \cite{JabZho23}). However, many works have studied the mean-field limit of spatially-structured networks of neurons that are not of integrate-and-fire type: The case of Kuramoto oscillators was treated in \cite{Med14,Med14b,KalMed18,OliRei19}, the cases of FitzHugh-Nagumo or Hodgkin-Huxley neurons in \cite{Tou14, LucSta14, Cre19, MehSch20, Luc20}, and the case of nonlinear Hawkes processes in \cite{CheDua19,Aga22}. Our general non-exchangeable system model, presented in Sec.~\ref{sec:non-exchangeable}, differs from the models considered in the aforementioned works in two respects. First, we do not assume any prescribed spatial structure; and second, we consider integrate-and-fire-type neurons, which have discontinuous dynamics that poses specific mathematical difficulties.

Our result and proof method mainly build on three recent works. (I) Our result is a generalization of the mean-field convergence proof for exchangeable systems of interacting integrate-and-fire neurons with escape noise \cite{FouLoe16} to the case of non-exchangeable systems. Notably, we use the same probabilistic modeling framework as in \cite{FouLoe16}, which involves Poisson random measures, and our neuron model is almost identical to theirs. Also, our treatment of individual neuron trajectories in the limit (Corollary~\ref{cor:trajectory}) uses some of their methods. (II) Our result is complementary to another work by two of the authors \cite{JabZho23}, where a recent method for proving the convergence in the mean-field limit of non-exchangeable systems \cite{JabPoy24} is applied to networks of integrate-and-fire neurons. In \cite{JabZho23}, it is proven that non-exchangeable systems of interacting integrate-and-fire neurons with both diffusive noise and escape noise (soft threshold) converge to a spatially-extended mean-field limit. While the result in \cite{JabZho23} is closely related to the result presented here, there is an important difference that makes the two results complementary: In \cite{JabZho23}, the mean-field convergence is not proven through the empirical measure of the networks but through a tree-indexed hierarchy of macroscopic ``observables'' which can be seen as a generalization of the classical BBGKY hierarchy for exchangeable systems. The advantage of the approach in \cite{JabZho23} is that it enables the study the mean-field limit of networks with synaptic weight scalings that are more general than $\mathcal{O}(1/N)$ synaptic weight scaling. However, the disadvantage of the ``observables'' of \cite{JabZho23} is that, by their very definition, they destroy all information about the dynamics of individual neurons or subsets of neurons. In contrast, the spatially-extended empirical measure we consider here keeps track of the dynamics of individual neurons and subsets of neurons, offering a much more detailed and transparent view of the system as we take the limit. The cost of working directly on the empirical measure is that we have to require synaptic weights to scale as $\mathcal{O}(1/N)$, which means that our result only applies to the limit of dense networks. Interestingly, diffusive noise is necessary for the proof method in \cite{JabZho23}, whereas the proof presented here relies on the \textit{absence} of diffusive noise, which adds to the complementary nature of the two works. (III) The weak-* convergence of the empirical measure of exchangeable systems of interacting integrate-and-fire neurons was studied in \cite{FlaPri19} (using a slightly different neuron model). We will adapt the notion of weak-* convergence of the empirical measure to the case of non-exchangeable systems; in particular, we will show that in the case of non-exchangeable systems, a \textit{spatially-extended} empirical measure is the natural object to consider.

We end this section by briefly outlining the general mathematical context of the use of graph limits for the study of limit nonlinear dynamics in large systems. The two extremes of the theory of graph limits are, on the one hand, the theory of graphons for dense graphs \cite{LovSze06, Lov12}, i.e., graphs with degrees growing linearly with the number of nodes $N$ as $N\to\infty$, and, on the other hand, the Benjamini-Schramm local convergence theory for sparse graphs \cite{BenSch01}, i.e., graphs with finite degrees as $N\to\infty$. In between, there are various graph limit theories for graphs with intermediate sparsity, e.g., $L^p$ graphons \cite{BorCha19, BorCha18}, extended graphons \cite{JabPoy24}, and graphops \cite{BacSze22} (note that graphops also apply to the limit of some specific sparse graphs). These different graph limits can then be applied to the study of the limit nonlinear dynamics of systems as the number of agents or particles tends to infinity. The theory of graphons for dense graphs was first applied to the mean-field limit of deterministic systems, notably in \cite{Med14,Med14b,KalMed18,ChiMed19,PauTre22}, and then to the mean-field limit of interacting diffusions, e.g., in \cite{Luc20}, and of interacting point processes \cite{Aga22}. At the other end of the sparsity spectrum, the Benjamini-Schramm local convergence for sparse graphs was used to characterize the limit nonlinear dynamics of interacting diffusions; see, for instance, \cite{Mac18,OliRei20,LacRam21,LacRam23,LacRam23b}. For systems interacting on graphs of intermediate sparsity, the mean-field limit of deterministic systems has been studied in \cite{KalMed17, JabPoy24, GkoKue22}, of interacting diffusions in \cite{OliRei19, CopDie20}, and, as mentioned above, of networks of integrate-and-fire-type neurons in \cite{JabZho23}.

\section{Preliminaries}

\subsection{Kinetic theory formalism for the mean-field limit of the empirical measure}\label{sec:exchangeable}

Most methods currently used in mathematical neuroscience to study mean-field limits of networks of neurons have their origin in mathematical kinetic theory, a field initiated by Kac who first introduced the notion of propagation of chaos in a seminal paper \cite{Kac56}. Propagation of chaos is a pivotal notion that allows us to link the microscopic dynamics of interacting particles to mean-field PDEs describing the evolution of the density of particles, as the number of particles $N$ tends to infinity. In very simple terms, a system of interacting particles (or agents, neurons) satisfy propagation of chaos if, as $N\to\infty$, particles behave like independent processes with the same law (see \cite{ChaDie22a, ChaDie22b} for a recent and comprehensive two-part review). In the case of exchangeable systems---systems where the joint law of particles is invariant with respect to permutations of particles---mean-field limits and propagation of chaos, are practically synonymous. However, in the case of non-exchangeable systems, e.g., spatially-structured networks, there can be mean-field limits without propagation of chaos (see references in the previous sections). Here, we say that a system converges to a mean-field limit if, as $N \to \infty$, its dynamics can be described by a deterministic equation where each particle interacts with an infinite number of other particles. 

Since deterministic methods will play an important role in this work, we briefly recall below a formalism used in kinetic theory to study the convergence of the empirical measure of deterministic systems \cite{BraHep77, Dob79}. This formalism will guide our proof strategy. The following exposition is adapted from \cite{Jab14, Gol16}, where deterministic methods for the study of mean-field limits are reviewed.

Let us consider the toy model described by the following system of $N$ ordinary differential equations: For $i \in \{1, \dots, N\}$,
\begin{equation}\label{eq:toy}
\begin{aligned}
\frac{\rd}{\rd t} X^{i;N}(t) =\;& b(X^{i;N}(t)) + \frac{w_0}{N} \sum_{j=1}^N  f(X^{j;N}(t)),
\quad
X^{i;N}(0) = X^{i;N}_0.
\end{aligned}
\end{equation}
The toy model~\eqref{eq:toy} is deterministic and, in the context of neural networks, can be interpreted as a network of $N$ rate-units with uniform synaptic weights $w_0/N$. If the initial data $\{\bm{X}^{i;N}_0\}_{i=1}^N$ are assumed to be $\textit{i.i.d.}$ with law $\rho_0(\rd x)$ (in this work, random variables and random functions are in bold), the system~\eqref{eq:toy} becomes an exchangeable system, and it is a well-known result that the empirical measure
\begin{equation*}
\begin{aligned}
\bm{\rho}^N(t,\rd x) = \frac{1}{N} \sum_{i=1}^N\delta_{\bm{X}^{i;N}(t)}(\rd x), \quad \forall t > 0,
\end{aligned}
\end{equation*}
converges the the solution $\rho(t,\rd x)$ to the mean-field PDE
\begin{equation}\label{eq:toy_PDE}
\begin{aligned}
\partial_t \rho(t,\rd x) + \partial_x  \bigg[\bigg(b(x) + w_0 \int_\R f(y)\rho(t,\rd y)\bigg)\rho(t,\rd x)\bigg] = 0, 
\quad
\rho(0,\rd x) = \rho_0(\rd x).
\end{aligned}
\end{equation}

Deterministic methods offer an elegant way to quantify the convergence of the random empirical measure $\bm{\rho}^N$ to the deterministic solution $\rho$ to the mean-field PDE~\eqref{eq:toy_PDE}. First, going back to deterministic initial data $\{X^{i;N}_0\}_{i=1}^N$ in the system~\eqref{eq:toy}, let us notice that the deterministic empirical measure $\rho^N$ given by the system~\eqref{eq:toy} for initial data $\{X^{i;N}_0\}_{i=1}^N$ is already a solution to the mean-field PDE~\eqref{eq:toy_PDE}, in the sense of distributions, if the initial datum $\rho_0(\rd x)$ is replaced by the empirical measure $\rho^N_0(\rd x) := \frac{1}{N} \sum_{i=1}^N\delta_{\bm{X}^{i;N}_0}(\rd x)$. This observation suggests that proving the convergence of the empirical measure to the mean-field limit can be reduced to a deterministic problem: Can we find an appropriate weak metric on the space of probability measure $\mathcal{P}(\R)$ such that the distance between the empirical measure $\rho^N(t,\rd x)$ and the limit $\rho(t,\rd x)$ at time $t>0$ can be controlled by the distance between $\rho^N(0,\rd x)$ and $\rho_0(\rd x)$?

An example of such weak metric is the negative Sobolev norm $W^{-1,1}(\R)$, defined as 
\begin{equation*} 
\|\rho_1 - \rho_2\|_{W^{-1,1}(\R)} \defeq \sup_{ \| \phi \|_{W^{1,\infty}(\R)} \leq 1} \int_{\R} \phi(x) \Big( \rho_1(\rd x) - \rho_2(\rd x) \Big), \;\forall (\rho_1, \rho_2)\in \mathcal{M}(\R)\times\mathcal{M}(\R),
\end{equation*}
where $W^{1,\infty}(\R)$ denotes the Sobolev space with $k=1$ and $p=\infty$
and $\mathcal{M}(\R)$ denotes the space of signed Borel measures with finite total variation norm on $\R$.\setcounter{footnote}{0}\footnote{Note that this negative Sobolev norm is similar, but not equivalent, to the Wasserstein distance $W_1(\R)$
\begin{equation*} 
\|\rho_1 - \rho_2\|_{W_1(\R)} \defeq \sup_{ \| \phi \|_{{\mathrm{Lip}}(\R)} \leq 1} \int_{\R} \phi(x) \Big( \rho_1(\rd x) - \rho_2(\rd x) \Big), \qquad \forall (\rho_1, \rho_2)\in \mathcal{P}(\R)\times\mathcal{P}(\R),
\end{equation*}
which more commonly used in the study of mean-field limits; cf. \cite{Szn91,Jab14,Gol16}.} The weak norm $W^{-1,1}(\R)$ is convenient for two reasons. First, it metrizes the weak-* convergence of probability measures on $\R$, in the sense that, given a sequence of probability measures $\{\varrho^N\}$, $\|\varrho^N -\varrho\|_{W^{-1,1}(\R)} \to 0$ as $N\to\infty$ if and only if $\varrho^N \wsto \varrho$ ($\wsto$ denotes weak-* convergence) and if the sequence $\{\varrho^N\}$ is tight (which is typically guaranteed via some uniform moment bound).\footnote{For more detailed statements, we refer the reader to \cite[Sec.~3.3]{Jab14}.} Second, the weak norm $W^{-1,1}(\R)$ can be used to obtain a stability estimate for the mean-field PDE~\eqref{eq:toy_PDE},
\begin{equation}\label{eq:Dobrushin}
\begin{aligned}
\|\rho^N(t,\cdot) - \rho(t,\cdot)\|_{W^{-1,1}(\R)} \leq C(t) \|\rho^N_0 - \rho_0\|_{W^{-1,1}(\R)},
\end{aligned}
\end{equation}
where $C(t)$ is a function depending on $\|b\|_{W^{1,\infty}}$ and $\|f\|_{W^{1,\infty}}$. This stability estimate indeed implies that, taking \textit{i.i.d} initial data $\{\bm{X}^{i;N}_0\}_{i=1}^N$ with law $\rho_0$, the convergence of the random empirical measure $\boldsymbol{\rho}^N(t,\cdot)$ to the deterministic $\rho(t,\cdot)$ at time $t > 0$ is directly controlled by the convergence of the initial datum $\boldsymbol{\rho}^N_0$ to the law $\rho_0$ (which only require some moment bound). 

The main arguments for the proof of the stability estimate~\eqref{eq:Dobrushin} can be found, e.g., in \cite[Theorem~2.2]{Jab14}. Crucially, the proof involves the propagation of the regularity of the test function $\bar{\phi}\in W^{1,\infty}(\R)$ along the dual-backward equation
\begin{equation*}
\begin{aligned}
\partial_t \phi(t,x) + \left(b(x) + w_0 \int_\R f(y)\rho(t,\rd y)\right)\partial_x  \phi(t, x) = 0, 
\quad
\phi(t, x) = \bar{\phi}(x).
\end{aligned}
\end{equation*}

To some extent, the proof of our main result is inspired by this formalism. Notably, the proof will also require a propagation of regularity along the dual-backward equation of a mean-field PDE. Of course, the toy model~\eqref{eq:toy} is much simpler than the network model~\eqref{eqn:system}, which is non-exchangeable and involves stochastic jumps. Therefore, the non-exchangeable system of interacting neurons model~\eqref{eqn:system} will require a much more general framework, as well as a combination of deterministic and probabilistic methods.

\subsection{Limits of dense synaptic weight matrices}\label{sec:graphon}
The auxiliary variable $\xi \in [0,1]$ we use in the definition of the extended empirical measure~\eqref{eqn:extended_empirical_measure} is motivated by the theory of graphons \cite{Lov12}. This theory provides a topological framework for studying the limits of large dense graphs and comes with a standard representation for graphs and their large size limits.

The theory of graphons is primarily concerned with the study of the large size limits of simple graphs that are undirected and unweighted.
A measurable function $w : I \times I \to \R$, where $I = [0,1]$, is a \textit{graphon} if it is symmetric, non-negative and bounded by $1$.
Let $I^{1;N},\dots,I^{N;N}$ be intervals of equal length that form a partition of $I = [0,1]$. A simple graph $G$ with $N$ nodes can be associated with a graphon $W_G$ defined as
\begin{equation*}
    W_G(\xi,\zeta) := \begin{cases} 1 \quad\text{if } (\xi,\zeta) \in I^{i;N} \times I^{j;N}$ \text{ and } $(i,j)\text{ is an edge of }G, \\
    0 \quad\text{otherwise}.  
\end{cases}
\end{equation*}
This definition can be extended to graphs that are directed and weighted. Dropping the requirement for symmetry and non-negativity, a measurable function $w \in L^\infty([0,1]^2)$ is referred to as \textit{digraphon} in the graphon theory literature.
A weight matrix $\{w^N_{i,j}\}_{i,j = 1}^N \in \R^{N \times N}$ is associated with the digraphon $w^{I;N} \in L^\infty([0,1]^2)$, defined as
\begin{equation*}
\begin{aligned}
w^{I;N}(\xi,\zeta) \defeq \sum_{i,j=1}^N w^N_{i,j}\mathbbm{1}_{I^{i;N}}(\xi)\mathbbm{1}_{I^{j;N}}(\zeta).
\end{aligned}
\end{equation*}

Graphons and digraphons can be considered as $L^p([0,1]) \to L^q([0,1])$ operators for any $1 \leq p,q \leq \infty$. One of the essential findings of graphon theory is that the most convenient norm for studying dense graph limits is, arguably, the operator norm $\|\cdot\|_{L^\infty([0,1]) \to L^1([0,1])}$, which is weaker that than the $L^p$ norms $\|\cdot\|_{L^p([0,1]^2)}$. It is important to mention that the operator norm $\|\cdot\|_{L^\infty([0,1]) \to L^1([0,1])}$ is equivalent to the \textit{cut norm} $\|\cdot\|_{\square}$ more commonly used in graphon theory, as explained in \cite[Chapter~8]{Lov12}.

When considering graphs with unlabeled nodes, the order of the nodes is arbitrary. Similarly, considering our non-exchangeable system~\eqref{eqn:system}, for any permutation $\sigma$ of the neuron indices $\{1,\dots,N\}$, $\{\bm{X}^{\sigma(i);N}\}_{i = 1}^N$ is a solution of \eqref{eqn:system} with weight matrix $\{w^N_{\sigma(i),\sigma(j)}\}_{i,j = 1}^N$ if and only if $\{\bm{X}^{i;N}\}_{i = 1}^N$ is a solution of \eqref{eqn:system} with weight matrix $\{w^N_{i,j}\}_{i,j = 1}^N$. This motivates the definition of a distance between graphs (or systems) up to node (or neuron) re-orderings. Such a distance is given by the \textit{cut distance} $\delta_\square$, which is the infimum of the operator norm $\|\cdot\|_{L^\infty([0,1]) \to L^1([0,1])}$ over all possible re-orderings: 
\begin{equation*}
\begin{aligned}
\delta_{\square}(W_1, W_2) \defeq \inf_{\psi} \|W_1^{\psi} - W_2\|_{L^\infty \to L^1}, \quad \forall (W_1,W_2)\in L^\infty([0,1]^2)^2,
\end{aligned}
\end{equation*}
where $\psi$ ranges over all measure-preserving maps from $I$ to $I$, and, for any measure-preserving map $\psi$, $W^{\psi}(\xi,\zeta) = W(\psi(\xi),\psi(\zeta))$. Actually, the cut distance $\delta_{\square}$ is only a pseudo-distance since any $W \in L^\infty([0,1]^2)$ has distance $0$ with its ``re-ordered'' version $W^\psi$.





The cut distance $\delta_{\square}$ is remarkable because is leads to a compactness result on the space of equivalent classes it induces.
\begin{lemma} \label{lem:graphon_compactness}
For any sequence $W^N$, $N\to\infty$, of elements in $L^\infty([0,1]^2)$ satisfying the uniform boundedness condition $\sup_N \|W^N\|_{L^\infty} < \infty$, we can extract a subsquence $\{W^{n(N)}\}_{N=1}^\infty$ such that there exists a digraphon $W \in L^\infty([0,1]^2)$ such that $\delta_{\square}(W^{n(N)},W) \to 0$ as $N\to\infty$.
\end{lemma}

The proof of the lemma, in the case of standard graphons, can be found, for example, in \cite[Chapter~9]{Lov12}; the proof in the more general case of digraphons and their extensions can be found in \cite{LovSze10, LovVes13}.

For a sequence of digraphons corresponding to weight matrices, the following lemma is a straightforward consequence of Lemma~\ref{lem:graphon_compactness}.




\begin{lemma} \label{lem:graphon_compactness_2}
For any sequence $\{w^N_{i,j}\}_{1\leq i,j\leq N}$, $N\to\infty$, of weight matrices with increasing sizes $N\times N$ satisfying the uniform boundedness condition
\begin{equation*}
\begin{aligned}
\sup_N \max_{1 \leq i,j \leq N} |w^N_{i,j}| < \infty,
\end{aligned}
\end{equation*}
we can extract a subsequence $\{w^{n(N)}_{i,j}\}_{1\leq i,j\leq n(N)}$ such that there exist a digraphon $w \in L^\infty([0,1]^2)$ and almost everywhere partitions $\{\{E^{i;n(N)}\}_{i = 1}^{n(N)}\}_{N=1}^\infty$ such that the digraphons
\begin{equation*}
\begin{aligned}
w^{n(N)}(\xi,\zeta) \defeq w^{E;n(N)}(\xi,\zeta) = \sum_{i,j=1}^{n(N)} w^{n(N)}_{i,j}\mathbbm{1}_{E^{i;n(N)}}(\xi)\mathbbm{1}_{E^{j;n(N)}}(\zeta)
\end{aligned}
\end{equation*}
satisfy $\| w^{n(N)} - w\|_{L^\infty \to L^1} \to 0$ as $N\to\infty$.
\end{lemma}

\subsection{Defining a kernel-dependent weak metric for the empirical measure}\label{sec:metric}

Passing to the limit in the weights $w^N_{i,j}$ to obtain a limit kernel $w$ as $N\to\infty$ is only the first step towards passing to the limit in the empirical measure, $\bm{\mu}^N \wsto \mu$.
Our proof strategy relies on identifying a proper metric for measures allowing us to pass to the limit from the non-exchangeable system~\eqref{eqn:system} with weights $w^N$ to the limit mean-field PDE~\eqref{eqn:Vlasov} with the kernel $w$. 
To this end, the metric should correspond to the weak-* convergence discussed in Section~\ref{sec:exchangeable} and the form of the limit equation~\eqref{eqn:Vlasov}.
Note that the naive choice of the $W^{-1,1}([0,1] \times \R)$-norm does not appear to work as
proving the propagation of this norm seems to require the kernel $w$ to be Lipschitz, while we only assume $w$ to be in $L^\infty([0,1]^2)$.

\begin{subequations}

We introduce a metric that is \textit{adapted} to the limit kernel $w$. For any $\mu_1, \mu_2 \in L^\infty([0,1] ; \mathcal{M} (\R))$, we define the $\Phi_w^{-1}$ metric between them through the duality formula
\begin{equation} \label{eqn:weak_distance}
\|\mu_1 - \mu_2\|_{\Phi_w^{-1}} \defeq \sup_{ \| \phi \|_{\Phi_w} \leq 1} \int_{[0,1] \times \R} \phi(\xi,x) \Big( \mu_1(\xi,x) - \mu_2(\xi,x) \Big) \; \rd \xi \rd x,
\end{equation}
\begin{equation} \label{eqn:dual_norm}
\begin{aligned}
\| \phi \|_{\Phi_w} &\defeq \max \bigg( \|\phi\|_{L^\infty_{\xi,x}} , \|\partial_x \phi\|_{L^\infty_{\xi}L^1_x}, \|\partial_x \phi\|_{L^\infty_{\xi,x}},
\\
& \hphantom{\defeq }
\inf \bigg\{ C > 0 :  \forall h \in (-1,1) , \; \int_{[0,1]} \sup_{x \in \R} \Big| \phi(\xi - h,x) - \phi(\xi,x) \Big| \;\rd \xi \leq C \epsilon_w(|h|) \bigg\} \bigg),
\end{aligned}
\end{equation}
where $\epsilon_w : [0,1] \to [0,\infty)$ is a non-decreasing function providing a modulus of continuity of $w$ in $L^1$:
\begin{equation} \label{eqn:adapted_regularity}
\begin{aligned}
& \lim_{r \to 0^+} \epsilon_w(r) = \epsilon_w(0) = 0,
\\
& \int_{[0,1]^2} \big| w(\xi - h, \zeta) - w(\xi, \zeta) \big| \;\rd \zeta \rd \xi \leq \epsilon_w(|h|), \quad \forall h \in \R,
\\
& \int_{[0,1]^2} \big| w(\xi, \zeta - h) - w(\xi, \zeta) \big| \;\rd \zeta \rd \xi \leq \epsilon_w(|h|), \quad \forall h \in \R.
\end{aligned}
\end{equation}
\end{subequations}
In \eqref{eqn:weak_distance}, we make a small abuse of notation, writing $(\mu_1(\xi,x)-\mu_2(\xi,x))\,\rd \xi\,\rd x$ instead of $(\mu_1(\rd \xi,\rd x)-\mu_2(\rd \xi,\rd x))$. This abuse of notation will be useful to indicate the order of integration. 

In \eqref{eqn:adapted_regularity}, $\xi - h$ may exit the domain $[0,1]$. To address this, we adopt the following convention: Any $w \in L^\infty([0,1]^2)$ is identified with its natural periodic extension to $L^\infty(\T^2)$, where $\T^2$ denotes the torus in dimension~$2$. Similarly, any $\mu \in L^\infty([0,t_*] \times [0,1]; \mathcal{M}(\R))$ is extended to $L^\infty([0,t_*] \times \T; \mathcal{M}(\R))$, where $\T$ denotes the torus in dimension~$1$ (the circle). This convention is to be applied to all subsequent functions involving the variable $\xi \in [0,1]$.

The theory of $L^p$ spaces tells us that for any $w \in L^\infty([0,1]^2)$, there always exists a modulus of continuity $\epsilon_w$ satisfying~\eqref{eqn:adapted_regularity}.\footnote{Indeed, $w\in L^\infty([0,1]^2)$ implies $w\in L^1([0,1]^2)$, and the space continuous functions $C([0,1]^2)$ is dense in $L^1([0,1]^2)$.}
This can be interpreted as a weak form of regularity, adapted to the particular choice of $w$, which is nevertheless stronger than mere $L^\infty$ regularity.
Exploiting this slightly stronger regularity, we will be able to obtain a stability analysis with respect to the $\Phi_w^{-1}$ metric between the system~\eqref{eqn:system} and the limit equation~\eqref{eqn:Vlasov}.


\begin{theorem} \label{thm:main_metric}
Grant Assumptions~\ref{assumption:b_f}. Let $\{\bm{X}^{i;N}\}_{i = 1}^N$, $N \to \infty$, be a sequence of solutions to~\eqref{eqn:system} with weight matrices $\{w^N_{i,j}\}_{i,j = 1}^N$ satisfying the uniform boundedness condition~\eqref{eq:boundedness_condition} and the uniform a priori moment bound
\begin{equation}\label{eq:th2_moment_bound}
\begin{aligned}
\sup_{N} \frac{1}{N} \sum_{i = 1}^N \E \big[ \eta \big(\bm{X}^{i;N}_t \big) \big] =
\sup_{N} \E\Big[ \int_{[0,1]} \int_{\R} \eta(x) \bm{\mu}^N(t,\xi,\rd x) \rd \xi \Big]
\leq C(t) < \infty,
\end{aligned}
\end{equation}
for some moment function $\eta:\R\to\R_+$ satisfying $\|\eta^{-1}\|_{L^1} < \infty$.
Assume that $w$ is a digraphon and $w^N$ are digraphons corresponding to $\{w^N_{i,j}\}_{i,j = 1}^N$ and some almost everywhere partition $\{E^{i;N}\}_{i = 1}^N$ such that $\|w^N - w\|_{L^\infty \to L^1} \to 0$.
Moreover, assume that $\mu \in L^\infty([0,t_*] \times [0,1]; \mathcal{M}(\R))$ is a solution to the mean-field PDE~\eqref{eqn:Vlasov} with kernel $w$, and $\bm{\mu}^N$ are extended empirical measures corresponding to $\{\bm{X}^{i;N}\}_{i = 1}^N$ and $\{E^{i;N}\}_{i = 1}^N$ such that
\begin{equation}\label{eq:th2_mu0}
\begin{aligned}
\E\Big[ \| \bm{\mu}^N(0,\cdot,\cdot) - \mu(0,\cdot,\cdot) \|_{\Phi_w^{-1}} \Big] \to 0,  \quad \text{ as } N \to 0.
\end{aligned}
\end{equation}
Then, for all $t_* > 0$,
\begin{equation}\label{eq:th2_Phi-convergence}
\begin{aligned}
\sup_{t \in [0,t_*]} \E\Big[ \| \bm{\mu}^N(t,\cdot,\cdot) - \mu(t,\cdot,\cdot) \|_{\Phi_w^{-1}} \Big] \to 0, \quad \text{ as } N \to 0.
\end{aligned}
\end{equation}
\end{theorem}

Theorem~\ref{thm:main_metric} can be seen as the core of Theorem~\ref{thm:main}. Below, we show that the assumptions of Theorem~\ref{thm:main} imply the assumptions of Theorem~\ref{thm:main_metric} and, conversely, the result of Theorem~\ref{thm:main_metric} implies the result of Theorem~\ref{thm:main}. The relationship between Theorems~\ref{thm:main}~and~\ref{thm:main_metric} can be decomposed in three propositions. Since their proofs follow standard arguments, they are postponed to Appendix~\ref{appen:supplementary}.


\begin{proof}[Proof of Theorem~\ref{thm:main} via Theorem~\ref{thm:main_metric}]

To start with, by Lemma~\ref{lem:graphon_compactness_2}, we can extract a subsequence and choose almost everywhere partitions $\{E^{i;N}\}_{i = 1}^N$ such that $\|w^N - w\|_{L^\infty \to L^1} \to 0$. The first proposition states that the initial moment bound~\eqref{eq:th1_moment_bound} in Theorem~\ref{thm:main} can be propagated through time to obtain the uniform moment bound~\eqref{eq:th2_moment_bound} in Theorem~\ref{thm:main_metric} (choosing $\eta(x) = 1 + x^2$).
\begin{proposition} \label{prop:moments_propagate}

Let $\{\bm{X}^{i;N}\}_{i = 1}^N$ be a solution to \eqref{eqn:system} with weight matrix $(w^N_{i,j})_{i,j=1}^N$ and with initial moment bound
\begin{equation*}
\begin{aligned}
\;& \frac{1}{N} \sum_{i = 1}^N \E \big[ \big|\bm{X}^{i;N}_0 \big|^2 \big] = E_0
< \infty.
\end{aligned}
\end{equation*}
Then, for all $t \in [0,\infty)$,
\begin{equation*}
\begin{aligned}
\;& \frac{1}{N} \sum_{i = 1}^N \E \big[ \big|\bm{X}^{i;N}_t \big|^2 \big] = C(t) < \infty,
\end{aligned}
\end{equation*}
where $C(t)$ may depend on $E_0$, $\|b\|_{L^\infty}$, $\|f\|_{L^\infty}$, and $\max_{1 \leq i,j \leq N}|w^N_{i,j}|$.
\end{proposition}

The second proposition establishes the equivalence between the expected $\Phi_w^{-1}$-convergence in Theorem~\ref{thm:main_metric} and the weak-* convergence in Theorem~\ref{thm:main}.
\begin{proposition} \label{prop:equiv_topology_strong}
Let $\mu^N$, $N \to \infty$, be a sequence of measure on $[0,1] \times \R$, and $\mu$ be a measure on $[0,1] \times \R$. Assume that $\mu^N(\xi,\cdot) \in \mathcal{P}(\R)$ for a.e. $\xi \in [0,1]$ and $\mu(\xi,\cdot) \in \mathcal{P}(\R)$ for a.e. $\xi \in [0,1]$. Further, let $\epsilon_w : (0,1) \to (0,\infty)$ be a non-decreasing function satisfying $\lim_{r \to 0+} \epsilon_w(r) = 0$.
The following statements are equivalent:
\begin{itemize}
\item[(i)] 
$\|\mu^N - \mu\|_{\Phi_w^{-1}} \to 0$ as $N \to \infty$,
where the metric $\Phi_w^{-1}$ is defined in~\eqref{eqn:weak_distance}-\eqref{eqn:dual_norm}.
\item[(ii)] 
For all $\varphi \in C_c([0,1] \times \R)$,
\begin{equation*}
\begin{aligned}
\int_{[0,1] \times \R} \varphi ( \mu^N - \mu ) \;\rd \xi \rd x \to 0, \ \text{ as } \ N \to \infty.
\end{aligned}
\end{equation*}

\end{itemize}

\end{proposition}
Statement (ii) in Proposition~\ref{prop:equiv_topology_strong} defines the weak-* convergence in $\mathcal{M}([0,1]\times \R)$ we use in Theorem~\ref{thm:main}.
The third proposition states that the assumptions on the initial data in Theorem~\ref{thm:main} imply the expected $\Phi^{-1}_w$-convergence of the empirical measure at time $0$, Eq.~\eqref{eq:th2_mu0}, in Theorem~\ref{thm:main_metric}.

\begin{proposition} \label{prop:passing_limit_initial}

Assume that the sequence of initial data $\{\bm{X}^{i;N}_0\}_{i = 1}^N$, $N \to \infty$, are made of independent random variables with $\textnormal{Law}(\bm{X}^{i;N}_0) = \mu^{i;N}_0 \in \mathcal{P}(\R)$ and moments are bounded by
\begin{equation*}
\begin{aligned}
\sup_{N \in \N} \sup_{1 \leq i \leq N} \E\big[|\bm{X}^{i;N}_0|^2\big] < \infty.
\end{aligned}
\end{equation*}
Moreover, assume that $\{E^{i;N}\}_{i = 1}^N$, $N \to \infty$, are almost everywhere partitions.

Then, up to extraction of a subsequence, there exists $\mu_0 \in L^\infty([0,1]; \mathcal{M}(\R))$ such that $\mu_0(\xi,\cdot) \in \mathcal{P}(\R)$ for a.e. $\xi \in [0,1]$ and
\begin{equation*}
\begin{aligned}
\bigg\| \bigg( \sum_{i = 1}^N \mu^{i;N}_0 \mathbbm{1}_{E^{i;N}} \bigg) - \mu_ 0 \bigg\|_{\Phi_w^{-1}} \to 0, \ \text{ as } \ N \to \infty.
\end{aligned}
\end{equation*}
Further, for this subsequence, the extended empirical measure $\bm{\mu}^N_0 = \sum_{i = 1}^N \delta_{\bm{X}^{i;N}_0} \mathbbm{1}_{E^{i;N}}$ satisfies
\begin{equation*}
\begin{aligned}
\E \big[ \big\| \bm{\mu}^N_0 - \mu_0 \big\|_{\Phi_w^{-1}} \big] \to 0,  \ \text{ as } \ N \to \infty.
\end{aligned}
\end{equation*}

\end{proposition}

In summary, Propositions~\ref{prop:moments_propagate} and \ref{prop:passing_limit_initial} guarantee that all the assumptions of Theorem~\ref{thm:main_metric} are satisfied under the assumptions of Theorem~\ref{thm:main}. Theorem~\ref{thm:main_metric} establishes the convergence in terms of the expectation of the $\Phi_w^{-1}$ metric. Finally, Proposition~\ref{prop:equiv_topology_strong} shows that the convergence stated in Theorem~\ref{thm:main_metric} implies the weak-* convergence in Theorem~\ref{thm:main}. Hence, to prove Theorem~\ref{thm:main}, it suffices to prove Theorem~\ref{thm:main_metric}.
\end{proof}

\section{Proof of the main result}\label{sec:stability}

\subsection{Well-posedness} \label{subsec:well_posedness}
Before starting the proof of Theorem~\ref{thm:main_metric}, we first describe the solutions to the system of SDEs~\eqref{eqn:system} and to the mean-field PDE~\eqref{eqn:Vlasov}.
The following proposition is about the mean-field PDE~\eqref{eqn:Vlasov}.
In addition to its solution $\mu$, we also introduce the measure $\mu_\#$, which is the pushforward measure of $\mu$ by the measurable function $(t,\xi,x)\mapsto (t,\xi,x-H(t,\xi))$. Formally, the measure $\mu_\#$ can be seen as the result of a change of variable that is rough with respect to the $\xi$ variable.
The motivation for introducing this change of variable will be explained in Section~\ref{sec:more_on}.

\begin{proposition} \label{prop:extended_empirical_measure}

Grant Assumptions~\ref{assumption:b_f}. In addition, assume that $w \in L^\infty([0,1]^2)$ and $\mu_0 \in L^\infty([0,1]; \mathcal{M}(\R))$. Then, for any $t_* > 0$, there exists a unique $\mu \in L^\infty([0,t_*] \times [0,1]; \mathcal{M}(\R))$ such that $\mu$ is a solution to the mean-field PDE~\eqref{eqn:Vlasov}, restated here
\begin{subequations} \notag
\begin{align}\label{eq:Vlasov_prop4}
&\partial_t  \mu(t,\xi,\rd x) + \partial_x \Big[\big(b(x) + h(t,\xi)\big) \mu(t,\xi,\rd x)\Big] + f(x) \mu(t,\xi,\rd x) - r(t,\xi)\delta_0(\rd x) = 0, \tag{\ref{eqn:Vlasov_field}}
\\
& r(t,\xi) = \int_{\R}f(x) \mu(t,\xi,\rd x), \quad h(t,\xi) = \int_{[0,1]} w(\xi,\zeta) r(t,\zeta) \rd \zeta, \tag{\ref{eqn:Vlasov_r_h}}
\\
&\mu(0,\xi,\rd x) = \mu_0(\xi,\rd x), \notag
\end{align}
\end{subequations}
in the sense of characteristics.
Moreover, define
\begin{equation*}
\begin{aligned}
& H(t,\xi) \defeq \int_0^t h(s,\xi) \;\rd s = \int_0^t \int_{[0,1]} w(\xi,\zeta) \int_{\R}f(x) \mu(s,\zeta,\rd x)  \rd \zeta \rd s,
\\
& \mu_\#(t,\xi,\rd x) \defeq \mu(t,\xi,\rd x+H(t,\xi)),
\end{aligned}
\end{equation*}
together with
\begin{equation*}
\left\{
\begin{aligned}
&b_\#(t,\xi,x) \defeq b(x+H(t,\xi)),
\\
&f_\#(t,\xi,x) \defeq f(x+H(t,\xi)),
\\
&\delta_\#(t,\xi,\rd x) \defeq \delta_0(\rd x+H(t,\xi)).
\end{aligned} \right.
\end{equation*}
Then, $H \in L^\infty([0,1]; C([0,t_*]))$ and the equation
\begin{multline}\label{eq:Vlasov_3}
\partial_t \mu_\#(t,\xi,\rd x) + \partial_x\Big[b_\#(t,\xi,x) \mu_\#(t,\xi,\rd x)\Big] \\+ f_\#(t,\xi,x) \mu_\#(t,\xi,\rd x) - r(t,\xi) \delta_\#(t,\xi,\rd x) = 0 
\end{multline}
holds in the following distributional sense:  $\forall \varphi \in L^\infty([0,1]; C_b^1([0,t_*] \times \R)), t \in [0,t_*]$,
\begin{equation} \label{eqn:Vlasov_4}
\begin{aligned}
& \int_{[0,1] \times \R} \varphi(t,\xi,x) \mu_\#(t,\xi,\rd x) \;\rd \xi 
-
\int_{[0,1] \times \R} \varphi(0,\xi,x) \mu_\#(0,\xi,\rd x) \;\rd \xi 
\\
& \quad = \int_0^t \int_{[0,1] \times \R} \left[\partial_s \varphi(s,\xi,x) + b_\#(s,\xi,x) \partial_x \varphi(s,\xi,x)\right] \mu_\#(s,\xi,\rd x) \;\rd \xi \rd s
\\
& \quad \quad + \int_0^t \int_{[0,1] \times \R} \varphi(s,\xi,x) \big[ - f_\#(s,\xi,x) \mu_\#(s,\xi,\rd x) + r(s,\xi) \delta_\#(s,\xi,\rd x) \big]
\;\rd \xi \rd x \rd s.
\end{aligned}
\end{equation}
\end{proposition}

For clarity, let us explicitly state what we mean, above, by ``solution in the sense of characteristics.''
Given $\mu \in L^\infty([0,t_*] \times [0,1]; \mathcal{M}(\R))$, consider the unique flow solving the system of characteristic equations
\begin{equation*}
\left\{
\begin{aligned} 
\frac{\rd X_\mu}{\rd t}(t,s,\xi,x) & = b\big(X_\mu (t,s,\xi,x)\big) + h(t,\xi),
\\
X_\mu(s,s,\xi,x) & = x,
\end{aligned} \right.
\end{equation*}
for a.e. $\xi \in [0,1]$. 
(Note that $h(t,\xi)$, defined in \eqref{eqn:Vlasov_field}, depends on $\mu$.)
The measure $\mu$ is a solution to the mean-field PDE~\eqref{eqn:Vlasov} in the sense of characteristics, if for any $t \in [0,t_*]$ and a.e. $\xi \in [0,1]$, 
\begin{equation*}
\begin{aligned} 
\mu(t,\xi,\rd x) = \mu_0\big(\xi, X_\mu(0,t,\xi,\rd x) \big) \partial_x X_\mu(0,t,\xi,x).
\end{aligned}
\end{equation*}

\begin{proof}
    The well-posedness of the nonlinear mean-fiel PDE~\eqref{eqn:Vlasov} can be obtained by adapting a standard contraction argument for population equations without spatial extension \cite{Che17,CorTan20,Sch22,FonSch22} to the spatially-extended case; the spatial extension $\xi$ in \eqref{eqn:Vlasov} does not cause any difficulty under the assumption that $w\in L^\infty([0,1]^2)$. 
    
    Eqs.~\eqref{eq:Vlasov_3} and \eqref{eqn:Vlasov_4} simply stem from the change of variable $\mu_\#(t,\xi,\rd x) \defeq \mu(t,\xi,\rd x+H(t,\xi))$.
\end{proof}

The next proposition establishes the well-posedness of the system of SDEs \eqref{eqn:system}, together with a change of variable similar to the one we used for $\mu_\#$ in the previous proposition.

\begin{proposition} \label{prop:SDE_strong_solution}
Grant Assumptions~\ref{assumption:b_f}. For all $N \geq 1$, let $\{\bm{X}^{i;N}_0\}_{i = 1}^N$ be $\mathscr{F}_0$-measurable random variables. 
\begin{itemize}
\item[(i)] There exists a pathwise unique càdlàg strong solution $\{\bm{X}^{i;N}\}_{i = 1}^N$ to the system of SDEs~\eqref{eqn:system}, in the sense that $\{\bm{X}^{i;N}\}_{i = 1}^N$ is a càdlàg $(\mathscr{F}_t)_{t \in \R_+}$-adapted process solving, a.s.,
\begin{multline*} 
\bm{X}^{i;N}(t) = 
\bm{X}^{i;N}_0 + \int_0^t b(\bm{X}^{i;N}(s)) \rd s \\
+ \frac{1}{N} \sum_{j=1}^N \int_{[0,t]\times \R_+} w^N_{i,j} \mathbbm{1}_{\{z\leq f(\bm{X}^{j;N}(s-))\}}\bm{\Pi}^{j;N}(\rd s, \rd z)
\\
- \int_{[0,t]\times \R_+} \bm{X}^{i;N}(s-) \mathbbm{1}_{\{z\leq f(\bm{X}^{i;N}(s-))\}}\bm{\Pi}^{i;N}(\rd s, \rd z), \quad \forall t \geq 0.
\end{multline*}

\item[(ii)] Moreover, define the random processes $\{\bm{H}^{i;N}\}_{i = 1}^N$, $\{\bm{Y}^{i;N}\}_{i = 1}^N$ as
\begin{equation*}
\begin{aligned}
\bm{Y}^{i;N}(t) &:= \bm{X}^{i;N}(t) - \bm{H}^{i;N}(t),
\\
\bm{H}^{i;N}(t) &:= \frac{1}{N} \sum_{j=1}^N \int_{[0,t]\times \R_+} w^N_{i,j} \mathbbm{1}_{\{z\leq f(\bm{X}^{j;N}(s-))\}} \bm{\Pi}^{j;N}(\rd s, \rd z).
\end{aligned}
\end{equation*}
Then, a.s., for all $\varphi^{i;N}\in C^1_b([0,t_*] \times \R)$, $i \in \{1, \dots, N\}$,

\begin{equation} \label{eqn:SPDE_2_Ito}
\begin{aligned}
& \sum_{i = 1}^N \varphi^{i;N}(t,\bm{Y}^{i;N}(t)) - \sum_{i = 1}^N \varphi^{i;N}(0,\bm{Y}^{i;N}(0))
\\
& \quad = \sum_{i = 1}^N \int_0^t \partial_s \varphi^{i;N}(s,\bm{Y}^{i;N}(s)) + b(\bm{Y}^{i;N}(s) + \bm{H}^{i;N}(s)) \partial_x \varphi^{i;N}(s,\bm{Y}^{i;N}(s)) \;\rd s
\\
& \quad \quad + \sum_{i = 1}^N \int_{[0,1] \times \R} \big[ \varphi^{i;N}(s,-\bm{H}^{i;N}(s-)) - \varphi^{i;N}(s,\bm{Y}^{i;N}(s-)) \big]
\\
& \hspace{7cm} \mathbbm{1}_{\{z\leq f(\bm{Y}^{i;N}(s-) + \bm{H}^{i;N}(s-))\}}\bm{\Pi}^{i;N}(\rd s, \rd z).
\end{aligned}
\end{equation}
\end{itemize}
\end{proposition}
\begin{proof}
    The well-posedness (i) has been established in \cite{FouLoe16} and, as discussed in \cite{GalLoe20,Sch22}, a.s., a solution can be constructed via a standard thinning procedure \cite{LewShe79,Oga81,DelFou16}. For (ii), since the solution $\{\bm{X}^{i;N}\}_{i = 1}^N$ has bounded variation a.s., $\{\bm{Y}^{i;N}\}_{i = 1}^N$ also has bounded variation a.s. Hence, for any realization of $\{\bm{Y}^{i;N}\}_{i = 1}^N$, a.s., \eqref{eqn:SPDE_2_Ito} holds for all $\varphi^{i;N}\in C_b^1([0,t_*] \times \R)$ by simple differentiation of $\varphi^{i;N}(s,\bm{Y}^{i;N}(s))$ in the distributional sense.
\end{proof}

We would like to emphasize that Proposition~\ref{prop:SDE_strong_solution} states the identity \eqref{eqn:SPDE_2_Ito} in a stronger sense than It\^o's lemma usually does.
The most widely known version of It\^o's lemma deals with semi-martingales and only asserts that for any test function $\varphi$, the It\^o formula (which is \eqref{eqn:SPDE_2_Ito} in our case) is satisfied a.s. 
In contrast, Proposition~\ref{prop:SDE_strong_solution} 
states that for a.e. realization $\omega \in \Omega$, \eqref{eqn:SPDE_2_Ito} holds for \emph{all} $\varphi^{i;N} \in C^1_b([0,t_*] \times \R)$.
This is possible because $(\bm{Y}^{i;N})_{i = 1}^N$ and $(\bm{H}^{i;N})_{i = 1}^N$ are bounded variation processes, and the integral in \eqref{eqn:SPDE_2_Ito} can be understood in the Stieltjes sense.
To see that these two types of statements are not equivalent, consider \eqref{eqn:SPDE_2_Ito} where $\varphi^{i;N}$ is replaced by any random $\bm{\varphi}^{i;N}$ such that $\bm{\varphi}^{i;N} \in C_b^1([0,t_*] \times \R)$, $i \in \{1,\dots, N\}$ a.s. Proposition~\ref{prop:SDE_strong_solution} guarantees that the identity remains true a.s., while a standard It\^o lemma-type statement would be insufficient to guarantee the identity in this case.


The identity \eqref{eqn:SPDE_2_Ito} can be reformulated as an expression similar to \eqref{eqn:Vlasov_4} by the notion of extended empirical measure. Given any almost everywhere partition $\{E^{i;N}\}_{i = 1}^N$, we define
\begin{equation*}
\begin{aligned}
\bm{H}^N(t,\xi) \defeq \;& \textstyle \sum_{i=1}^N\bm{H}^{i;N}(t)\mathbbm{1}_{E^{i;N}}(\xi),
\\
\bm{\mu}_\#^N(t,\xi,\rd x) \defeq \;& \textstyle \sum_{i=1}^N\delta_{\bm{Y}^{i;N}(t)}(\rd x)\mathbbm{1}_{E^{i;N}}(\xi),
\end{aligned}
\end{equation*}
in accordance with the previous definition of $\bm{\mu}^N$, restated here
\begin{equation*}
\begin{aligned}
\bm{\mu}^N(t,\xi,\rd x) \defeq \;& \textstyle \sum_{i=1}^N\delta_{\bm{X}^{i;N}(t)}(\rd x)\mathbbm{1}_{E^{i;N}}(\xi).
\end{aligned}
\end{equation*}
Then, a.s., for all $\varphi \in L^\infty([0,1]; C_b^1([0,t_*] \times \R), t \in [0,t_*]$,
\begin{equation} \label{eqn:SPDE_2}
\begin{aligned}
& \int_{[0,1] \times \R} \varphi(t,\xi,x) \bm{\mu}_\#^N(t,\xi,\rd x) \;\rd \xi 
-
\int_{[0,1] \times \R} \varphi(0,\xi,x) \bm{\mu}_\#^N(0,\xi,\rd x) \;\rd \xi 
\\
& \quad = \int_0^t \int_{[0,1] \times \R} \left[\partial_s \varphi(s,\xi,x) + \bm{b}_\#^N(s,\xi,x) \partial_x \varphi(s,\xi,x)\right] \bm{\mu}_\#^N(s,\xi,\rd x) \;\rd \xi \rd s
\\
& \quad \quad + \int_{[0,1] \times \R} \int_{[0,t]\times\R_+} \varphi(s,\xi,x) \big[ - \bm{\mu}_\#^N(s-,\xi,\rd x) + \bm{\delta}_\#^N(s-,\xi,\rd x) \big]
\\
& \quad\quad\quad\quad\quad\quad\quad\quad
\sum_{i=1}^N\mathbbm{1}_{E^{i;N}}(\xi) \mathbbm{1}_{\{z\leq \bm{r}^N(s-,\xi)\}} \bm{\Pi}^{i;N}(\rd s, \rd z)
\;\rd \xi,
\end{aligned}
\end{equation}
where the random fields and measures are defined as
\begin{equation*}
\left\{
\begin{aligned}
& \bm{b}_\#^N(t,\xi,x) \defeq b(x+\bm{H}^N(t,\xi)),
\\
& \bm{f}_\#^N(t,\xi,x) \defeq f(x+\bm{H}^N(t,\xi)),
\\
& \bm{\delta}_\#^N(t,\xi,\rd x) \defeq \delta_0(\rd x+\bm{H}^N(t,\xi)),
\end{aligned} \right.
\end{equation*}
and
\begin{equation*}
\begin{aligned}
\bm{r}^N(s,\xi) \defeq \sum_{i=1}^N f(\bm{X}^{i;N}(s)) \mathbbm{1}_{E^{i;N}}(\xi) = \int_{\R} f(x) \bm{\mu}^N(s,\xi,\rd x) = \int_{\R} \bm{f}_\#^N(s,\xi,x) \bm{\mu}_\#^N(s,\xi,\rd x).
\end{aligned}
\end{equation*}
To see that \eqref{eqn:SPDE_2} is a reformulation of \eqref{eqn:SPDE_2_Ito}, it suffices to notice that $\varphi^{i;N}(t,x) = \int_{E^{i;N}} \varphi(t,\xi,x)\rd \xi$.

\subsection{On the change of variable involving the intergrated postsynaptic input}\label{sec:more_on}
In this subsection, we explain the usefulness of the aforementioned change of variable. In the original system \eqref{eqn:system}, the SDEs governing $\bm{X}^{i;N}$ with distinct indices $i \in \{1,\dots, N\}$ are coupled via the term
\begin{equation*}
\begin{aligned}
\frac{1}{N}\sum_{j=1}^N w^N_{i,j}\rd \bm{Z}^{j;N}(t).
\end{aligned}
\end{equation*}
This coupling term is interpreted as the postsynaptic input to the $i$-th neuron and is singular in time.
Intuitively, as $N \to \infty$, we expect this term to converge to the drift term in the mean-field PDE~\eqref{eqn:Vlasov}
\begin{equation*}
\begin{aligned}
h(t,\xi) = \int_{[0,1]} w(\xi,\zeta) r(t,\zeta) \rd \zeta,
\end{aligned}
\end{equation*}
which is continuous in time. However, the singularity in time poses serious challenges if we try to directly pass to the limit from the system~\eqref{eqn:system} to the mean-field PDE~\eqref{eqn:Vlasov}. To address this difficulty, we separate the coupling term from the rest of the dynamics.
In Proposition~\ref{prop:extended_empirical_measure}, we have defined the ``integrated drift term'' $H(t,\xi)$ as
\begin{equation*}
\begin{aligned}
H(t,\xi) \defeq \int_0^t h(s,\xi) \;\rd s = \int_0^t \int_{[0,1]} w(\xi,\zeta) r(s,\zeta) \rd \zeta \rd s,
\end{aligned}
\end{equation*}
and adopted the change of variable $\mu_\#(t,\xi,\rd x) \defeq \mu(t,\xi,\rd x + H(t,\xi))$.
Analogously, in Proposition~\ref{prop:SDE_strong_solution} we have defined the ``integrated postsynaptic input'' $\bm{H}^{i;N}(t)$ as
\begin{equation*}
\begin{aligned}
\bm{H}^{i;N}(t) \defeq \;& \frac{1}{N} \sum_{j=1}^N \int_0^t w^N_{i,j} \rd \bm{Z}^{j;N}(s) = \frac{1}{N} \sum_{j=1}^N \int_{[0,t]\times \R_+} w^N_{i,j} \mathbbm{1}_{\{z\leq f(\bm{X}^{j;N}(s-))\}}\bm{\Pi}^{j;N}(\rd s, \rd z),
\end{aligned}
\end{equation*}
and applied a subsequent change of variable $\bm{Y}^{i;N}(t) \defeq  \bm{X}^{i;N}(t) - \bm{H}^{i;N}(t)$.
The change of variable prevents the direct appearance of the singular-in-time coupling terms in the dynamics, as we can see in \eqref{eqn:SPDE_2_Ito}, or equivalently \eqref{eqn:SPDE_2}. The cost of this change of variable is that it makes field $\bm{b}_\#^N$ in \eqref{eqn:SPDE_2} time-dependent and random.

The following proposition shows that distances \textit{after} the change of variable control distances \textit{before} the change of variable.
%
\begin{proposition} \label{prop:input_separation}

For any $t_* > 0$, let $\{\bm{X}^{i;N}\}_{i = 1}^N$ and $\{\bm{H}^{i;N}\}_{i = 1}^N$ be $\R^N$-valued random processes on $[0,t_*]$. Given an almost everywhere partition $\{E^{i;N}\}_{i = 1}^N$ and taking the definition $\bm{Y}^{i;N} := \bm{X}^{i;N} - \bm{H}^{i;N}$,  define, for all $t\in[0,t^*]$,
\begin{equation*}
\begin{aligned}
\bm{H}^N(t,\xi) \defeq \;& \textstyle \sum_{i=1}^N\bm{H}^{i;N}(t)\mathbbm{1}_{E^{i;N}}(\xi),
\\
\bm{\mu}^N(t,\xi,\rd x) \defeq \;& \textstyle \sum_{i=1}^N\delta_{\bm{X}^{i;N}(t)}(\rd x)\mathbbm{1}_{E^{i;N}}(\xi),
\\
\bm{\mu}_\#^N(t,\xi,\rd x) \defeq \;& \textstyle \sum_{i=1}^N\delta_{\bm{Y}^{i;N}(t)}(\rd x)\mathbbm{1}_{E^{i;N}}(\xi),
\end{aligned}
\end{equation*}
Further, let $\mu \in L^\infty([0,t_*] \times [0,1]; \mathcal{M}(\R))$ and $H \in L^\infty([0,t_*] \times [0,1])$. Then, for all $t\in[0,t^*]$,
\begin{equation*} 
\begin{aligned}
\;& \bm{\mu}_\#^N(t,\xi,\rd x) = \bm{\mu}^N(t,\xi,\rd x+\bm{H}^N(t,\xi)), \quad \forall t\geq 0,\; \text{a.e. } \xi \in [0,1],
\\
\;& \|\bm{\mu}^N(t,\cdot,\cdot) - \mu(t,\cdot,\cdot) \|_{\Phi_w^{-1}} \leq \|\bm{\mu}_\#^N(t,\cdot,\cdot) - \mu_\#(t,\cdot,\cdot) \|_{\Phi_w^{-1}} + \|\bm{H}^N(t,\cdot) - H(t,\cdot)\|_{L^1}.
\end{aligned}
\end{equation*}
\end{proposition}
\begin{proof}
For a.e. $\xi \in [0,1]$ there exists a unique index $i(\xi) \in \{1,\dots,N\}$ such that $\xi \in E^{i(\xi);N}$ and $\xi \notin E^{j;N}$ for all $j$ distinct from $i(\xi)$. Given this uniqueness,
\begin{equation*} 
\begin{aligned}
\bm{\mu}_\#^N(t,\xi,\rd x) = \;& \delta_{\bm{Y}^{i(\xi);N}(t)}(\rd x) = \delta_{\bm{X}^{i(\xi);N}(t) - \bm{H}^{i(\xi);N}(t)}(\rd x)
\\
= \;& \delta_{\bm{X}^{i(\xi);N}(t)}(\rd x + \bm{H}^{i(\xi);N}(t)) = \bm{\mu}^N(t,\xi,\rd x+\bm{H}^N(t,\xi)),
\end{aligned}
\end{equation*}
which proves the first identity a.e.

Next, using the definition of the metric $\Phi_w^{-1}$ and incorporating the identity just proved yields
\begin{equation*} 
\begin{aligned}
& \|\bm{\mu}^N(t,\cdot,\cdot) - \mu(t,\cdot,\cdot) \|_{\Phi_w^{-1}} = \sup_{ \| \phi \|_{\Phi_w} \leq 1} \int_{[0,1] \times \R} \phi(\xi,x) \Big( \bm{\mu}^N(t,\xi,\rd x) - \mu(t,\xi,\rd x) \Big) \;\rd \xi 
\\
& \quad = \sup_{ \| \phi \|_{\Phi_w} \leq 1} \int_{[0,1] \times \R} \phi(\xi,x) \Big( \bm{\mu}_\#^N(t,\xi,\rd x - \bm{H}^N(t,\xi)) - \mu_\#(t,\xi,\rd x - H(t,\xi)) \Big) \;\rd \xi.
\end{aligned}
\end{equation*}
Note that $\bm{\mu}^N$, $\bm{\mu}_\#^N$, and $\bm{H}^N$ are random. Hence, the function $\phi$ used in the supremum and the integration is also random.

Apply the change of variable, $x \mapsto x + H(t,\xi)$, we get
\begin{equation*} 
\begin{aligned}
& \|\bm{\mu}^N(t,\cdot,\cdot) - \mu(t,\cdot,\cdot) \|_{\Phi_w^{-1}}
\\
& \quad = \sup_{ \| \phi \|_{\Phi_w} \leq 1} \int_{[0,1] \times \R} \phi(\xi,x + H(t,\xi))
\\
& \hspace{3cm} \Big( \bm{\mu}_\#^N(t,\xi,\rd x - \bm{H}^N(t,\xi) + H(t,\xi)) - \mu_\#(t,\xi,\rd x) \Big) \;\rd \xi 
\\
& \quad \leq \sup_{ \| \phi \|_{\Phi_w} \leq 1} \int_{[0,1] \times \R} \phi(\xi,x + H(t,\xi)) \Big( \bm{\mu}_\#^N(t,\xi,\rd x) - \mu_\#(t,\xi,\rd x) \Big) \;\rd \xi 
\\
& \quad \quad + \sup_{ \| \phi \|_{\Phi_w} \leq 1} \int_{[0,1] \times \R} \phi(\xi,x + H(t,\xi))
\\
& \hspace{3cm} \Big( \bm{\mu}_\#^N(t,\xi,\rd x - \bm{H}^N(t,\xi) + H(t,\xi)) - \bm{\mu}_\#^N(t,\xi,\rd x) \Big) \;\rd \xi 
\\
& \quad \eqdef L_\# + L_H.
\end{aligned}
\end{equation*}
By definition, $L_\# = \|\bm{\mu}_\#^N(t,\cdot,\cdot) - \mu_\#(t,\cdot,\cdot) \|_{\Phi_w^{-1}}$. We can further simplify $L_H$ with another change of variable,
\begin{equation*}
\begin{aligned}
L_H
= \;& \sup_{ \| \phi \|_{\Phi_w} \leq 1} \int_{[0,1] \times \R} \Big( \phi(\xi,x+\bm{H}^N(t,\xi)) - \phi(\xi,x + H(t,\xi)) \Big) \bm{\mu}_\#^N(t,\xi,\rd x) \;\rd \xi 
\\
\leq \;& \int_{[0,1] \times \R} | \bm{H}^N(t,\xi) - H(t,\xi)| \bm{\mu}_\#^N(t,\xi,\rd x) \;\rd \xi 
\\
= \;& \|\bm{H}^N(t,\cdot) - H(t,\cdot)\|_{L^1},
\end{aligned}
\end{equation*}
where, in the inequality we used the fact that $\|\phi\|_{\Phi_w} \leq 1$ implies that $\phi$ is $1$-Lipschitz in the $x$-direction, and in the final equality, we integrated with respect to $x$ then $\xi$.
\end{proof}

Proposition~\ref{prop:input_separation} tells us that to control $\E[\|\bm{\mu}^N(t,\cdot,\cdot) - \mu(t,\cdot,\cdot) \|_{\Phi_w^{-1}}]$ (and prove Theorem~\ref{thm:main_metric}), it suffices to control $\E[\|\bm{H}^N(t,\cdot) - H(t,\cdot)\|_{L^1}]$ and $\E[\|\bm{\mu}_\#^N(t,\cdot,\cdot) - \mu_\#(t,\cdot,\cdot) \|_{\Phi_w^{-1}}]$, which are more tractable.
In Section~\ref{subsec:stability_input}, we derive the bound for $\|\bm{H}^N(t,\cdot) - H(t,\cdot)\|_{L^1}$, which illustrates the main ideas of our proof; the estimation of $\|\bm{\mu}_\#^N(t,\cdot,\cdot) - \mu_\#(t,\cdot,\cdot) \|_{\Phi_w^{-1}}$, which is more intricate, is treated in Sections~\ref{subsec:stability_backward}-\ref{subsec:stability_Poisson}. Finally, these estimates are combined in Section~\ref{subsec:stability_conclude} to conclude the proof.

\subsection{Quantifying the difference between the integrated postsynaptic inputs} \label{subsec:stability_input}

This subsection is devoted to the bound of the term $\E[\|\bm{H}^N(t,\cdot) - H(t,\cdot)\|_{L^1}]$. Let us restate here that
\begin{equation*}
\begin{aligned}
H(t,\xi) \defeq \int_0^t\int_{[0,1]} w(\xi,\zeta) \int_{\R} f(x)\mu(t,\zeta,\rd x)  \,\rd\zeta \rd s,
\end{aligned}
\end{equation*}
and
\begin{equation*}
\begin{aligned}
\bm{H}^N(t,\xi) & \defeq \textstyle \sum_{i=1}^N\bm{H}^{i;N}(t)\mathbbm{1}_{E^{i;N}}(\xi),
\\
\bm{H}^{i;N}(t) & \defeq \frac{1}{N} \sum_{j=1}^N \int_{[0,t]\times \R_+} w^N_{i,j} \mathbbm{1}_{\{z\leq f(\bm{X}^{j;N}(s-))\}}\bm{\Pi}^{j;N}(\rd s, \rd z),
\end{aligned}
\end{equation*}
where $\{E^{i;N}\}_{i = 1}^N$ is an almost everywhere partition.

For the forthcoming steps, we also introduce the auxiliary functions
\begin{equation*} 
\begin{aligned}
\bar{\bm{H}}^N(t,\xi) = \;& \int_0^t\int_{[0,1]} w^N(\xi,\zeta) \int_{\R} f(x)\bm{\mu}^N(t,\zeta,\rd x)  \, \rd\zeta \rd s,
\\
\widetilde{\bm{H}}^N(t,\xi) = \;& \int_0^t\int_{[0,1]} w(\xi,\zeta) \int_{\R} f(x)\bm{\mu}^N(t,\zeta,\rd x)  \, \rd\zeta \rd s,
\end{aligned}
\end{equation*}
where
\begin{equation*}
\begin{aligned}
w^N(\xi,\zeta) \defeq \;& \textstyle \sum_{i,j=1}^N w^N_{i,j}\mathbbm{1}_{E^{i;N}}(\xi)\mathbbm{1}_{E^{j;N}}(\zeta),
\\
\bm{\mu}^N(t,\xi,\rd x) \defeq \;& \textstyle \sum_{i=1}^N\delta_{\bm{X}^{i;N}(t)}(\rd x)\mathbbm{1}_{E^{i;N}}(\xi).
\end{aligned}
\end{equation*}
The definition of these auxiliary functions is closely related to that of $H$, the only difference being that $\bm{\mu}^N$ substitutes $\mu$ and that $w^N$ substitutes $w$. 
\begin{lemma} \label{lem:difference_flow_measure}

The following estimate holds:
\begin{equation*} \label{eqn:difference_flow_measure}
\begin{aligned}
\|\bar {\bm H}^N(t,\cdot) - H(t,\cdot)\|_{L^1_\xi} \leq & \kappa_1 \int_0^t \|\bm{\mu}^N(s,\cdot,\cdot) - \mu(s,\cdot,\cdot)\|_{\Phi_w^{-1}} \rd s + \|w^N - w \|_{L^\infty \to L^1} \|f\|_{L^\infty} t,
\\
\kappa_1 \defeq \;& \Big( \max \big( \|f\|_{L^\infty} , \|\partial_x f\|_{L^\infty} , \|\partial_x f\|_{L^1} \big) \|w\|_{L^\infty L^1} + \|f\|_{L^\infty} \Big).
\end{aligned}
\end{equation*}

\end{lemma}

\begin{proof}

We first examine the distance between $\widetilde{\bm{H}}^N$ and $H$, which differ only through the measures $\bm{\mu}^N$ and $\mu$.
\begin{equation*}
\begin{aligned}
& \int_0^1 \bigg|\widetilde{\bm{H}}^N(t,\xi) - H(t,\xi)\bigg| \;\rd\xi
\\
\quad & = \int_0^1 \bigg|\int_0^t \int_{[0,1]} w(\xi, \zeta) \int_{\R} f(x) \big[ \bm{\mu}^N(s,\zeta,\rd x) - \mu(s,\zeta,\rd x) \big] \, \rd \zeta \rd s \bigg| \;\rd\xi
\\
\quad & = \sup_{\psi: \|\psi\|_{L^\infty_\xi} \leq 1} \int_0^1 \int_0^t \int_{[0,1]} \psi(\xi) w(\xi, \zeta) \int_{\R} f(x) \big[ \bm{\mu}^N(s,\zeta,\rd x) - \mu(s,\zeta,\rd x) \big] \, \rd \zeta \rd s \;\rd\xi
\\
\quad &= \sup_{\psi: \|\psi\|_{L^\infty_\xi} \leq 1}  \int_0^t \int_{[0,1]} \int_0^1\psi(\xi) w(\xi, \zeta) \rd\xi \int_{\R} f(x) \big[ \bm{\mu}^N(s,\zeta,\rd x) - \mu(s,\zeta,\rd x) \big] \,\rd\zeta \rd s
\\
\quad &\leq \sup_{\psi: \|\psi\|_{L^\infty_\xi} \leq 1} \int_0^t \bigg\|\int_0^1 \psi(\xi) w (\xi,\cdot) f(\cdot) \rd\xi \bigg\|_{\Phi_w} \|\bm{\mu}^N(s,\cdot,\cdot) - \mu(s,\cdot,\cdot)\|_{\Phi_w^{-1}}  \rd s,
\end{aligned}
\end{equation*}
where in the last inequality, we used the definition of the metric $\Phi^{-1}_w$. Moreover, 
\begin{equation*}
\begin{aligned}
& \sup_{\psi: \|\psi\|_{L^\infty_\xi} \leq 1} \bigg\|\int_0^1 \psi(\zeta) w (\zeta,\cdot) f(\cdot) \rd\zeta \bigg\|_{\Phi_w}
\\
& \quad \leq \sup_{\psi: \|\psi\|_{L^\infty_\xi} \leq 1} \max \bigg( \bigg\|\int_0^1 \psi(\zeta) w (\zeta,\cdot) f(\cdot) \rd\zeta \bigg\|_{L^\infty_{\xi,x}} , \bigg\|\int_0^1 \psi(\zeta) w (\zeta,\cdot) \partial_x f(\cdot) \rd\zeta \bigg\|_{L^\infty_{\xi,x}} ,
\\
& \hspace{1cm} \bigg\|\int_0^1 \psi(\zeta) w (\zeta,\cdot) \partial_x f(\cdot) \rd\zeta \bigg\|_{L^\infty_{\xi}L^1_x},
\\
& \hspace{1cm} 
\sup_{h} \bigg\{ \epsilon_w(|h|)^{-1} \int_{[0,1]} \sup_{x} \bigg| \int_{[0,1]} \psi(\zeta) w (\zeta,\xi) f(x) - \psi(\zeta) w (\zeta,\xi - h) f(x) \;\rd\zeta \bigg| \; \rd \xi \bigg\} \bigg)
\\
& \quad \leq \max\bigg( \max\big( \|f\|_{L^\infty} , \|\partial_x f\|_{L^\infty} , \|\partial_x f\|_{L^1} \big) \sup_{\psi: \|\psi\|_{L^\infty_\xi} \leq 1} \bigg\|\int_0^1 \psi(\zeta) w (\zeta,\cdot) \rd\zeta \bigg\|_{L^\infty_{\xi}},
\\
& \quad \hphantom{\leq \max\bigg(}
\|f\|_{L^\infty} \sup_{h} \bigg\{ \epsilon_w(|h|)^{-1} \int_{[0,1]^2} \big| w (\zeta,\xi) - w (\zeta,\xi - h) \big| \; \rd \zeta \rd \xi \bigg\}\bigg)
\\
& \quad \leq \max \big( \|f\|_{L^\infty} , \|\partial_x f\|_{L^\infty} , \|\partial_x f\|_{L^1} \big) \|w\|_{L^\infty L^1} + \|f\|_{L^\infty}.
\end{aligned}
\end{equation*}
Therefore, we have that
\begin{equation*}
\begin{aligned}
\int_0^1 \bigg|\widetilde{\bm{H}}^N(t,\xi) - H(t,\xi)\bigg| \;\rd\xi &\leq \Big( \max \big( \|f\|_{L^\infty} , \|\partial_x f\|_{L^\infty} , \|\partial_x f\|_{L^1} \big) \|w\|_{L^\infty L^1} + \|f\|_{L^\infty} \Big)
\\
& \quad \int_0^t \|\mu^N(s,\cdot,\cdot) - \mu(s,\cdot,\cdot)\|_{\Phi_w^{-1}}  \rd s.
\end{aligned}
\end{equation*}

Now, we turn to the distance between $\widetilde{\bm{H}}^N$ and $\bar {\bm H}^N$, which differ only through the kernels $w^N$ and $w$. 
By considering $w^N$ and $w$ as $L^\infty \to L^1$ operators (from the $\zeta$ domain to the $\xi$ domain), we directly have that
\begin{equation*}
\begin{aligned}
& \int_0^1 \bigg|\bar {\bm H}^N(t,\xi) - \widetilde{\bm{H}}^N(t,\xi)\bigg| \;\rd\xi 
\\
& \quad = \int_{[0,1]} \bigg| \int_0^t\int_{[0,1]} \big[ w^N(\xi,\zeta) - w (\xi,\zeta) \big] \int_{\R} f(x)\bm{\mu}^N(t,\zeta,\rd x) \, \rd\zeta \rd s \bigg| \;\rd \xi
\\
& \quad \leq \big\|w^N - w \big\|_{L^\infty \to L^1} \big\|{\textstyle \int_0^t \int_{\R} f(x)\bm{\mu}^N(t,\cdot,\rd x)\, \rd s} \big\|_{L^\infty}
\\
& \quad \leq \big\|w^N - w \big\|_{L^\infty \to L^1} \|f\|_{L^\infty} t.
\end{aligned}
\end{equation*}

A triangular inequality  concludes the proof.
\end{proof}

Next, we compare $\bar{\bm{H}}^N$ and $\bm{H}^N$; this is where the Poisson random measures $\{\bm{\Pi}^{i;N}(\rd t, \rd z)\}_{i=1}^N$ come into play. Establishing a deterministic or a.s. bound on the distance would be difficult and inefficient. However, because of the $1/N$ scaling in the definition of $\bm{H}^N$, the expectation of the distance should be controllable by some sort of law of large numbers argument. This is what is shown in the following lemma, which uses It\^o isometry. 

\begin{lemma} \label{lem:difference_flow_SPDE_VPDE}

The following estimate holds:
\begin{equation*} \label{eqn:difference_flow_SPDE_VPDE}
\begin{aligned}
\E\Big[\|\bm{H}^N(t,\cdot) - \bar{\bm{H}}^N(t,\cdot)\|_{L^1_\xi}\Big] \leq \frac{1}{\sqrt{N}}\left(\|w^N\|_{L^\infty}^2 \|f\|_{L^\infty} t\right)^{1/2}.
\end{aligned}
\end{equation*}

\end{lemma}

\begin{proof}

Let us fully expand the definition of $\bm{H}^N$
\begin{equation*}
\begin{aligned}
\bm{H}^N(t,\xi) = \sum_{i=1}^N \bigg( \frac{1}{N} \sum_{j=1}^N \int_{[0,t]\times \R_+} w^N_{i,j} \mathbbm{1}_{\{z\leq f(\bm{X}^{j;N}(s-))\}}\bm{\Pi}^{j;N}(\rd s, \rd z) \bigg) \mathbbm{1}_{E^{i;N}}(\xi).
\end{aligned}
\end{equation*}
Similarly, $\bar{\bm{H}}^N$ can be written as
\begin{equation*}
\begin{aligned}
\bar{\bm{H}}^N(t,\xi) = \sum_{i=1}^N \bigg( \frac{1}{N} \sum_{j=1}^N \int_{[0,t]\times \R_+} w^N_{i,j} \mathbbm{1}_{\{z\leq f(\bm{X}^{j;N}(s-))\}}\rd s \rd z \bigg) \mathbbm{1}_{E^{i;N}}(\xi).
\end{aligned}
\end{equation*}
Hence, the $L^1_\xi$ distance between these two functions can be reduced to a discrete $\ell^1$ distance on $1 \leq i \leq N$:
\begin{equation*}
\begin{aligned}
&\int_0^1 \bigg|\bm{H}^N(t,\xi) - \bar{\bm{H}}^N(t,\xi)\bigg| \;\rd\xi
\\
& \quad = \frac{1}{N}\sum_{i=1}^N \bigg| \frac{1}{N} \sum_{j=1}^N \int_{[0,t]\times \R_+} w^N_{i,j} \mathbbm{1}_{\{z\leq f(\bm{X}^{j;N}(s-))\}} \big[ \bm{\Pi}^{j;N}(\rd s, \rd z) - \rd s \rd z \big] \bigg|.
\end{aligned}
\end{equation*}
The expectation can then be bounded by Jensen's inequality:
\begin{equation*}
\begin{aligned}
& \E\Big[\|\bm{H}^N(t,\cdot) - \bar{\bm{H}}^N(t,\cdot)\|_{L^1_\xi}\Big]
\\
& \quad = \E \Bigg[ \frac{1}{N}\sum_{i=1}^N \bigg| \frac{1}{N} \sum_{j=1}^N \int_{[0,t]\times \R_+} w^N_{i,j} \mathbbm{1}_{\{z\leq f(\bm{X}^{j;N}(s-))\}} \big[ \bm{\Pi}^{j;N}(\rd s, \rd z) - \rd s \rd z \big] \bigg| \Bigg]
\\
& \quad \leq \frac{1}{N}\sum_{i=1}^N \E \Bigg[ \bigg| \frac{1}{N} \sum_{j=1}^N \int_{[0,t]\times \R_+} w^N_{i,j} \mathbbm{1}_{\{z\leq f(\bm{X}^{j;N}(s-))\}} \big[ \bm{\Pi}^{j;N}(\rd s, \rd z) - \rd s \rd z \big] \bigg|^2 \Bigg]^{\frac{1}{2}}.
\end{aligned}
\end{equation*}
Since $\{\bm{X}^{i;N}\}_{i = 1}^N$ is an adapted process, we can apply the It\^o isometry for Poisson random measures to simplify the expectations above, and get
\begin{equation*}
\begin{aligned}
& \E\Big[\|\bm{H}^N(t,\cdot) - \bar{\bm{H}}^N(t,\cdot)\|_{L^1_\xi}\Big] 
\\
& \quad \leq \frac{1}{N}\sum_{i=1}^N \E \Bigg[ \bigg| \frac{1}{N} \sum_{j=1}^N \int_{[0,t]\times \R_+} w^N_{i,j} \mathbbm{1}_{\{z\leq f(\bm{X}^{j;N}(s-))\}} \big[ \bm{\Pi}^{j;N}(\rd s, \rd z) - \rd s \rd z \big] \bigg|^2 \Bigg]^{\frac{1}{2}}
\\
& \quad = \frac{1}{N}\sum_{i=1}^N \E \Bigg[ \sum_{j=1}^N \int_{[0,t]\times \R_+} \Big( \frac{1}{N} w^N_{i,j} \mathbbm{1}_{\{z\leq f(\bm{X}^{j;N}(s-))\}} \Big)^2 \rd s \rd z \Bigg]^{\frac{1}{2}}
\\
& \quad \leq \frac{1}{N}\sum_{i=1}^N \Bigg( \sum_{j=1}^N \Big( \frac{1}{N} w^N_{i,j} \Big)^2 \|f\|_{L^\infty} t \Bigg)^{\frac{1}{2}}
\\
& \quad \leq \frac{1}{\sqrt{N}}\left(\|w^N\|_{L^\infty}^2 \|f\|_{L^\infty} t\right)^{1/2},
\end{aligned}
\end{equation*}
which concludes the proof.
\end{proof}

Combining the two lemmas, we have that
\begin{equation*} 
\begin{aligned}
& \E \Big[ \|\bm{H}^N(t,\cdot) - H(t,\cdot)\|_{L^1} \Big]
\\
& \quad \leq \kappa_1 \int_0^t \E \Big[ \|\bm{\mu}^N(s,\cdot,\cdot) - \mu(s,\cdot,\cdot)\|_{\Phi_w^{-1}} \Big] \rd s + \|w^N - w \|_{L^\infty \to L^1} \|f\|_{L^\infty} t
\\
& \qquad + \frac{1}{\sqrt{N}}\left(\|w^N\|_{L^\infty}^2 \|f\|_{L^\infty} t\right)^{1/2}.
\end{aligned}
\end{equation*}
This concludes the first part of the main stability estimate.

\subsection{Stability via duality} \label{subsec:stability_backward}

We now consider the term $\E[\|\bm{\mu}_\#^N(t,\cdot,\cdot) - \mu_\#(t,\cdot,\cdot) \|_{\Phi_w^{-1}}]$.
Let $\{\bm{X}^{i;N}\}_{i = 1}^N$ be a strong solution to the system~\eqref{eqn:system} with $\{w^N_{i,j}\}_{i,j=1}^N \in \R^{N \times N}$, as provided by Proposition~\ref{prop:SDE_strong_solution}, and let $\mu$ be a solution to the mean-field PDE~\eqref{eqn:Vlasov} in the sense of characteristics with $w \in L^\infty([0,1])$, as provided by Proposition~\ref{prop:extended_empirical_measure}.
Then, $\bm{\mu}_\#^N$ and $\mu_\#$, the measures after the change of variable, satisfy \eqref{eqn:SPDE_2} and \eqref{eqn:Vlasov_4} respectively.
By subtracting \eqref{eqn:Vlasov_4} from \eqref{eqn:SPDE_2}, we immediately derive the following lemma.
\begin{lemma} \label{lem:backward}
Grant all the assumptions and take all the definitions from Propositions~\ref{prop:extended_empirical_measure} and \ref{prop:SDE_strong_solution}. Define the linear operator $D^H$ for $\varphi \in L^\infty([0,1]; C_b^1([0,t_*] \times \R))$ as
\begin{equation} \label{eqn:backward}
\begin{aligned}
D^H \varphi(s,\xi,x) \defeq b_\#(s,\xi,x) \partial_x \varphi(s,\xi,x) + f_\#(s,\xi,x)\left(\varphi(s,\xi,- H(s,\xi)) - \varphi(s,\xi,x)\right).
\end{aligned}
\end{equation}
Then, a.s., for all $\varphi \in L^\infty([0,1]; C_b^1([0,t_*] \times \R))$ and $t \in [0,t_*]$,
\begin{equation*}
\begin{aligned}
& \int_{[0,1] \times \R} \varphi(t,\xi,x) \Big( \bm{\mu}_\#^N(t,\xi,\rd x) - \mu_\#(t,\xi,\rd x) \Big) \;\rd \xi 
\\
&\quad = \int_{[0,1] \times \R} \varphi(0,\xi,x) \Big( \bm{\mu}_\#^N(0,\xi,\rd x) - \mu_\#(0,\xi,\rd x) \Big) \;\rd \xi  +
\bm{J_D}(\varphi) + \bm{J_H}(\varphi) + \bm{J_F}(\varphi),
\end{aligned}
\end{equation*}
where
\begin{equation} \label{eqn:backward_remainder}
\begin{aligned}
\bm{J_D}(\varphi) & \defeq \int_{0}^t \int_{[0,1] \times \R} (\partial_s + D^H) \varphi(s,\xi,x) \Big( \bm{\mu}_\#^N(s,\xi,\rd x) - \mu_\#(s,\xi,\rd x) \Big) \;\rd \xi  \rd s,
\\
\bm{J_H}(\varphi) & \defeq \int_{0}^t \int_{[0,1] \times \R} \Big[ \Big( \bm{b}_\#^N(s,\xi,x) - b_\#(s,\xi,x) \Big) \partial_x \varphi(s,\xi,x)
\\
\;& \quad\quad\quad\quad - \Big( \bm{f}_\#^N(s,\xi,x) - f_\#(s,\xi,x) \Big) \varphi(s,\xi,x) + \bm{f}_\#^N(s,\xi,x) \varphi(s,\xi,- \bm{H}^N(s,\xi)) 
\\
& \hspace{5cm} - f_\#(s,\xi,x) \varphi(s,\xi,- H(s,\xi)) \Big] \bm{\mu}_\#^N(s,\xi,\rd x) \;\rd \xi  \rd s,
\\
\bm{J_F}(\varphi) & \defeq \int_{(\xi,x)\in[0,1] \times \R} \int_{(s,z)\in[0,t]\times\R_+} \varphi(s,\xi,x) \big[ - \bm{\mu}_\#^N(s-,\xi,\rd x) + \bm{\delta}_\#^N(s-,\xi,\rd x) \big]
\\
& \quad\quad\quad\quad\quad\quad\quad\quad
\sum_{i=1}^N\mathbbm{1}_{E^{i;N}}(\xi) \mathbbm{1}_{\{z\leq \bm{r}^N(s-,\xi)\}} \big[ \bm{\Pi}^{i;N}(\rd s, \rd z) - \rd s \rd z \big]
\;\rd \xi.
\end{aligned}
\end{equation}

\end{lemma}

Through Lemma~\ref{lem:backward} and the definition of the metric $\Phi_w^{-1}$, we can readily obtain that, a.s.,
\begin{equation} \label{eqn:backward_remainder_optimize}
\begin{aligned}
& \|\bm{\mu}_\#^N(t,\cdot,\cdot) - \mu_\#(t,\cdot,\cdot) \|_{\Phi_w^{-1}}
\\
&= \sup_{\| \phi \|_{\Phi_w} \leq 1, \phi \in L^\infty_\xi C^1_{x}} \int_{[0,1] \times \R} \phi(\xi,x) \Big( \bm{\mu}_\#^N(t,\xi,\rd x) - \mu_\#(t,\xi,\rd x) \Big) \; \rd \xi 
\\
&= \sup_{\| \phi \|_{\Phi_w} \leq 1, \phi \in L^\infty_\xi C^1_{x}} \inf_{\varphi(t,\cdot,\cdot) = \phi, \varphi \in L^\infty_\xi C^1_{t,x}} \Bigg[ \int_{[0,1] \times \R} \varphi(0,\xi,x) \Big( \bm{\mu}_\#^N(0,\xi,\rd x) - \mu_\#(0,\xi,\rd x) \Big) \;\rd \xi 
\\
& \quad \hphantom{= \sup_{\| \phi \|_{\Phi_w} \leq 1, \phi \in L^\infty_\xi C^1_{x}} \inf_{\varphi(t,\cdot,\cdot) = \phi, \varphi \in L^\infty_\xi C^1_{t,x}} \Bigg[}
+ \bm{J_D}(\varphi) + \bm{J_H}(\varphi) + \bm{J_F}(\varphi) \Bigg].
\end{aligned}
\end{equation}
Note that in the original definition of the norm $\|.\|_{\Phi_w^{-1}}$, we do not impose that $\phi \in L^\infty_\xi C^1_{x}$. Of course, there exist $\phi$ that are not $C^1_x$-differentiable but still satisfy $\| \phi \|_{\Phi_w} \leq 1$. Nonetheless, it is easy to verify that the subset of $C^1_x$ functions is dense in $\Phi_w$. Therefore, it is sufficient to take the supremum over this dense subset. 

The decomposition in \eqref{eqn:backward_remainder} and \eqref{eqn:backward_remainder_optimize} can be interpreted as follows: $\bm{J_D}(\varphi)$ corresponds to the definition of the distributional solution presented in \eqref{eqn:Vlasov_4}, and vanishes if $\varphi$ is a $C^1_{t,x}$-solution of the adjoint equation $(\partial_t + D^H) \varphi(t,\xi,x) = 0$ (which we will later call the dual-backward equation). 
On the other hand,  $\bm{\mu}_\#^N$ does not exactly solve \eqref{eqn:Vlasov_4}, leading to the additional terms $\bm{J_H}(\varphi)$ and $\bm{J_F}(\varphi)$. Like in Section~\ref{subsec:stability_input}, we put in $\bm{J_H}(\varphi)$ the parts that can be bounded a.s., and leave in $\bm{J_F}(\varphi)$ the effect of the random resets of the membrane potentials due to the Poisson random measures.

\subsection{Propagation of regularity along the dual-backward equation}

In this subsection, we show the propagation of the $\Phi_w$ norm along the adjoint equation backward in time. The ${\Phi_w^{-1}}$-stability of measures follows from a primal-dual argument on the right-hand side of \eqref{eqn:backward_remainder_optimize}.

\begin{lemma} \label{lem:propagation_of_regularity_backward}
For any $H \in L^\infty([0,1]; C([0,t_*]))$, $\phi \in L^\infty([0,1];C_b^1(\R))$, and any $t \in [0,t_*]$, there exists a unique solution to the dual-backward equation in the sense of characteristics,
\begin{equation} \label{eqn:backward_equation}
\begin{aligned}
(\partial_s + D^H) \varphi(s,\xi,x) & = 0, \quad s \in [0,t],
\\
\varphi(t,\xi,x) & = \phi(\xi,x),
\end{aligned}
\end{equation}
such that $\varphi \in  L^\infty([0,1];C_b^1([0,t] \times \R))$. Assume that $\|\phi\|_{\Phi_w} \leq 1$. Then, by choosing 
\begin{equation*}
\begin{aligned}
& \kappa_2(t) \defeq \exp\Big[ \Big(2\|f\|_{L^\infty} + 2\|\partial_x f\|_{L^1} + 2 \|\partial_x f\|_{L^\infty} + \|\partial_x b\|_{L^\infty} \Big) t\Big],
\\
& \kappa_3(t) \defeq \exp \Big( 2 \|f\|_{L^\infty} t \Big) \bigg\{ 1 + \frac{t^2}{2} \|f\|_{L^\infty} \Big( \|\partial_x b\|_{L^\infty} + 2 \|\partial_x f\|_{L^\infty} + \|f\|_{L^\infty} \Big) \kappa_2(t) \bigg\},
\end{aligned}
\end{equation*}
we have 
\begin{subequations}
\begin{align} \label{eqn:regularity_x}
\sup_{s \in [0,t]} \max \Big( \|\varphi(s,\cdot,\cdot)\|_{L^\infty_{\xi,x}} , \|\partial_x \varphi(s,\cdot,\cdot)\|_{L^\infty_{\xi}L^1_x} , \|\partial_x \varphi(s,\cdot,\cdot)\|_{L^\infty_{\xi,x}} \Big) \leq \kappa_2(t ),
\\ \label{eqn:regularity_xi}
\int_{[0,1]} \sup_{s \in [0,t]} \sup_{x \in \R} \Big| \varphi(s,\xi - h,x) - \varphi(s,\xi,x) \Big| \;\rd \xi \leq \epsilon_w(|h|) \kappa_3(t).
\end{align}
\end{subequations}
As a consequence,
\begin{equation*}
\begin{aligned}
\sup_{s \in [0,t]} \| \varphi(s,\cdot,\cdot) \|_{\Phi_w} \leq \kappa(t) \defeq \max\big( \kappa_2(t) , \kappa_3(t) \big),
\end{aligned}
\end{equation*}
and
\begin{equation*}
\begin{aligned}
\| \partial_t \varphi \|_{L^\infty_{t,\xi,x}} \leq \kappa_4(t) \defeq (\|b\|_{L^\infty} + 2\|f\|_{L^\infty}) \kappa_2(t).
\end{aligned}
\end{equation*}

\end{lemma}

\begin{proof}

For a.e. $\xi \in [0,1]$, we can associate a unique and well-defined flow solving the system of characteristic equations
\begin{equation*}
\left\{
\begin{aligned}
& \frac{\rd}{\rd s} Y_H(t,s,\xi,x) = b_\#\big(s,\xi,Y_H(t,s,\xi,x)\big),
\\
& Y_H(t,t,\xi,x) = x.
\end{aligned} \right.
\end{equation*}
By the method of characteristics, we can derive the identity
\begin{equation*}
\begin{aligned}
\varphi(\tau, \xi, x)
& = \phi\big(\xi,Y_H(\tau,t,\xi,x)\big) \exp \bigg( - \int_\tau^t f_\#\big(s,\xi,Y_H(\tau,s,\xi,x)\big) \;\rd s \bigg)
\\
& \quad + \int_\tau^t \varphi\big( r,\xi,-H(r,\xi) \big) \exp \bigg( - \int_\tau^r f_\#\big(s,\xi,Y_H(\tau,s,\xi,x)\big) \;\rd s \bigg) \;\rd r,
\end{aligned}
\end{equation*}
from which we can prove the uniqueness of a solution $\varphi \in L^\infty([0,1]; C_b^1([0,t_*] \times \R))$ through a contraction argument, using the iteration
\begin{equation*}
\begin{aligned}
\varphi^{(k)}(\tau, \xi, x)
& = \phi\big(\xi,Y_H(\tau,t,\xi,x)\big) \exp \bigg( - \int_\tau^t f_\#\big(s,\xi,Y_H(\tau,s,\xi,x)\big) \;\rd s \bigg)
\\
& \quad + \int_\tau^t \varphi^{(k-1)}\big( r,\xi,-H(r,\xi) \big) \exp \bigg( - \int_\tau^r f_\#\big(s,\xi,Y_H(\tau,s,\xi,x)\big) \;\rd s \bigg) \;\rd r.
\end{aligned}
\end{equation*}

While the detailed regularity estimate in the $x$-direction can also be deduced via the method of characteristics, we only present here a formal calculation due to the length a rigorous argument would take; turning the following calculation into a rigorous argument would not pose any major difficulty, since $\varphi \in L^\infty([0,1]; C_b^1([0,t_*] \times \R))$.

We omit the $t,\xi$ variables for brevity. The test function $\varphi$ satisfies the adjoint equation,

\begin{equation*}
\begin{aligned}
\partial_t \varphi(x) + b_\#(x) \partial_x \varphi(x) - f_\#(x) \varphi(x) + f_\#(x) \varphi(- H) = 0,
\end{aligned}
\end{equation*}
and its $x$-derivative satisfies
\begin{equation*}
\begin{aligned}
\partial_t \varphi_x(x) + \partial_x \big[ b_\#(x) \varphi_x(x) \big] - \partial_x f_\#(x) \big[ \varphi(x) - \varphi(- H) \big] - f_\#(x) \varphi_x(x) = 0.
\end{aligned}
\end{equation*}
It is easy to verify that
\begin{equation*}
\begin{aligned}
\bigg| \frac{\rd}{\rd t} \|\varphi(t,\cdot,\cdot)\|_{L^\infty_{\xi,x}} \bigg| & \leq 2\|f\|_{L^\infty}\|\varphi(t,\cdot,\cdot)\|_{L^\infty_{\xi,x}}
\\
\bigg| \frac{\rd}{\rd t} \|\partial_x \varphi(t,\cdot,\cdot)\|_{L^\infty_{\xi}L^1_x} \bigg| & \leq 2\|\partial_x f\|_{L^1} \|\varphi(t,\cdot,\cdot)\|_{L^\infty_{\xi,x}} + \|f\|_{L^\infty} \|\partial_x \varphi(t,\cdot,\cdot)\|_{L^\infty_{\xi}L^1_x}
\\
\bigg| \frac{\rd}{\rd t} \|\partial_x \varphi(t,\cdot,\cdot)\|_{L^\infty_{\xi,x}} \bigg| & \leq \|\partial_x b\|_{L^\infty} \|\partial_x \varphi(t,\cdot,\cdot)\|_{L^\infty_{\xi,x}} + 2\|\partial_x f\|_{L^\infty} \|\varphi(t,\cdot,\cdot)\|_{L^\infty_{\xi,x}}
\\
& \quad + \|f\|_{L^\infty} \|\partial_x \varphi(t,\cdot,\cdot)\|_{L^\infty_{\xi,x}},
\end{aligned}
\end{equation*}
which yields the regularity \eqref{eqn:regularity_x} in the $x$-direction.

Next, we consider the regularity in the $\xi$-direction. For a.e. $\xi \in [0,1]$,
\begin{equation*}
\begin{aligned}
& \bigg| \frac{\rd}{\rd t} \sup_{x \in \R} \Big| \varphi(t,\xi - h,x) - \varphi(t,\xi,x) \Big| \bigg|
\\
& \quad \leq \sup_{x \in \R} \Big| \big[ b_\#(t,\xi - h,x) - b_\#(t,\xi,x) \big] \partial_x \varphi(t,\xi,x) \Big|
\\
& \qquad + \sup_{x \in \R} \Big| f_\#(t,\xi-h,x) \varphi(t,\xi-h,x) - f_\#(t,\xi,x) \varphi(t,\xi,x) \Big|
\\
& \quad \quad + \sup_{x \in \R} \Big| f_\#(t,\xi-h,x) \varphi(t,\xi-h,- H(t,\xi-h)) - f_\#(t,\xi,x) \varphi(t,\xi,- H(t,\xi)) \Big|
\\
& \quad \leq 2 \|f\|_{L^\infty} \sup_{x \in \R} \Big| \varphi(t,\xi - h,x) - \varphi(t,\xi,x) \Big| 
\\
& \quad \quad + \Big( \|\partial_x b\|_{L^\infty} \|\partial_x \varphi(t,\cdot,\cdot)\|_{L^\infty_{\xi,x}} + 2 \|\partial_x f\|_{L^\infty} \|\varphi(t,\cdot,\cdot)\|_{L^\infty_{\xi,x}} + \|f\|_{L^\infty} \|\partial_x \varphi(t,\cdot,\cdot)\|_{L^\infty_{\xi,x}} \Big)
\\
& \hspace{2cm} |H(t,\xi-h) - H(t,\xi)|.
\end{aligned}
\end{equation*}
Taking the supremum over the interval $[0,t]$, integrating over $\xi \in [0,1]$, and then applying Gronwall's lemma, we obtain
\begin{multline}\label{eq:lem6_gronwall}
\int_{[0,1]} \sup_{s \in [0,t]} \sup_{x \in \R} \Big| \varphi(s,\xi - h,x) - \varphi(s,\xi,x) \Big| \;\rd \xi
\\
\leq \bigg\{ \int_{[0,1]} \sup_{x \in \R} \Big| \varphi(t,\xi - h,x) - \varphi(t,\xi,x) \Big| \;\rd \xi + \Big( \|\partial_x b\|_{L^\infty} \|\partial_x \varphi(t,\cdot,\cdot)\|_{L^\infty_{\xi,x}}
\\
+ 2 \|\partial_x f\|_{L^\infty} \|\varphi(t,\cdot,\cdot)\|_{L^\infty_{\xi,x}} + \|f\|_{L^\infty} \|\partial_x \varphi(t,\cdot,\cdot)\|_{L^\infty_{\xi,x}} \Big) \\
\times \int_0^t \int_{[0,1]} |H(s,\xi-h) - H(s,\xi)| \;\rd \xi \rd s \bigg\}\exp \Big( 2 \|f\|_{L^\infty} t \Big).
\end{multline}
The difference $|H(t,\xi-h) - H(t,\xi)|$ can be bounded by the inherent regularity of $w$ characterized by the modulus of continuity $\epsilon_w : [0,1] \to [0,\infty)$. To establish this, we first observe that
\begin{equation*}
\begin{aligned}
H(t,\xi-h) - H(t,\xi) & = \int_0^t \int_{[0,1]} \big[ w(\xi - h, \zeta) - w(\xi, \zeta) \big] \int_{\R} f(x) \mu(s,\zeta,\rd x) \; \rd \zeta \rd s
\\
& \leq t \|f\|_{L^\infty} \int_{[0,1]} \big| w(\xi - h, \zeta) - w(\xi, \zeta) \big| \;\rd \zeta.
\end{aligned}
\end{equation*}
Integrating over the domain $[0,t] \times [0,1]$, we encounter a double integral spanning $[0,1]^2$, which allows us to use the definition of $\epsilon_w$, yielding
\begin{equation*}
\begin{aligned}
& \int_0^t \int_{[0,1]} |H(s,\xi-h) - H(s,\xi)| \;\rd \xi \rd s 
\\
& \quad \leq \int_0^t \int_{[0,1]} s \|f\|_{L^\infty} \int_{[0,1]} \big| w(\xi - h, \zeta) - w(\xi, \zeta) \big| \;\rd \zeta \rd \xi \rd s
\\
& \quad = \frac{t^2}{2} \|f\|_{L^\infty} \int_{[0,1]^2} \big| w(\xi - h, \zeta) - w(\xi, \zeta) \big| \;\rd \zeta \rd \xi
\\
& \quad \leq \frac{t^2}{2} \|f\|_{L^\infty} \epsilon_w(|h|).
\end{aligned}
\end{equation*}
Plugging this bound into the Gronwall inequality~\eqref{eq:lem6_gronwall}, we obtain the regularity in the $\xi$-direction, as expressed in \eqref{eqn:regularity_xi}:
\begin{equation*}
\begin{aligned}
& \int_{[0,1]} \sup_{s \in [0,t]} \sup_{x \in \R} \Big| \varphi(s,\xi - h,x) - \varphi(s,\xi,x) \Big| \;\rd \xi
\\
& \quad \leq
\bigg\{ \int_{[0,1]} \sup_{x \in \R} \Big| \varphi(t,\xi - h,x) - \varphi(t,\xi,x) \Big| \;\rd \xi + \Big( \|\partial_x b\|_{L^\infty} + 2 \|\partial_x f\|_{L^\infty} + \|f\|_{L^\infty} \Big) \times
\\
& \quad \quad \max\Big(\|\varphi(t,\cdot,\cdot)\|_{L^\infty_{\xi,x}} , \|\partial_x \varphi(t,\cdot,\cdot)\|_{L^\infty_{\xi,x}} \Big) \frac{t^2}{2} \|f\|_{L^\infty} \epsilon_w(|h|) \bigg\} \exp \Big( 2 \|f\|_{L^\infty} t \Big)
\\
& \quad \leq \epsilon_w(|h|) \kappa_3(t) \| \varphi(t,\cdot,\cdot) \|_{\Phi_w}
\end{aligned}
\end{equation*}

Finally, combining the two regularity propagation estimates, we conclude that
\begin{equation*}
\begin{aligned}
\sup_{s \in [0,t]} \| \varphi(s,\cdot,\cdot) \|_{\Phi_w} \leq \max \big( \kappa_2(t) , \kappa_3(t) \big) \| \varphi(t,\cdot,\cdot) \|_{\Phi_w}.
\end{aligned}
\end{equation*}
The last estimate,
\begin{equation*}
\begin{aligned}
\| \partial_t \varphi \|_{L^\infty_{t,\xi,x}} \leq (\|b\|_{L^\infty} + 2\|f\|_{L^\infty}) \kappa_2(t) \| \varphi(t,\cdot,\cdot) \|_{\Phi_w},
\end{aligned}
\end{equation*}
can also easily be derived from the dual-backward equation~\eqref{eqn:backward_equation}.
\end{proof}

To estimate the sup-inf expression on the right-hand side of \eqref{eqn:backward_remainder_optimize}, given any $\phi \in L^\infty([0,1];C_b^1(\R))$ such that $\|\phi\|_{\Phi_w} \leq 1$, we choose $\varphi$ as the solution to the dual-backward equation~\eqref{eqn:backward_equation}. With this choice, we ensure that $\bm{J_D}(\varphi) = 0$ and have an a priori estimates provided by Lemma~\ref{lem:propagation_of_regularity_backward}, which gives
\begin{equation*}
\begin{aligned}
& \sup_{\| \phi \|_{\Phi_w} \leq 1, \phi \in L^\infty_\xi C^1_{x}, \; \varphi \textnormal{ solves } \eqref{eqn:backward_equation}} \Bigg[ \int_{[0,1] \times \R} \varphi(0,\xi,x) \Big( \bm{\mu}_\#^N(0,\xi,\rd x) - \mu_\#(0,\xi,\rd x) \Big) \;\rd \xi \Bigg]
\\
& \quad \leq \kappa(t) \|\bm{\mu}_\#^N(0,\cdot,\cdot) - \mu_\#(0,\cdot,\cdot) \|_{\Phi_w^{-1}}.
\end{aligned}
\end{equation*}
In addition, we can bound $\bm{J_H}(\varphi)$ as
\begin{equation*}
\begin{aligned}
&|\bm{J_H}(\varphi)| \\
&\quad\leq \int_0^t \int_{[0,1] \times \R} \Big[ \Big| \bm{b}_\#^N(s,\xi,x) - b_\#(s,\xi,x) \Big| |\partial_x \varphi(s,\xi,x)| \\
&\quad\quad- \Big| \bm{f}_\#^N(s,\xi,x) - f_\#(s,\xi,x) \Big| |\varphi(s,\xi,x)|
\\
&\quad\quad + \big|\bm{f}_\#^N(s,\xi,x) \varphi(s,\xi,- \bm{H}^N(s,\xi)) - f_\#(s,\xi,x) \varphi(s,\xi,- H(s,\xi))\big| \Big] \bm{\mu}_\#^N(s,\xi,\rd x) \;\rd \xi \rd s
\\
& \quad \leq \int_0^t \int_{[0,1] \times \R} \big(\|\partial_x b\|_{L^\infty} + \|f\|_{L^\infty} + 2 \|\partial_x f\|_{L^\infty} \big) \kappa_2(t) | \bm{H}^N(s,\xi) - H(s,\xi) | \\
& \qquad\qquad\qquad\qquad\qquad\qquad\qquad\qquad\qquad\qquad\qquad\qquad\qquad\qquad\,\times\bm{\mu}_\#^N(s,\xi,\rd x) \;\rd \xi \rd s.
\end{aligned}
\end{equation*}
Integrate first over $x$ and then over $\xi$, we conclude that
\begin{equation*}
\begin{aligned}
\sup_{\| \phi \|_{\Phi_w} \leq 1, \phi \in L^\infty_\xi C^1_{x}, \; \varphi \textnormal{ solves } \eqref{eqn:backward_equation}} |\bm{J_H}(\varphi)| & \leq \kappa_5(t) \int_0^t \| \bm{H}^N(s,\cdot) - H(s,\cdot)\|_{L^1_\xi} \;\rd s,
\\
\kappa_5(t) & \defeq \big(\|\partial_x b\|_{L^\infty} + \|f\|_{L^\infty} + 2 \|\partial_x f\|_{L^\infty} \big) \kappa_2(t).
\end{aligned}
\end{equation*}
(Recall that a bound for $\E\left[\| \bm{H}^N(s,\cdot) - H(s,\cdot)\|_{L^1_\xi}\right]$ has already been derived in Lemmas~\ref{lem:difference_flow_measure}~and~\ref{lem:difference_flow_SPDE_VPDE}.)

\subsection{Effect of the random resets of the membrane potentials} \label{subsec:stability_Poisson}

It remains to estimate $\bm{J_F}(\varphi)$ in \eqref{eqn:backward_remainder_optimize}, which summarizes the effect of the random resets of the membrane potentials due to the randomness of the Poisson random measures.
Recall that for all $\varphi \in L^\infty([0,1]; C_b^1([0,t_*] \times \R))$,
\begin{equation} \label{eqn:backward_remainder_poisson}
\begin{aligned}
\bm{J_F}(\varphi) & \defeq \int_{(\xi,x)\in[0,1] \times \R} \int_{(s,z)\in[0,t]\times\R_+} \varphi(s,\xi,x) \big[ - \bm{\mu}_\#^N(s-,\xi,\rd x) + \bm{\delta}_\#^N(s-,\xi,\rd x) \big]
\\
& \quad\quad\quad\quad\quad\quad\quad\quad
\sum_{i=1}^N\mathbbm{1}_{E^{i;N}}(\xi) \mathbbm{1}_{\{z\leq \bm{r}^N(s-,\xi)\}} \big[ \bm{\Pi}^{i;N}(\rd s, \rd z) - \rd s \rd z \big]
\;\rd \xi,
\\
\end{aligned}
\end{equation}
where 
\begin{equation*}
\begin{aligned}
\varphi^{i;N}(s,x) = \int_{E^{i;N}}\varphi(s,\xi,x)\;\rd \xi.
\end{aligned}
\end{equation*}
We introduce the \textit{compensated} jump process,
\begin{equation*}
    \bar{Z}^{i;N}_t(\xi) = \int_{[0,t]\times\R_+} \mathbbm{1}_{\{z\leq \bm{r}^N(s-,\xi)\}} \big[ \bm{\Pi}^{i;N}(\rd s, \rd z) - \rd s \rd z \big],
\end{equation*}
so that we can write
\begin{multline*}
    \bm{J_F}(\varphi) \\
    = \int_{[0,t] \times [0,1] \times \R} \varphi(s,\xi,x) \big[ - \bm{\mu}_\#^N(s-,\xi,\rd x) + \bm{\delta}_\#^N(s-,\xi,\rd x) \big]\sum_{i=1}^N\mathbbm{1}_{E^{i;N}}(\xi) \rd \bar{Z}^{i;N}_s(\xi)\rd \xi.
\end{multline*}

Note that, for all $\xi\in[0,1]$ and $t\geq 0$, $\E[\bar{Z}^{i;N}_t(\xi)]=0$, which implies that $\E[\bm{J_F}(\varphi)]=0$. In analogy to Lemma~\ref{lem:difference_flow_SPDE_VPDE}, we can expect that $\bm{J_F}(\varphi)$ vanishes as $N\to\infty$ because of some sort of law of large numbers. In the following, we prove this intuition, using a mollification procedure and It\^o isometry. For notational convenience, we will write
\begin{equation*}
\bm{J_F}(\varphi) = \int_{[0,t] \times [0,1] \times \R} \varphi(s,\xi,x) \bm{F}^N(\rd s, \xi, \rd x) \;\rd \xi,
\end{equation*}
where
\begin{equation*}
    \bm{F}^N(\rd s, \xi, \rd x)\rd\xi = \big[ - \bm{\mu}_\#^N(s-,\xi,\rd x) + \bm{\delta}_\#^N(s-,\xi,\rd x) \big]
\sum_{i=1}^N\mathbbm{1}_{E^{i;N}}(\xi) \rd \bar{Z}^{i;N}_s(\xi)\rd\xi.
\end{equation*}

In order to interpolate between the ``oscillation'' of $\varphi$ and the randomness of $\bm{F}^N$, we introduce the mollification kernel $\chi$, which is a normalized indicator function with support size $[0, \Delta t] \times [-\Delta \xi, \Delta \xi] \times [-\Delta x, \Delta x]$:
\begin{equation} \label{eqn:indicator_kernel}
\begin{aligned}
\chi(s,\xi,x) & \defeq \chi_1(s) \chi_2(\xi) \chi_3(x) \defeq \frac{1}{\Delta t} \mathbbm{1}_{[0, \Delta t]}(s) \frac{1}{2 \Delta \xi }\mathbbm{1}_{[-\Delta \xi, \Delta \xi]}(\xi) \frac{1}{2\Delta x}\mathbbm{1}_{[-\Delta x, \Delta x]}(x),
\\
\chi^\dagger(s,\xi,x) & \defeq \chi(-s,-\xi,-x).
\end{aligned}
\end{equation}
Notably, $\chi \star \varphi$ involves a convolution along both the $\xi$- and $t$-directions, thereby requiring our embedding of $[0,1]$ in the torus (see Sec.~2.3) and an extended definition for $\varphi$ when $t < 0$.
As the definition of the dual operator $D^H$ only requires $H \in L^\infty([0,1]; C([0,t_*]))$, we can take the formal extension $H(s,\xi) = 0$ for negative time $s \in [-\Delta t,0]$ and solve \eqref{eqn:backward_equation} over $[-\Delta t, t]$. With this extension, all a priori estimates for $\varphi$ in Lemma~\ref{lem:propagation_of_regularity_backward} apply, with the sole modification that $\kappa(t)$ has to be replaced by $\kappa(t+ \Delta t)$. We can further decompose $\bm{J_F}(\varphi)$ into two parts by using $\varphi = [\varphi - (\chi \star \varphi)] + (\chi \star \varphi)$. The first part is bounded by
\begin{equation*}
\begin{aligned}
\bm{J_F}(\varphi - (\chi \star \varphi))
&= \int_{[0,t] \times [0,1] \times \R} \big[\varphi(s,\xi,x) - (\chi \star \varphi) (s,\xi,x) \big] \bm{F}^N(\rd s, \xi, \rd x) \;\rd \xi
\\
&\leq \bigg\{ \int_{[0,1]} \bigg( \sup_{s \in [0,t], x \in \R} \big|\varphi(s,\xi,x) - \big(\chi \star \varphi\big)(s,\xi,x) \big| \bigg)^2 \rd \xi \bigg\}^{\frac{1}{2}}
\\
& \qquad \bigg\{ \int_{[0,1]} \bigg( \int_{[0,t]\times\R} \bm{F}^N(\rd s, \xi, \rd x) \bigg)^2 \rd \xi \bigg\}^{\frac{1}{2}}.
\end{aligned}
\end{equation*}
For any moment function $\eta: \R \to \R_+$ such that $\|\eta^{-1}\|_{L^1} < \infty$, the second part is bounded by
\begin{equation*}
\begin{aligned}
& \bm{J_F}(\chi \star \varphi) = \int_{[0,t] \times [0,1] \times \R} (\chi \star \varphi) (s,\xi,x) \bm{F}^N(\rd s, \xi, \rd x) \rd \xi
\\
& \quad = \int_{-\infty}^{+\infty} \int_{[0,1] \times \R} \varphi(s,\xi,x) \mathbbm{1}_{[-\Delta t, t]}(s) \Big[\chi^\dagger \star \big(\bm{F}^N \mathbbm{1}_{[0,t] \times [0,1] \times \R} \big)\Big] (s,\xi,x) \;\rd \xi \rd x \rd s
\\
& \quad\leq \bigg\{ \int_{-\infty}^{+\infty} \int_{[0,1] \times \R} \big( \varphi(s,\xi,x) \big)^2 \mathbbm{1}_{[-\Delta t, t]}(s) \eta^{-1}(x) \;\rd \xi \rd x \rd s \bigg\}^{\frac{1}{2}}
\\
& \quad \quad \bigg\{ \int_{-\infty}^{+\infty} \int_{[0,1] \times \R} \bigg( \Big[\chi^\dagger \star \big(\bm{F}^N \mathbbm{1}_{[0,t] \times [0,1] \times \R} \big)\Big] (s,\xi,x) \bigg)^2 \eta(x) \;\rd \xi \rd x \rd s \bigg\}^{\frac{1}{2}}.
\end{aligned}
\end{equation*}
The two factors involving $\varphi$ above can be respectively bounded by
\begin{equation*}
\begin{aligned}
& \sup_{\| \phi \|_{\Phi_w} \leq 1, \,\phi \in L^\infty_\xi C^1_{x}, \; \varphi \textnormal{ solves } \eqref{eqn:backward_equation}} \bigg\{ \int_{[0,1]} \bigg( \sup_{s \in [0,t], x \in \R} \big|\varphi(s,\xi,x) - \big(\chi \star \varphi\big)(s,\xi,x) \big| \bigg)^2 \rd \xi \bigg\}^{\frac{1}{2}}
\\
& \quad \leq \sup_{\dots} \bigg[ \int_{[0,1]} \sup_{s,x} \big| \varphi(s,\xi,x) - \big(\chi \star \varphi\big)(s,\xi,x) \big| \;\rd \xi \bigg]^{\frac{1}{2}} \sup_{s,\xi,x} \big| \varphi(s,\xi,x) - \big(\chi \star \varphi\big)(s,\xi,x) \big|^{\frac{1}{2}}
\\
& \quad \leq \sup_{\dots} \bigg[ \|\partial_t \varphi\|_{L^\infty} \Delta t + \|\partial_x \varphi\|_{L^\infty} \Delta x + \|\varphi\|_{\Phi_w} \epsilon_w(\Delta \xi) \bigg]^{\frac{1}{2}} \big( 2\|\varphi\|_{L^\infty} \big)^{\frac{1}{2}}
\\
& \quad \leq \max \big( \kappa(t + \Delta t) , \kappa_4(t + \Delta t) \big) \sqrt{2(\Delta t + \Delta x + \epsilon_w(\Delta \xi) )},
\end{aligned}
\end{equation*}
and
\begin{equation*}
\begin{aligned}
& \sup_{\| \phi \|_{\Phi_w} \leq 1, \phi \in L^\infty_\xi C^1_{x}, \; \varphi \textnormal{ solves } \eqref{eqn:backward_equation}} \bigg\{ \int_{-\infty}^{+\infty} \int_{[0,1] \times \R} \big( \varphi(s,\xi,x) \big)^2 \mathbbm{1}_{[-\Delta t, t]}(s) \eta^{-1}(x) \;\rd \xi \rd x \rd s \bigg\}^{\frac{1}{2}}
\\
& \quad \leq \sup_{\dots} \Big\{ \|\varphi\|_{L^\infty}^2 \|\eta^{-1}\|_{L^1} (t + \Delta t) \Big\}^{\frac{1}{2}}
\\
& \quad \leq \kappa(t + \Delta t) \Big( \|\eta^{-1}\|_{L^1}(t + \Delta t) \Big)^{\frac{1}{2}}.
\end{aligned}
\end{equation*}

The mollification procedure perfomed above reduces the proof of the convergence of $\bm{J_F}(\varphi)$ to zero as $N\to\infty$ to the following lemma, which relies on It\^o isometry.  
\begin{lemma} \label{lem:difference_renewal_SPDE_VPDE_2}
\begin{subequations}

The following estimate holds:
\begin{equation} \label{eqn:Ito_F_1}
\begin{aligned}
\E \Bigg[ \bigg\{ \int_{[0,1]} \bigg( \int_{[0,t]\times\R} \bm{F}^N(\rd s, \xi, \rd x) \bigg)^2 \rd \xi \bigg\}^{\frac{1}{2}} \Bigg] \leq \Big( 16\|f\|_{L^\infty}^2 t^2 + 4 \|f\|_{L^\infty} t \Big)^{\frac{1}{2}}.
\end{aligned}
\end{equation}
Moreover, given the a priori moment bound
\begin{equation} \label{eqn:moment_bound_Ito}
\begin{aligned}
\sup_{s \in [0,t]} \frac{1}{N} \sum_{i=1}^N \E\bigg[ \eta\big( |\bm{X}^{i;N}(s-) - \bm{H}^{i;N}(s-)|+ \Delta x \big) + \eta\big( |\bm{H}^{i;N}(s-) | + \Delta x \big) \bigg] \leq M_\eta(t),
\end{aligned}
\end{equation}
the following estimate holds:
\begin{multline} \label{eqn:Ito_F_2}
\E \Bigg[ \bigg\{ \int_{\R\times[0,1]\times\R} \bigg( \Big[\chi^\dagger \star \big(\bm{F}^N \mathbbm{1}_{[0,t] \times [0,1] \times \R} \big)\Big] (s,\xi,x) \bigg)^2 \eta(x) \;\rd \xi \rd x \rd s \bigg\}^{\frac{1}{2}} \Bigg]
\\
\leq \frac{1}{\sqrt{N}} \Bigg( \frac{\|f\|_{L^\infty} t}{2 \Delta t \Delta \xi \Delta x} M_\eta(t) \Bigg)^{\frac{1}{2}}.
\end{multline}

\end{subequations}
\end{lemma}

\begin{proof}
Proof of \eqref{eqn:Ito_F_1}: It is easy to verify that
\begin{equation*}
\begin{aligned}
E_1 & \defeq \E \Bigg[ \bigg\{ \int_{[0,1]} \bigg( \int_{[0,t]\times\R} \bm{F}^N(\rd s, \xi, \rd x) \bigg)^2 \rd \xi \bigg\}^{\frac{1}{2}} \Bigg]
\\
& \leq \E \bigg[ \int_{[0,1]} \bigg( \int_{[0,t]\times\R} \Big| \bm{F}^N(\rd s, \xi, \rd x) \Big| \bigg)^2 \rd \xi \bigg]^{\frac{1}{2}}
\\
& \leq \E \bigg[ \int_{[0,1]} \bigg( \int_{[0,t]\times\R_+} \int_{\R} \big[ \bm{\mu}_\#^N(s-,\xi,\rd x) + \bm{\delta}_\#^N(s-,\xi,\rd x) \big]
\\
& \quad\quad\quad\quad\quad\quad\quad\quad
\sum_{i=1}^N\mathbbm{1}_{E^{i;N}}(\xi) \mathbbm{1}_{\{z\leq \bm{r}^N(s-,\xi)\}} \big[ \bm{\Pi}^{i;N}(\rd s, \rd z) + \rd s \rd z \big] \bigg)^2
\;\rd \xi \bigg]^{\frac{1}{2}}.
\end{aligned}
\end{equation*}
Since $\|\bm{r}^N\|_{L^\infty} \leq \|f\|_{L^\infty}$ holds a.s., we can consider the integration $\rd z$ over $z \in [0, \|f\|_{L^\infty}] \subset \R_+$.
Integrating first over $\rd x$ and then over $\rd z$, we obtain the following:
\begin{equation*}
\begin{aligned}
E_1 & \leq \E \bigg[ \int_{[0,1]} \bigg( \int_{[0,t] \times [0,\|f\|_{L^\infty}]} 2
\sum_{i=1}^N\mathbbm{1}_{E^{i;N}}(\xi) \big[ \bm{\Pi}^{i;N}(\rd s, \rd z) + \rd s \rd z \big] \bigg)^2
\;\rd \xi \bigg]^{\frac{1}{2}}
\\
& = \E \bigg[ \frac{1}{N} \sum_{i=1}^N \bigg(2 \int_{[0,t] \times [0,\|f\|_{L^\infty}]} \big[ \bm{\Pi}^{i;N}(\rd s, \rd z) + \rd s \rd z \big] \bigg)^2 \bigg]^{\frac{1}{2}}
\\
& = \Big( 16\|f\|_{L^\infty}^2 t^2 + 4 \|f\|_{L^\infty} t \Big)^{\frac{1}{2}},
\end{aligned}
\end{equation*}
where the last line is a straightforward computation of the moment of a Poisson random measure.

Proof of \eqref{eqn:Ito_F_2}: 
Replace the variables $(s,\xi,x)$ in the integration of \eqref{eqn:Ito_F_2} by $(\tau,\theta,y)$, we get

\begin{equation*}
\begin{aligned}
E_2 & \defeq \E \Bigg[ \bigg\{ \int_{-\infty}^{+\infty} \int_{[0,1] \times \R} \bigg( \Big[\chi^\dagger \star \big(\bm{F}^N \mathbbm{1}_{[0,t] \times [0,1] \times \R} \big)\Big] (\tau,\theta,y) \bigg)^2 \eta(y) \;\rd \theta \rd y \rd \tau \bigg\}^{\frac{1}{2}} \Bigg]
\\
& \leq \E \bigg[ \int_{-\infty}^{+\infty} \int_{[0,1] \times \R} \bigg( \Big[\chi^\dagger \star \big(\bm{F}^N \mathbbm{1}_{[0,t] \times [0,1] \times \R} \big)\Big] (\tau,\theta,y) \bigg)^2 \eta(y) \;\rd \theta \rd y \rd \tau \bigg]^{\frac{1}{2}}
\\
& = \E\bigg[ \int_{(-\infty,\infty) \times [0,1] \times \R} \bigg( \int_{[0,1]} \int_{(-\infty,\infty) \times \R_+} \int_{\R} \chi(s - \tau, \xi - \theta, x - y) \; \mathbbm{1}_{[0,t]}(s)
\\
& \hspace{7cm} \big[ \bm{\mu}_\#^N(s-,\xi,\rd x) - \bm{\delta}_\#^N(s-,\xi,\rd x) \big]
\\
& \hphantom{= \E\bigg[}
\sum_{i=1}^N\mathbbm{1}_{E^{i;N}}(\xi) \mathbbm{1}_{\{z\leq \bm{r}^N(s-,\xi)\}} \big[ \bm{\Pi}^{i;N}(\rd s, \rd z) - \rd s \rd z \big] \;\rd \xi \bigg)^2 \eta(y) \;\rd \tau \rd \theta \rd y \bigg]^{\frac{1}{2}}.
\end{aligned}
\end{equation*}
Here, the random measures have finite total variation with probability one. Hencey, Fubini's theorem for measures applies a.s., leading to
\begin{equation*}
\begin{aligned}
E_2 & \leq \E\bigg[ \int_{(-\infty,\infty) \times [0,1] \times \R} \bigg( \sum_{i=1}^N \int_{[0,t] \times \R_+}
\\
& \hphantom{\leq \E\bigg[}
\bigg(\int_{E^{i;N}} \chi\big(s - \tau, \xi - \theta, \bm{Y}^{i;N}(s-) - y\big) - \chi\big(s - \tau, \xi - \theta, - \bm{H}^{i;N}(s-) - y\big) \;\rd \xi \bigg)
\\
& \hphantom{\leq \E\bigg[}
\qquad\qquad\qquad\qquad\quad\mathbbm{1}_{\{z\leq f(\bm{X}^{i;N}(s-))\}} \big[ \bm{\Pi}^{i;N}(\rd s, \rd z) - \rd s \rd z \big] \bigg)^2 \eta(y) \;\rd \tau \rd \theta \rd y \bigg]^{\frac{1}{2}}.
\end{aligned}
\end{equation*}
By It\^o isometry, we have
\begin{equation*}
\begin{aligned}
E_2 & \leq \E\bigg[ \int_{(-\infty,\infty) \times [0,1] \times \R} \sum_{i=1}^N \int_{[0,t] \times \R_+}
\\
& \hphantom{\leq \E\bigg[}
\bigg(\int_{E^{i;N}} \chi\big(s - \tau, \xi - \theta, \bm{Y}^{i;N}(s-) - y\big) - \chi\big(s - \tau, \xi - \theta, - \bm{H}^{i;N}(s-) - y\big) \;\rd \xi \bigg)^2
\\
& \hphantom{\leq \E\bigg[}
\qquad\qquad\qquad\qquad\qquad\qquad\qquad\qquad\quad\mathbbm{1}_{\{z\leq f(\bm{X}^{i;N}(s-))\}} \rd s \rd z \; \eta(y) \;\rd \tau \rd \theta \rd y \bigg]^{\frac{1}{2}}.
\end{aligned}
\end{equation*}
Again, the integration in $z$ can be restricted over $z \in [0, \|f\|_{L^\infty}] \subset \R_+$ so that 
\begin{equation*}
\begin{aligned}
E_2 & \leq \E\bigg[ \int_{(-\infty,\infty) \times [0,1] \times \R} \sum_{i=1}^N \int_{[0,t]}
\\
& 
\big(\chi_1(s - \tau)\big)^2 \bigg(\int_{E^{i;N}} \chi_2\big(\xi - \theta\big) \;\rd \xi \bigg)^2 \Big(\chi_3\big(\bm{Y}^{i;N}(s-) - y\big) - \chi_3\big(-\bm{H}^{i;N}(s-) - y\big)\Big)^2
\\
& \hphantom{\leq \E\bigg[}
\qquad\qquad\qquad\qquad\qquad\qquad\qquad\qquad\qquad\qquad\qquad\|f\|_{L^\infty} \rd s \; \eta(y) \;\rd \tau \rd \theta \rd y \bigg]^{\frac{1}{2}}.
\end{aligned}
\end{equation*}
Taking the integral over $s$ last and separating the variables $\tau, \theta$ and $y$ in the integration, we see that
\begin{equation*}
\begin{aligned}
E_2 & \leq \E\bigg[ \|f\|_{L^\infty} \int_{[0,t]} \bigg( \int_{(-\infty,\infty)} \big(\chi_1(s - \tau)\big)^2 \; \rd \tau \bigg)
\\
& \hspace{2.9cm}
\sum_{i=1}^N \bigg( \int_{[0,1]} \bigg(\int_{E^{i;N}} \chi_2\big(\xi - \theta\big) \;\rd \xi \bigg)^2 \;\rd \theta \bigg)
\\
& \hspace{2.9cm} \bigg( \int_{\R} \Big(\chi_3\big(\bm{Y}^{i;N}(s-) - y\big) - \chi_3\big(-\bm{H}^{i;N}(s-) - y\big)\Big)^2 \eta(y) \rd y \bigg)
\;\rd s \bigg]^{\frac{1}{2}}.
\end{aligned}
\end{equation*}
By the definition of $\chi_1$, we have
\begin{equation*}
\begin{aligned}
 \int_{(-\infty,\infty)} \big(\chi_1(s - \tau)\big)^2 \; \rd \tau = \frac{1}{\Delta t}.
\end{aligned}
\end{equation*}
Noticing that $E^{i;N} = 1/N$, $1\leq i\leq N$, we also have
\begin{equation*}
\begin{aligned}
\int_{[0,1]} \bigg(\int_{E^{i;N}} \chi_2\big(\xi - \theta\big) \;\rd \xi \bigg)^2 \;\rd \theta
& \leq \int_{[0,1]} \frac{1}{N} \int_{E^{i;N}} \big( \chi_2\big(\xi - \theta\big) \big)^2 \;\rd \xi \;\rd \theta
\\
& = \frac{1}{N^2 (2 \Delta \xi)},
\end{aligned}
\end{equation*}
where the final bound is independent of the index $i\in\{1,\dots,N\}$. 

Finally, by the fact that $\chi_3$ has compact support $[-\Delta x, \Delta x]$, we obtain that
\begin{equation*}
\begin{aligned}
& \E\bigg[ \int_{[0,t]} \sum_{i=1}^N \int_{\R} \Big(\chi_3\big(\bm{Y}^{i;N}(s-) - y\big) - \chi_3\big(-\bm{H}^{i;N}(s-) - y\big)\Big)^2 \eta(y) \rd y \rd s\bigg]
\\
& \leq 2\sum_{i=1}^N \E\bigg[ \int_{[0,t]} \int_{\R} \Big(\chi_3\big(\bm{Y}^{i;N}(s-) - y\big)\Big)^2 \eta(y) \rd y \rd s
\\
&\hspace{3cm} + \int_{[0,t]} \int_{\R} \Big(\chi_3\big(-\bm{H}^{i;N}(s-) - y\big)\Big)^2 \eta(y) \rd y \rd s \bigg]
\\
& \leq 2\sum_{i=1}^N \E\bigg[ \int_{[0,t]} \eta\big( |\bm{Y}^{i;N}(s-)|+ \Delta x \big)  \int_{\R} \Big(\chi_3\big(\bm{Y}^{i;N}(s-) - y\big)\Big)^2 \rd y \rd s
\\
& \hphantom{2\sum_{i=1}^N \E\bigg[ }
+ \int_{[0,t]} \eta\big( |\bm{H}^{i;N}(s-) | + \Delta x \big)  \int_{\R} \Big(\chi_3\big(-\bm{H}^{i;N}(s-) - y\big)\Big)^2 \rd y \rd s \bigg]
\\
& = \frac{2}{(2\Delta x)} \sum_{i=1}^N \E\bigg[ \int_{[0,t]}\eta\big( |\bm{Y}^{i;N}(s-)|+ \Delta x \big) + \eta\big( |\bm{H}^{i;N}(s-) | + \Delta x \big) ds\bigg].
\end{aligned}
\end{equation*}
Applying the three inequalities to our previous calculation, we conclude that
\begin{equation*}
\begin{aligned}
&E_2  \defeq \E \Bigg[ \bigg\{ \int_{-\infty}^{+\infty} \int_{[0,1] \times \R} \bigg( \Big[\chi^\dagger \star \big(\bm{F}^N \mathbbm{1}_{[0,t] \times [0,1] \times \R} \big)\Big] (\tau,\theta,y) \bigg)^2 \eta(y) \;\rd \theta \rd y \rd \tau \bigg\}^{\frac{1}{2}} \Bigg]
\\
&\leq \frac{1}{\sqrt{N}} \Bigg( \frac{\|f\|_{L^\infty} t}{2 \Delta t \Delta \xi \Delta x} \sup_{s \in [0,t]} \frac{1}{N} \sum_{i=1}^N \E\bigg[ \eta\big( |\bm{Y}^{i;N}(s-)|+ \Delta x \big) + \eta\big( |\bm{H}^{i;N}(s-) | + \Delta x \big) \bigg] \Bigg)^{\frac{1}{2}},
\end{aligned}
\end{equation*}
which is equivalent to~\eqref{eqn:Ito_F_2} by substituting $\bm{Y}^{i;N}(s-) = \bm{X}^{i;N}(s-) - \bm{H}^{i;N}(s-)$.
\end{proof}


\subsection{Concluding the proof} \label{subsec:stability_conclude}


\begin{proof}[Proof of Theorem~\ref{thm:main_metric}]
Let $(\bm{X}^{i;N})_{i = 1}^N$ be a strong solution to the system~\eqref{eqn:system} with $(w^N_{i,j})_{i,j=1}^N \in \R^{N \times N}$, as provided by Proposition~\ref{prop:SDE_strong_solution}, and let $\mu$ be a solution to the mean-field PDE~\eqref{eqn:Vlasov} in the sense of characteristics with $w \in L^\infty([0,1])$, as provided by Proposition~\ref{prop:extended_empirical_measure}.
We recall the inequality established in Proposition~\ref{prop:input_separation},
\begin{equation*} 
\begin{aligned}
\|\bm{\mu}^N(t,\cdot,\cdot) - \mu(t,\cdot,\cdot) \|_{\Phi_w^{-1}} \leq \|\bm{\mu}_\#^N(t,\cdot,\cdot) - \mu_\#(t,\cdot,\cdot) \|_{\Phi_w^{-1}} + \|\bm{H}^N(t,\cdot) - H(t,\cdot)\|_{L^1_\xi}.
\end{aligned}
\end{equation*}
The first term on the right hand side a.s. satisfies \eqref{eqn:backward_remainder_optimize}, restated here,
\begin{equation*}
\begin{aligned}
& \|\bm{\mu}_\#^N(t,\cdot,\cdot) - \mu_\#(t,\cdot,\cdot) \|_{\Phi_w^{-1}}
\\
& \quad = \sup_{\| \phi \|_{\Phi_w} \leq 1, \phi \in L^\infty_\xi C^1_{t,x}} \int_{[0,1] \times \R} \phi(\xi,x) \Big( \bm{\mu}_\#^N(t,\xi,\rd x) - \mu_\#(t,\xi,\rd x) \Big) \; \rd \xi 
\\
& \quad = \sup_{\| \phi \|_{\Phi_w} \leq 1, \phi \in L^\infty_\xi C^1_{x}, \; \varphi \textnormal{ solves } \eqref{eqn:backward_equation}} \Bigg[ \int_{[0,1] \times \R} \varphi(0,\xi,x) \Big( \bm{\mu}_\#^N(0,\xi,\rd x) - \mu_\#(0,\xi,\rd x) \Big) \;\rd \xi
\\
& \quad \hphantom{= \sup_{\| \phi \|_{\Phi_w} \leq 1, \phi \in L^\infty_\xi C^1_{x}} \inf_{\varphi(t,\cdot,\cdot) = \phi, \varphi \in L^\infty_\xi C^1_{t,x}} \Bigg[} +
\bm{J_D}(\varphi) + \bm{J_H}(\varphi) + \bm{J_F}(\varphi) \Bigg].
\end{aligned}
\end{equation*}
By Lemma~\ref{lem:propagation_of_regularity_backward}, we have that, a.s.,
\begin{equation*}
\begin{aligned}
 \|\bm{\mu}_\#^N(t,\cdot,\cdot) - \mu_\#(t,\cdot,\cdot) \|_{\Phi_w^{-1}} & \leq \kappa(t) \|\bm{\mu}_\#^N(0,\cdot,\cdot) - \mu_\#(0,\cdot,\cdot) \|_{\Phi_w^{-1}}
\\
& \quad + \kappa_5(t) \int_0^t \| \bm{H}^N(s,\cdot) - H(s,\cdot)\|_{L^1_\xi} \;\rd s + L(\bm{F}^N; \chi),
\end{aligned}
\end{equation*}
where the choice of $\chi$ depends on $\Delta t, \Delta \xi, \Delta x$ and
\begin{multline*}
L(\bm{F}^N; \chi)\defeq \max \big( \kappa(t + \Delta t) , \kappa_4(t + \Delta t) \big)\sqrt{2(\Delta t + \Delta x + \epsilon_w(\Delta \xi) )} \\
\times\bigg\{ \int_{[0,1]} \bigg( \int_0^t \int_{\R} \bm{F}^N(\rd s, \xi, \rd x) \bigg)^2 \rd \xi \bigg\}^{\frac{1}{2}}
+ \kappa(t + \Delta t) \Big( \|\eta^{-1}\|_{L^1}(t + \Delta t) \Big)^{\frac{1}{2}}\\
\times\bigg\{ \int_{-\infty}^{+\infty} \int_{[0,1] \times \R} \bigg( \Big[\chi^\dagger \star \big(\bm{F}^N \mathbbm{1}_{[0,t] \times [0,1] \times \R} \big)\Big] (s,\xi,x) \bigg)^2 \eta(x) \;\rd \xi \rd x \rd s \bigg\}^{\frac{1}{2}}
\end{multline*}
(see Sec.~\ref{subsec:stability_Poisson}). 

The expectation of the norm of $\bm{H}^N(t,\cdot) - H(t,\cdot)$ is controlled by Lemmas~\ref{eqn:difference_flow_measure} and \ref{lem:difference_flow_SPDE_VPDE}:
\begin{multline*} 
\E \Big[ \|\bm{H}^N(t,\cdot) - H(t,\cdot)\|_{L^1_\xi} \Big] \\
\leq \kappa_1 \int_0^t \E \Big[ \|\bm{\mu}^N(s,\cdot,\cdot) - \mu(s,\cdot,\cdot)\|_{\Phi_w^{-1}} \Big] \rd s + \|w^N - w \|_{L^\infty \to L^1} \|f\|_{L^\infty} t
\\
+ \frac{1}{\sqrt{N}}\left(\|w^N\|_{L^\infty}^2 \|f\|_{L^\infty} t\right)^{1/2}.
\end{multline*}

Lastly, choosing $\Delta t = \Delta \xi = \Delta x = N^{-\frac{1}{4}}$ and applying Lemma~\ref{lem:difference_renewal_SPDE_VPDE_2}, we have that
\begin{equation*}
\begin{aligned}
\E \Big[ L(\bm{F}^N; \chi) \Big] \leq \kappa_6(t) \sqrt{\max\Big( \epsilon\big(N^{-\frac{1}{4}}\big), N^{-\frac{1}{4}} \Big)},
\end{aligned}
\end{equation*}
where
\begin{equation*}
\begin{aligned}
\kappa_6(t) &\defeq \Big( 16\|f\|_{L^\infty}^2 t^2 + 2 \|f\|_{L^\infty} t \Big)^{\frac{1}{2}} \max \big( \kappa(t + 1) , \kappa_4(t + 1) \big) \sqrt{6}
\\
& \quad + \kappa(t + 1) \Big( \|\eta^{-1}\|_{L^1}(t + 1) \Big)^{\frac{1}{2}}
\Bigg( \frac{\|f\|_{L^\infty} t}{2} M_\eta(t) \Bigg)^{\frac{1}{2}}.
\end{aligned}
\end{equation*}

Combining all estimates, we obtain that
\begin{equation} \label{eqn:stability_main}
\begin{aligned}
& \E [\| \bm{\mu}^N(t,\cdot,\cdot) - \mu(t,\cdot,\cdot) \|_{\Phi_w^{-1}} ]
\\
& \quad \leq \kappa(t) \E \bigg[ \| \bm{\mu}^N(0,\cdot,\cdot) - \mu(0,\cdot,\cdot) \|_{\Phi_w^{-1}} \bigg] + \kappa_6(t) \sqrt{\max\Big( \epsilon\big(N^{-\frac{1}{4}}\big), N^{-\frac{1}{4}} \Big)} 
\\
& \quad \quad + \kappa_5(t) \int_0^t \Bigg( \frac{1}{\sqrt{N}}\left(\|w^N\|_{L^\infty}^2 \|f\|_{L^\infty} \tau\right)^{1/2} + \|w^N - w \|_{L^\infty \to L^1} \|f\|_{L^\infty} \tau
\\
& \qquad \qquad + \kappa_1 \int_0^\tau \E\bigg[ \|\bm{\mu}^N(s,\cdot,\cdot) , \mu(s,\cdot,\cdot)\|_{\Phi_w^{-1}} \bigg] \rd s \Bigg) \;\rd \tau
\\
& \quad \quad \frac{1}{\sqrt{N}}\left(\|w^N\|_{L^\infty}^2 \|f\|_{L^\infty} t\right)^{1/2} + \|w^N - w \|_{L^\infty \to L^1} \|f\|_{L^\infty} t
\\
& \qquad + \kappa_1 \int_0^t \E\bigg[ \|\bm{\mu}^N(s,\cdot,\cdot) , \mu(s,\cdot,\cdot)\|_{\Phi_w^{-1}} \bigg] \rd s.
\end{aligned}
\end{equation}

By another Gronwall estimate, we conclude that \begin{equation*}
\begin{aligned}
\E [\| \bm{\mu}^N(t,\cdot,\cdot) - \mu(t,\cdot,\cdot) \|_{\Phi_w^{-1}} ] \to 0, \quad \forall t \geq 0,
\end{aligned}
\end{equation*}
provided that
\begin{equation*}
\begin{aligned}
N \to \infty, \quad \|w^N - w\|_{L^\infty \to L^1} \to 0, \quad \text{ and } \quad 
\E \bigg[ \| \bm{\mu}^N(0,\cdot,\cdot) - \mu(0,\cdot,\cdot) \|_{\Phi_w^{-1}} \bigg] \to 0.
\end{aligned}
\end{equation*}
\end{proof}

It only remains to prove Corollary~\ref{cor:trajectory} on the trajectories.
\begin{proof} [Proof of Corollary~\ref{cor:trajectory}]
Let us first restate the definitions of $H$ and $\bm{H}^N$, and introduce further notations simplifying our calculations.
Let $\{E^{i;N}\}_{i = 1}^N$ be an almost everywhere partition.
The integrated postsynaptic input $\bm{H}^{i;N}$ and the extension $\bm{H}^N$ to $\xi \in [0,1]$ are also defined as
\begin{equation*}
\begin{aligned}
\bm{H}^N(t,\xi) & \defeq \textstyle \sum_{i=1}^N\bm{H}^{i;N}(t)\mathbbm{1}_{E^{i;N}}(\xi),
\\
\bm{H}^{i;N}(t) & \defeq \frac{1}{N} \sum_{j=1}^N \int_{[0,t]\times \R_+} w^N_{i,j} \mathbbm{1}_{\{z\leq f(\bm{X}^{j;N}(s-))\}}\bm{\Pi}^{j;N}(\rd s, \rd z).
\end{aligned}
\end{equation*}
The limit input $H$ and its formal restriction $H^{i;N}$ to $i \in \{1, \dots, N\}$ can be defined as
\begin{equation*}
\begin{aligned}
& H(t,\xi) \defeq \int_0^t\int_{[0,1]} w(\xi,\zeta) \int_{\R} f(x)\mu(t,\zeta,\rd x) \, \rd\zeta \rd s, 
\\
& H^{i;N}(t) \defeq \int_{[0,1]} H(t,\xi) \ N \mathbbm{1}_{E^{i;N}}(\xi) \;\rd \xi. 
\end{aligned}
\end{equation*}

With the above notations, the equation for $\bm{X}^{i;N}$ and $\widetilde{\bm{X}}^{i;N}$ can be rewritten respectively as
\begin{equation*}
\begin{aligned}
\bm{X}^{i;N}(t) & = \bm{X}^{i;N}_0 + \int_0^t b(\bm{X}^{i;N}(s)) \rd s + \bm{H}^{i;N}(t)
\\
& \quad - \int_{[0,t]\times \R_+} \bm{X}^{i;N}(s-) \mathbbm{1}_{\{z\leq f(\bm{X}^{i;N}(s-))\}}\bm{\Pi}^{i;N}(\rd s, \rd z), \quad \text{ a.s. } \ \forall t \geq 0,
\\
\widetilde{\bm{X}}^{i;N}(t) & = \bm{X}^{i;N}_0 + \int_0^t b(\widetilde{\bm{X}}^{i;N}(s)) \rd s + H^{i;N}(t)
\\
& \quad - \int_{[0,t]\times \R_+} \widetilde{\bm{X}}^{i;N}(s-) \mathbbm{1}_{\{z\leq f(\widetilde{\bm{X}}^{i;N}(s-))\}}\bm{\Pi}^{i;N}(\rd s, \rd z), \quad \text{ a.s. } \ \forall t \geq 0.
\end{aligned}
\end{equation*}
Subtracting the two equations, we have that, a.s., $\forall t \geq 0$,
\begin{equation*}
\begin{aligned}
& \bm{X}^{i;N}(t) - \widetilde{\bm{X}}^{i;N}(t)
\\
& \quad = \int_0^t \Big( b(\bm{X}^{i;N}(s)) - b(\widetilde{\bm{X}}^{i;N}(s)) \Big) \rd s + \bm{H}^{i;N}(t) - H^{i;N}(t)
\\
& \quad \quad - \int_{[0,t]\times \R_+} \Big( \bm{X}^{i;N}(s-) - \widetilde{\bm{X}}^{i;N}(s-) \Big)
\\
& \qquad \qquad \mathbbm{1}\big\{z\leq \min\big( f(\bm{X}^{i;N}(s-)), f(\widetilde{\bm{X}}^{i;N}(s-)) \big)\big\} \ \bm{\Pi}^{i;N}(\rd s, \rd z)
\\
& \quad \quad - \int_{[0,t]\times \R_+} \bm{X}^{i;N}(s-) \mathbbm{1} \big\{f(\widetilde{\bm{X}}^{i;N}(s-)) < z\leq f(\bm{X}^{i;N}(s-)) \big\} \ \bm{\Pi}^{i;N}(\rd s, \rd z)
\\
& \quad \quad - \int_{[0,t]\times \R_+} - \widetilde{\bm{X}}^{i;N}(s-) \mathbbm{1} \big\{f(\bm{X}^{i;N}(s-)) < z\leq f(\widetilde{\bm{X}}^{i;N}(s-))\big\} \ \bm{\Pi}^{i;N}(\rd s, \rd z).
\end{aligned}
\end{equation*}

In the following, we will use an approximation of the absolute value function that is better suited for dealing with stochastic jumps. Its definition and properties are presented in the following lemma.
%
\begin{lemma}\label{prop:smoothing_1} 
For any $\epsilon > 0$, define the regularized absolute value function 
\begin{equation*} 
\begin{aligned} 
|x|_{\epsilon} \defeq \big(x^2 + \epsilon^2\big)^{1/2}. 
\end{aligned} 
\end{equation*} 
It satisfies 
\begin{equation*} 
\begin{aligned} 
|x| \leq |x|_{\epsilon} \leq |x| + \epsilon,\quad
|x|_{\epsilon}' = \frac{x}{\big(x^2 + \epsilon^2\big)^{1/2}}, 
\quad 
|x|_{\epsilon}'' = \frac{\epsilon^2}{\big(x^2 + \epsilon^2\big)^{3/2}}. 
\end{aligned} 
\end{equation*} 
In addition, define the regularized saturating absolute value function 
\begin{equation*} 
\begin{aligned} 
\lambda_{\epsilon}(x) \defeq \frac{(|x|_{\epsilon}+1) - \sqrt{(|x|_{\epsilon}-1)^2 + \epsilon^2}}{2} + \epsilon.
\end{aligned} 
\end{equation*} 
It satisfies  
\begin{equation*} 
\begin{aligned} 
\min(|x|,1) &\leq \lambda_{\epsilon}(x) \leq \min(|x|,1) + 3\epsilon,\\ \lambda_{\epsilon}'(x) & = \frac{1}{2} \bigg(1 - \frac{|x|_{\epsilon}-1}{\sqrt{(|x|_{\epsilon}-1)^2 + \epsilon^2}}\bigg) |x|_{\epsilon}'\leq |x|_{\epsilon}' \leq 1, \\ 
\lambda_{\epsilon}''(x) & = \frac{1}{2} \bigg(1 - \frac{|x|_{\epsilon}-1}{\sqrt{(|x|_{\epsilon}-1)^2 + \epsilon^2}}\bigg) |x|_{\epsilon}'' - \frac{1}{2} \frac{\epsilon^2}{\big( (|x|_{\epsilon}-1)^2 + \epsilon^2\big)^{3/2}} \big(|x|_{\epsilon}'\big)^2 \\
&\leq \frac{1}{2} \bigg(1 - \frac{|x|_{\epsilon}-1}{\sqrt{(|x|_{\epsilon}-1)^2 + \epsilon^2}}\bigg) |x|_{\epsilon}'' \leq |x|_{\epsilon}''. 
\end{aligned} 
\end{equation*} 
Moreover, the finite difference is bounded from above by 
\begin{equation*} 
\begin{aligned} 
\lambda_{\epsilon}(x+u) & = \lambda_{\epsilon}(x) + \lambda_{\epsilon}'(x) u + \int_0^u \lambda_{\epsilon}''(x+r) r \; \rd r 
\\ 
& \leq \lambda_{\epsilon}(x) + \lambda_{\epsilon}'(x) u + \frac{u^2}{2 \epsilon}. 
\end{aligned} 
\end{equation*} 
\end{lemma} 
We now come back to our proof but, instead of trying to directly bound $|\bm{X}^{i;N}(t) - \widetilde{\bm{X}}^{i;N}(t)|$, we work with $\lambda_{\epsilon}\Big(\bm{X}^{i;N}(t) - \widetilde{\bm{X}}^{i;N}(t)\Big)$.
Applying It\^o's lemma, we have that
\begin{equation*}
\begin{aligned}
& \lambda_{\epsilon}\Big(\bm{X}^{i;N}(t) - \widetilde{\bm{X}}^{i;N}(t)\Big) = \lambda_{\epsilon}(0) + \bm{B}^{i;N}(t) + \bm{I}^{i;N}(t) + \bm{R}^{i;N}_0(t) + \bm{R}^{i;N}_1(t),
\end{aligned}
\end{equation*}
where
\begin{equation*}
\begin{aligned}
\bm{B}^{i;N}(t) & = \int_0^t \Big( b(\bm{X}^{i;N}(s)) - b(\widetilde{\bm{X}}^{i;N}(s)) \Big) \lambda_{\epsilon}'\Big(\bm{X}^{i;N}(s) - \widetilde{\bm{X}}^{i;N}(s)\Big) \ \rd s,
\\
\bm{I}^{i;N}(t) & = \sum_{j=1}^N \int_{[0,t]\times \R_+} \mathbbm{1}_{\{z\leq f(\bm{X}^{j;N}(s-))\}} \Big[\lambda_{\epsilon}\Big(\bm{X}^{i;N}(s-) - \widetilde{\bm{X}}^{i;N}(s-) + w^N_{i,j}/N \Big) 
\\
& \hphantom{+ \bigg\{ \sum_{j=1}^N \int_{[0,t]\times \R_+} \mathbbm{1}_{\{z\leq f(\bm{X}^{j;N}(s-))\}} \Big[}
- \lambda_{\epsilon}\Big(\bm{X}^{i;N}(s-) - \widetilde{\bm{X}}^{i;N}(s-) \Big) \Big] \bm{\Pi}^{j;N}(\rd s, \rd z)
\\
& \quad \quad + \int_0^t - \bigg( \int_{[0,1]} \int_{[0,1]} N \mathbbm{1}_{E^{i;N}}(\xi) w(\xi,\zeta) \int_{\R} f(y)\mu({s},\zeta,y) \rd y \rd\zeta \rd \xi\bigg) \
\\
& \hspace{7cm} \lambda_{\epsilon}'\Big(\bm{X}^{i;N}(s) - \widetilde{\bm{X}}^{i;N}(s)\Big) \ \rd s,
\end{aligned}
\end{equation*}
and
\begin{equation*}
\begin{aligned}
\bm{R}^{i;N}_0(t) & = \int_{[0,t]\times \R_+} \mathbbm{1}_{\big\{z\leq \min\big( f(\bm{X}^{i;N}(s-)), f(\widetilde{\bm{X}}^{i;N}(s-)) \big)\big\}}
\\
& \hspace{4.5cm}\Big[ \lambda_{\epsilon}(0) - \lambda_{\epsilon} \Big( \bm{X}^{i;N}(s-) - \widetilde{\bm{X}}^{i;N}(s-) \Big)\Big]\ \bm{\Pi}^{i;N}(\rd s, \rd z),
\\
\bm{R}^{i;N}_1(t) & = \int_{[0,t]\times \R_+} \mathbbm{1}_{\big\{f(\widetilde{\bm{X}}^{i;N}(s-)) < z\leq f(\bm{X}^{i;N}(s-)) \big\}}
\\
& \hspace{2.5cm} \Big[ \lambda_{\epsilon}\Big(- \widetilde{\bm{X}}^{i;N}(s-)\Big) - \lambda_{\epsilon} \Big( \bm{X}^{i;N}(s-) - \widetilde{\bm{X}}^{i;N}(s-) \Big)\Big] \ \bm{\Pi}^{i;N}(\rd s, \rd z)
\\
& \quad + \int_{[0,t]\times \R_+} \mathbbm{1}_{\big\{f(\bm{X}^{i;N}(s-)) < z\leq f(\widetilde{\bm{X}}^{i;N}(s-))\big\}}
\\
& \hspace{2.9cm} \Big[ \lambda_{\epsilon}\Big(\bm{X}^{i;N}(s-)\Big) - \lambda_{\epsilon} \Big( \bm{X}^{i;N}(s-) - \widetilde{\bm{X}}^{i;N}(s-) \Big)\Big] \ \bm{\Pi}^{i;N}(\rd s, \rd z).
\end{aligned}
\end{equation*}
The estimates for $\bm{B}^{i;N}(t),$ $\bm{R}^{i;N}_0(t)$, and $\bm{R}^{i;N}_1(t)$ are easy to derive and are given by
\begin{equation*}
\begin{aligned}
\sup_{t \in [0,t_*]} \bm{B}^{i;N}(t) & \leq \|b\|_{W^{1,\infty}} \int_0^{t_*} \min\Big( \big| \bm{X}^{i;N}(s) - \widetilde{\bm{X}}^{i;N}(s) \big|, 1 \Big) \ \rd s,
\\
\sup_{t \in [0,t_*]} \bm{R}^{i;N}_0(t) & \leq 0,
\\
\sup_{t \in [0,t_*]} \bm{R}^{i;N}_1(t) & \leq {(1+3\epsilon)} \int_{[0,t_*]\times \R_+} \mathbbm{1}_{\big\{f(\widetilde{\bm{X}}^{i;N}(s-)) < z\leq f(\bm{X}^{i;N}(s-)) \big\}} 
\\
& \hspace{3cm} + \mathbbm{1}_{\big\{f(\bm{X}^{i;N}(s-)) < z\leq f(\widetilde{\bm{X}}^{i;N}(s-)) \big\}} \ \bm{\Pi}^{i;N}(\rd s, \rd z).
\end{aligned}
\end{equation*}
Taking expectations, we have (for sufficient small $\epsilon > 0$)
\begin{equation*}
\begin{aligned}
\E\Big[\sup_{t \in [0,t_*]} \bm{B}^{i;N}(t)\Big] & \leq \|b\|_{W^{1,\infty}} \int_0^{t_*} \E\Big[\min\Big( \big| \bm{X}^{i;N}(s) - \widetilde{\bm{X}}^{i;N}(s) \big|, 1 \Big)\Big] \ \rd s,
\\
\E\Big[\sup_{t \in [0,t_*]} \bm{R}^{i;N}_0(t)\Big] & \leq 0,
\\
\E\Big[\sup_{t \in [0,t_*]} \bm{R}^{i;N}_1(t)\Big] & \leq {2} \|f\|_{W^{1,\infty}} \int_0^{t_*} \E\Big[\min\Big( \big| \bm{X}^{i;N}(s) - \widetilde{\bm{X}}^{i;N}(s) \big|, 1 \Big)\Big] \ \rd s.
\end{aligned}
\end{equation*}

The term $\bm{I}^{i;N}(t)$ quantifies the influence of $\bm{H}^{i;N}(t) - H^{i;N}(t)$ on the dynamics.
Following the strategy in Section~\ref{subsec:stability_input}, we consider the martingale
\begin{equation*}
\begin{aligned}
\bm{I}^{i;N}_M(t) & = \sum_{j=1}^N \int_{[0,t]\times \R_+} \mathbbm{1}_{\{z\leq f(\bm{X}^{j;N}(s-))\}} \Big[\lambda_{\epsilon}\Big(\bm{X}^{i;N}(s-) - \widetilde{\bm{X}}^{i;N}(s-) + w^N_{i,j}/N \Big) 
\\
& \hspace{3cm}
- \lambda_{\epsilon}\Big(\bm{X}^{i;N}(s-) - \widetilde{\bm{X}}^{i;N}(s-) \Big) \Big] \big[ \bm{\Pi}^{j;N} (\rd s, \rd z) - \rd s \rd z \big],
\end{aligned}
\end{equation*}
together with the difference due to the finite jump sizes
\begin{equation*}
\begin{aligned}
\bm{I}^{i;N}_1(t) & = \sum_{j=1}^N \int_{[0,t]\times \R_+} \mathbbm{1}_{\{z\leq f(\bm{X}^{j;N}(s-))\}} \Big[\lambda_{\epsilon}\Big(\bm{X}^{i;N}(s-) - \widetilde{\bm{X}}^{i;N}(s-) + w^N_{i,j}/N \Big) 
\\
& \quad \hphantom{+ \sum_{j=1}^N \int_{[0,t]\times \R_+} \mathbbm{1}_{\{z\leq f(\bm{X}^{j;N}(s-))\}} \Big[}
- \lambda_{\epsilon}\Big(\bm{X}^{i;N}(s-) - \widetilde{\bm{X}}^{i;N}(s-) \Big) \Big] \ \rd s \rd z
\\
& \quad \quad - \int_0^t \bigg( \int_{[0,1]} \int_{[0,1]} N \mathbbm{1}_{E^{i;N}}(\xi) w^N(\xi,\zeta) \int_{\R} f(y)\bm{\mu}^N(s,\zeta,\rd y)\,  \rd\zeta \rd \xi\bigg) \ 
\\
& \hspace{6cm} \lambda_{\epsilon}'\Big(\bm{X}^{i;N}(s) - \widetilde{\bm{X}}^{i;N}(s)\Big) \ \rd s,
\end{aligned}
\end{equation*}
and the error created by advection 
\begin{equation*}
\begin{aligned}
\bm{I}^{i;N}_2(t) & = \int_0^t \bigg( \int_{[0,1]} \int_{[0,1]} N \mathbbm{1}_{E^{i;N}}(\xi) w^N(\xi,\zeta) \int_{\R} f(x)\bm{\mu}^N(s,\zeta,\rd x) \, \rd\zeta \rd \xi\bigg) \ 
\\
& \hspace{6cm} \lambda_{\epsilon}'\Big(\bm{X}^{i;N}(s) - \widetilde{\bm{X}}^{i;N}(s)\Big) \ \rd s
\\
& \quad \quad - \int_0^t \bigg( \int_{[0,1]} \int_{[0,1]} N \mathbbm{1}_{E^{i;N}}(\xi) w(\xi,\zeta) \int_{\R} f(x)\mu(s,\zeta,\rd x) \, \rd\zeta \rd \xi\bigg) \ 
\\
& \hspace{6cm} \lambda_{\epsilon}'\Big(\bm{X}^{i;N}(s) - \widetilde{\bm{X}}^{i;N}(s)\Big) \ \rd s.
\end{aligned}
\end{equation*}

By Doob's martingale inequality and It\^o isometry,
\begin{equation*}
\begin{aligned}
& \E\Big[\sup_{t \in [0,t_*]} \big| \bm{I}^{i;N}_M(t) \big|^2 \Big] 
\\
& \quad \leq 4 \E\Big[\big| \bm{I}^{i;N}_M(t_*) \big|^2 \Big]
\\
& \quad \leq \E\bigg[ \sum_{j=1}^N \int_{[0,t_*]\times \R_+} \mathbbm{1}_{\{z\leq f(\bm{X}^{j;N}(s-))\}} \Big[\lambda_{\epsilon}\Big(\bm{X}^{i;N}(s-) - \widetilde{\bm{X}}^{i;N}(s-) + w^N_{i,j}/N \Big) 
\\
& \qquad \hphantom{\E\bigg[ \sum_{j=1}^N \int_{[0,t]\times \R_+} \mathbbm{1}_{\{z\leq f(\bm{X}^{j;N}(s-))\}} \Big[}
- \lambda_{\epsilon}\Big(\bm{X}^{i;N}(s-) - \widetilde{\bm{X}}^{i;N}(s-) \Big) \Big]^2 \ \rd s \rd z \bigg]
\\
& \quad \leq \sum_{j=1}^N \int_{[0,t_*]\times \R_+} \mathbbm{1}_{\{z\leq \|f\|_{L^\infty}\}} \frac{\big(\max_{1 \leq i,j \leq N}|w^N_{i,j}|\big)^2}{N^2} \ \rd s \rd z
\\
& \quad \leq \frac{\|w^N\|_{L^\infty}^2}{N} \|f\|_{L^\infty} t_*.
\end{aligned}
\end{equation*}

Rewriting the extended variable $\xi \in [0,1]$ as the index $i,\dots,N$, we have
\begin{equation*}
\begin{aligned}
\bm{I}^{i;N}_1(t) & = \sum_{j=1}^N \int_0^t \ f(\bm{X}^{j;N}(s-)) \ \Big[\lambda_{\epsilon}\Big(\bm{X}^{i;N}(s-) - \widetilde{\bm{X}}^{i;N}(s-) + w^N_{i,j}/N \Big) 
\\
& \quad \hphantom{+ \sum_{j=1}^N \int_0^t \ f(\bm{X}^{j;N}(s-)) \ \Big[}
- \lambda_{\epsilon}\Big(\bm{X}^{i;N}(s-) - \widetilde{\bm{X}}^{i;N}(s-) \Big) \Big] \ \rd s
\\
& \quad \quad - \sum_{j=1}^N \int_0^t \ f(\bm{X}^{j;N}(s-)) \ \lambda_{\epsilon}'\Big(\bm{X}^{i;N}(s) - \widetilde{\bm{X}}^{i;N}(s)\Big) w^N_{i,j}/N \ \rd s.
\end{aligned}
\end{equation*}
The supremum over $t \in [0,t_*]$ is hence bounded by
\begin{equation*}
\begin{aligned}
\sup_{t \in [0,t_*]} \bm{I}^{i;N}_1(t) & \leq \sum_{j=1}^N \int_0^{t_*} \|f\|_{L^\infty} \frac{\big(\max_{1 \leq i,j \leq N}|w^N_{i,j}|\big)^2}{2\epsilon N^2} \ \rd s
\\
& \leq \frac{\|w^N\|_{L^\infty}^2}{2\epsilon N} \|f\|_{L^\infty} t_*,
\end{aligned}
\end{equation*}
where, in the first inequality, we used Lemma~\ref{prop:smoothing_1}.

Taking the maximum over $t \in [0,t_*]$ of $\bm{I}^{i;N}_2$,
\begin{equation*}
\begin{aligned}
\sup_{t \in [0,t_*]} \bm{I}^{i;N}_2(t) & = \int_0^{t_*} \bigg| \int_{[0,1]} \int_{[0,1]} N \mathbbm{1}_{E^{i;N}}(\xi) w^N(\xi,\zeta) \int_{\R} f(x)\bm{\mu}^N(s,\zeta,\rd x)\, \rd\zeta \rd \xi
\\
& \quad \quad \quad \quad - \int_{[0,1]} \int_{[0,1]} N \mathbbm{1}_{E^{i;N}}(\xi) w(\xi,\zeta) \int_{\R} f(x)\mu(s,\zeta,\rd x) \, \rd\zeta \rd \xi\bigg| \ \rd s.
\end{aligned}
\end{equation*}
Now, by averaging over the indices $i$, we find that
\begin{equation*}
\begin{aligned}
\frac{1}{N} \sum_{i=1}^N \sup_{t \in [0,t_*]} \bm{I}^{i;N}_2(t) & = \int_0^{t_*} \int_{[0,1]} \bigg|  \int_{[0,1]} w^N(\xi,\zeta) \int_{\R} f(x)\bm{\mu}^N(s,\zeta,\rd x) \rd\zeta
\\
& \quad \quad \quad \quad \quad \quad - \int_{[0,1]} w(\xi,\zeta) \int_{\R} f(x)\mu(s,\zeta,\rd x)\, \rd\zeta \bigg| \ \rd \xi \rd s.
\end{aligned}
\end{equation*}
Following the argument in the proof of Lemma~\ref{eqn:difference_flow_measure}, we deduce that
\begin{equation*}
\begin{aligned}
&\int_{[0,1]} \bigg| \int_{[0,1]} w^N(\xi,\zeta) \int_{\R} f(x)\bm{\mu}^N(s,\zeta,\rd x)\, \rd\zeta - \int_{[0,1]} w(\xi,\zeta) \int_{\R} f(x)\mu(s,\zeta,\rd x) \rd\zeta \bigg| \; \rd \xi
\\
& \quad \leq \kappa_1 \|\bm{\mu}^N(s,\cdot,\cdot) - \mu(s,\cdot,\cdot)\|_{\Phi_w^{-1}} + \|w^N - w \|_{L^\infty \to L^1} \|f\|_{L^\infty},
\end{aligned}
\end{equation*}
where $\kappa_1 \defeq ( \max \big( \|f\|_{L^\infty} , \|\partial_x f\|_{L^\infty} , \|\partial_x f\|_{L^1} \big) \|w\|_{L^\infty L^1} + \|f\|_{L^\infty} )$.

Thus, we have the bound
\begin{equation*}
\begin{aligned}
& \E\bigg[\frac{1}{N} \sum_{i=1}^N \sup_{t \in [0,t_*]} \bm{I}^{i;N}_2(t)\bigg]
\\
& \quad \leq \kappa_1 \int_0^{t_*} \E \big[ \|\bm{\mu}^N(s,\cdot,\cdot) - \mu(s,\cdot,\cdot)\|_{\Phi_w^{-1}} \big] \rd s + \|w^N - w \|_{L^\infty \to L^1} \|f\|_{L^\infty} t_*.
\end{aligned}
\end{equation*}
Finally,
\begin{equation*}
\begin{aligned}
& \E\bigg[\frac{1}{N} \sum_{i=1}^N \sup_{t \in [0,t_*]} \min\Big( \big| \bm{X}^{i;N}(t) - \widetilde{\bm{X}}^{i;N}(t) \big|, 1 \Big) \bigg]
\\
& \quad \leq \E\bigg[\frac{1}{N} \sum_{i=1}^N \sup_{t \in [0,t_*]} \lambda_{\epsilon}\Big(\bm{X}^{i;N}(t) - \widetilde{\bm{X}}^{i;N}(t)\Big) \bigg]
\\
& \quad = \lambda_{\epsilon}(0) + \E\bigg[\frac{1}{N} \sum_{i=1}^N \sup_{t \in [0,t_*]} \Big( \bm{B}^{i;N}(t) + \bm{I}^{i;N}(t) + \bm{R}^{i;N}_0(t) + \bm{R}^{i;N}_1(t) \Big) \bigg]
\\
& \quad \leq 3\epsilon + \big( \|b\|_{W^{1,\infty}} + \|f\|_{W^{1,\infty}} \big) \int_0^{t_*} \E\bigg[\frac{1}{N} \sum_{i=1}^N \min\Big( \big| \bm{X}^{i;N}(s) - \widetilde{\bm{X}}^{i;N}(s) \big|, 1 \Big) \bigg] \ \rd s
\\
& \quad \quad + \bigg( \frac{\|w^N\|_{L^\infty}^2}{N} \|f\|_{L^\infty} t_* \bigg)^{\frac{1}{2}} + \frac{\|w^N\|_{L^\infty}^2}{2\epsilon N} \|f\|_{L^\infty} t_*
\\
& \quad \quad + \kappa_1 \int_0^{t_*} \E\big[ \|\bm{\mu}^N(s,\cdot,\cdot) - \mu(s,\cdot,\cdot)\|_{\Phi_w^{-1}} \big] \rd s + \|w^N - w \|_{L^\infty \to L^1} \|f\|_{L^\infty} t_*.
\end{aligned}
\end{equation*}
We know that $\E \big[ \|\bm{\mu}^N(s,\cdot,\cdot) - \mu(s,\cdot,\cdot)\|_{\Phi_w^{-1}} \big] \to 0$ and $\|w^N - w \|_{L^\infty \to L^1} \to 0$ as $N \to 0$.
Taking the scaling $\epsilon \sim N^{-1/2}$ and applying Gronwall's lemma, we have that, for all $t_* \geq 0$, 
\begin{equation*}
\begin{aligned}
& \E\bigg[\frac{1}{N} \sum_{i=1}^N \sup_{t \in [0,t_*]} \min\Big( \big| \bm{X}^{i;N}(t) - \widetilde{\bm{X}}^{i;N}(t) \big|, 1 \Big) \bigg] \to 0, \ \text{ as } \ N \to \infty.
\end{aligned}
\end{equation*}
\end{proof}

\begin{appendix}
\section{} \label{appen:supplementary}
\subsection{Proof of Proposition~\ref{prop:moments_propagate} (Propagation of moments)}

By choosing $\eta(x) = x^2$ and applying It\^o's lemma, we have
\begin{equation} \label{eqn:Ito_polynomial}
\begin{aligned}
\E \big[ \big|\bm{X}^i_t\big|^2 \big]
& = \E\Bigg[ \big|\bm{X}^i_0\big|^2  + \int_0^t 2 b(\bm{X}^i_{s^-}) \bm{X}^i_{s^-} \;\rd s
\\
& \quad \hphantom{= \E\Bigg[}
+ \sum_{j \neq i} \int_0^t f(\bm{X}^j_{s^-}) \big( \big|\bm{X}^i_{s^-} + w^N_{i,j}/N\big|^2 - \big|\bm{X}^i_{s^-}\big|^2\big) \; \rd s
\\
& \quad \hphantom{= \E\Bigg[}
+ \int_0^t f(\bm{X}^i_{s^-}) \big( - \big|\bm{X}^i_{s^-}\big|^2\big) \; \rd s
\Bigg].
\end{aligned}
\end{equation}
Notice that the last term is non-positive, and $|x+w|^2 - |x|^2 \leq 2|w| |x| + |w|^2$.
Hence, we have
\begin{equation*}
\begin{aligned}
\E \big[ \big|\bm{X}^i_t\big|^2 \big]
& \leq \E \big[ \big|\bm{X}^i_0\big|^2 \big]  + \int_0^t 2\|b\|_{L^\infty} \E \big[ \big|\bm{X}^i_{s^-}\big| \big] \;\rd s
\\
& \quad + \int_0^t \|f\|_{L^\infty} \sum_{j \neq i} \Big( 2 \big|w^N_{i,j}/N\big| \E \big[ \big|\bm{X}^i_{s^-}\big| \big] + \big|w^N_{i,j}/N\big|^2 \Big) \;\rd s
\\
& \leq \E \big[ \big|\bm{X}^i_0\big|^2 \big]  + \int_0^t \Big( \|b\|_{L^\infty}^2 + \E \big[ \big|\bm{X}^i_{s^-}\big|^2 \big] \Big) \;\rd s
\\
& \quad + \int_0^t \|f\|_{L^\infty} \sum_{j \neq i} \Big( 2 \big|w^N_{i,j}/N\big| \E \big[ \big|\bm{X}^i_{s^-}\big| \big] + \big|w^N_{i,j}/N\big|^2 \Big) \;\rd s.
\end{aligned}
\end{equation*}
Taking the average over all indices $i$, we obtain
\begin{equation*}
\begin{aligned}
& \frac{1}{N} \sum_{i = 1}^N \E \big[ \big|\bm{X}^i_t\big|^2 \big]
\\
& \quad \leq \frac{1}{N} \sum_{i = 1}^N \E \big[ \big|\bm{X}^i_0\big|^2 \big]  + \int_0^t \bigg( \|b\|_{L^\infty}^2 + \frac{1}{N} \sum_{i = 1}^N \E \big[ \big|\bm{X}^i_{s^-}\big|^2 \big] \bigg) \;\rd s
\\
& \qquad + \int_0^t \|f\|_{L^\infty} \Big( \max_{1 \leq i,j \leq N}|w^N_{i,j}| \ \frac{1}{N} \sum_{i = 1}^N \E \big[ \big|\bm{X}^i_{s^-}\big|^2 \big] + \big(\max_{1 \leq i,j \leq N}|w^N_{i,j}|\big)^2 /N \Big) \;\rd s.
\end{aligned}
\end{equation*}
We can conclude by Gronwall's lemma, which yields a constant $C(t)$ that only depends on $E_0 = \frac{1}{N} \sum_{i = 1}^N \E \big[ \big|\bm{X}^i_0\big|^2 \big]$, $\|b\|_{L^\infty}$, $\|f\|_{L^\infty}$, $\max_{1 \leq i,j \leq N}|w^N_{i,j}|$, and $t$.


\subsection{Proof of Proposition~\ref{prop:equiv_topology_strong} (Weak-* convergence)} \label{sec:topology}
In this subsection, we prove Proposition~\ref{prop:equiv_topology_strong}, which establishes the equivalence between the weak-* convergence in $\mathcal{M}([0,1] \times \R)$ we use in Theorem~\ref{thm:main} and the convergence in $\Phi_w^{-1}$ defined in Sec.~\ref{sec:metric} and used in Theorem~\ref{thm:main_metric}.
Also, in order to treat the convergence of the initial data in Appendix~\ref{app:initial}, we introduce another metric on $\mathcal{P}([0,1] \times \R)$---a negative Sobolev norm---and show that the three topologies are all equivalent in our setting.

To define the negative Sobolev norm, let us introduce the kernel
$\Lambda: \R \to \R_+$,
\begin{equation*}
\begin{aligned}
\Lambda(x) \defeq \frac{1}{2}\exp(-|x|), \quad \forall x \in \R.
\end{aligned}
\end{equation*}
Then, the negative Sobolev norm, $\|\cdot\|_{H^{-1}(\R)}$, is defined by convolution:
\begin{equation*}
\begin{aligned}
\|f\|_{H^{-1}(\R)}^2 \defeq \int_{\R\times\R} f(x) \Lambda (x - y) f(y) \ \rd x \rd y.
\end{aligned}
\end{equation*}
Since the Fourier transform of $\Lambda$ is
\begin{equation*}
\begin{aligned}
\mathcal{F}(\Lambda)(v) = \int_{\R} \Lambda(x) e^{-2\pi i x v} \, \rd x = \frac{1}{1 + 4\pi^2 v^2}, \quad \forall v \in \R,
\end{aligned}
\end{equation*}
we see that our definition of $H^{-1}(\R)$ is equivalent to the usual definition of $H^{s}(\R)$ in the Fourier sense, i.e.,
\begin{equation*}
\begin{aligned}
\|f\|_{H^{s}(\R)}^2 \defeq \int_{\R} (1 + 4\pi^2 v^2)^s |\mathcal{F}(f)(v)|^2 \ \rd v.
\end{aligned}
\end{equation*}
Note that the convolution by $\Lambda$ can also be seen as the \emph{inverse} operator of $(1 - \partial_x^2)$, since, in the sense of distributions,
\begin{equation*}
\begin{aligned}
(1 - \partial_x^2) \Lambda = \delta, \quad (1 - \partial_x^2) (\Lambda \star f) = f, \quad \forall f \in \mathscr{D}'(\R),
\end{aligned}
\end{equation*}
where $\delta$ denotes the Dirac delta function and $\mathscr{D}'(\R)$ denotes the distributions on $\R$.

All of the above extend to Sobolev spaces on the circle, which we denote $H^{-1}([0,1])$ and $H^1([0,1])$. Replacing the kernel $\Lambda$ by $\tilde \Lambda: [0,1] \to \R_+$,
\begin{equation*}
\begin{aligned}
\tilde \Lambda(\xi) \defeq \frac{1}{2} \sum_{n = -\infty}^{+\infty} \exp(-|\xi + n|), \quad \forall \xi \in [0,1],
\end{aligned}
\end{equation*}
which has a natural periodic extension on $\R$, we define the negative Sobolev norm, $H^{-1}([0,1])$, by convolution:
\begin{equation*}
\begin{aligned}
\|f\|_{H^{-1}([0,1])}^2 \defeq \int_{[0,1]\times[0,1]} f(\xi) \tilde \Lambda (\xi - \zeta) f(\zeta) \ \rd \xi \rd \zeta.
\end{aligned}
\end{equation*}
We recall that, in this work, the interval $[0,1]$ is identified with the circle $\T$. Hence, the Sobolev spaces on the circle, $H^{-1}([0,1])$ and $H^1([0,1])$, should not be confused with Sobolev spaces defined on bounded domains (despite our slight abuse of notation).

Finally, we define the tensorized negative Sobolev structure, $H^{-1}([0,1])\otimes H^{-1}(\R)$, for functions on the product domain $[0,1] \times \R$, through the norm
\begin{equation*}
\begin{aligned}
\|f\|_{H^{-1}([0,1]) \otimes H^{-1}(\R)}^2 & \defeq \int_{[0,1] \times \R} f \ \ [(\tilde \Lambda \otimes \Lambda) \star f] \ \rd \xi \rd x
\\
& = \int_{[0,1] \times \R} \int_{[0,1] \times \R} f(\xi,x) \tilde \Lambda(\xi - \zeta) \Lambda(x-y) f(\zeta,y) \; \rd \xi \rd x \rd \zeta \rd y
\\
& = \sum_{k = -\infty}^{+\infty} \int_{\R} \frac{1}{\big(1 + 4\pi^2 k^2\big)\big(1 + 4\pi^2 v^2\big)} |\mathcal{F}(f)(k,v)|^2 \; \rd v,
\end{aligned}
\end{equation*}
where
\begin{equation*}
    \mathcal{F}(f)(k,v) := \int_{[0,1]\times\R} f(\xi,x) e^{-2\pi i (\xi k + x v)} \, \rd \xi\rd x.
\end{equation*}
In Appendix~\ref{app:initial}, we will also need the norm
\begin{equation*}
\begin{aligned}
\|f\|_{H^{1}([0,1]) \otimes H^{1}(\R)}^2 & \defeq \int_{[0,1] \times \R} \big[ (1 - \partial_x^2) (1 - \partial_\xi^2) f \big] f \ \rd \xi \rd x
\\
& = \sum_{k = -\infty}^{+\infty} \int_{\R} \big(1 + 4\pi^2 k^2\big)\big(1 + 4\pi^2 v^2\big) |\mathcal{F}(f)(k,v)|^2 \; \rd v.
\end{aligned}
\end{equation*}

The norm $H^{-1}([0,1]) \otimes H^{-1}(\R)$ will serve as the third notion of convergence in $\mathcal{P}([0,1] \times \R)$. 
The following proposition, which encompasses Proposition~\ref{prop:equiv_topology_strong}, states that the three types of convergence in $\mathcal{P}([0,1] \times \R)$ considered in this work, namely, the weak-* convergence, the $\Phi^{-1}_w$-convergence, and the convergence in the negative Sobolev norm $H^{-1}([0,1]) \otimes H^{-1}(\R)$ are all equivalent.

\begin{proposition} \label{prop:equiv_topology_strong_plus}

Let $\mu^N$, $N \to \infty$, be a sequence of measures on $[0,1] \times \R$, and $\mu$ be a measure on $[0,1] \times \R$. Assume that $\mu^N(\xi,\cdot) \in \mathcal{P}(\R)$ for a.e. $\xi \in [0,1]$ and $\mu(\xi,\cdot) \in \mathcal{P}(\R)$ for a.e. $\xi \in [0,1]$. Further, let $\epsilon_w : (0,1) \to (0,\infty)$ be a non-decreasing function satisfying $\lim_{r \to 0+} \epsilon_w(r) = 0$.
The following statements are equivalent:

\begin{itemize}
\item[(i)]
For all $\varphi \in C_c([0,1] \times \R)$,
\begin{equation*}
\begin{aligned}
\int_{[0,1] \times \R} \varphi ( \mu^N - \mu ) \;\rd \xi \rd x \to 0, \ \text{ as } \ N \to \infty.
\end{aligned}
\end{equation*}

\item[(ii)] $\|\mu^N - \mu\|_{\Phi_w^{-1}} \to 0$ as $N \to \infty$,
where the metric $\Phi_w^{-1}$ is defined in~\eqref{eqn:weak_distance}-\eqref{eqn:dual_norm}.

\item [(iii)] $\| \mu^N - \mu \|_{H^{-1}([0,1]) \otimes H^{-1}(\R)} \to 0$ as $N \to \infty$.

\end{itemize}

\end{proposition}

The proposition can be deduced from a general lemma:
\begin{lemma} \label{lem:equiv_topology}
Assume $\mu^N$, $N \to \infty$, and $\mu$ as in Proposition~\ref{prop:equiv_topology_strong_plus}. Let $E$ be a Banach space such that $C_c^\infty([0,1] \times \R) \subset E \subset L^\infty([0,1]; C(\R))$ and such that $E$ is compactly embedded in $L^1([0,1]; C(K))$ for any compact set $K$. 
Then, (i) is equivalent to the following:
\begin{itemize}
\item[(i')]
As $N\to\infty$,
\begin{equation*}
\begin{aligned}
\|\mu^N - \mu\|_{E^{-1}} 
\defeq \sup \bigg\{ \int_{[0,1] \times \R} \phi(\xi,x) \Big( \mu^N(\xi,x) - \mu(\xi,x) \Big) \rd \xi \rd x\,\bigg|\; \|\phi\|_{E} \leq 1 \bigg\} \to 0.
\end{aligned}
\end{equation*}
\end{itemize}
\end{lemma}
Lemma~\ref{lem:equiv_topology} directly implies Proposition~\ref{prop:equiv_topology_strong_plus}. Indeed, the equivalence between (i) and (iii) follows immediately from an extension of Morrey's inequality, which implies that $H^1([0,1])\otimes H^1(\R)$ is embedded in $C^{\frac{1}{2}}([0,1]\times \R)$.
The equivalence with (ii) is proved by a Fr\'echet-Kolmogorov argument once one observes that a bounded set in $\Phi_w$ is uniformly bounded, Lipschitz in $x$, and also equicontinuous in $\xi$ by the definition of $\epsilon_w$.

\begin{proof}[Proof of Lemma~\ref{lem:equiv_topology}]

$\text{(i')} \implies \text{(i)}$: The proof is immediate since for any $\varphi \in C_c([0,1] \times \R) \subset E$,
\begin{equation*}
\begin{aligned}
\bigg| \int_{[0,1] \times \R} \varphi ( \mu^N - \mu ) \;\rd \xi \rd x \bigg| \leq \|\varphi\|_{E} \|\mu^N - \mu\|_{E^{-1}}, \ \text{ as } \ N \to \infty.
\end{aligned}
\end{equation*}

$\text{(i)} \implies \text{(i')}$: Up to the extraction of a subsequence (still indexed by $N$), there exists $\{\phi^N\}_{N=1}^\infty$ such that $\| \phi^N \|_{E} \leq 1$ and 
\begin{equation*} 
\begin{aligned}
\lim_{N \to \infty} \int_{[0,1] \times \R} \phi^N ( \mu^N - \mu ) \; \rd \xi \rd x = \limsup_{N \to \infty} \|\mu^N - \mu\|_{E^{-1}}.
\end{aligned}
\end{equation*}
Since $E$ is compactly embedded in $L^1([0,1]; C(K))$ for any compact $K$, we can extract a further subsequence (still indexed by $N$ for readability) and there a $\phi$ such that, for any $R>0$, $\phi^N\to \phi$ in $L^1([0,1];C([-R,R]))$. Since $\mu^N$ and $\mu$ are probability measures for a.e. $\xi$, we have that
\begin{multline*}
\int_{[0,1] \times \R} (\phi^N-\phi) ( \mu^N - \mu ) \; \rd \xi \rd x 
\\
\leq 2 \|\phi^N-\phi\|_{L^1([0,1];C([-R,R]))} + \|\phi^N-\phi\|_{L^\infty([0,1];C(\R))} \int_{[0,1] \times (\R \setminus [-R,R])} (\mu^N + \mu).
\end{multline*}
Passing to the limit $N \to \infty$, the first term vanishes by the convergence $\phi^N \to \phi$.
By Prokhorov’s theorem, if (i) is true, then the sequence $\{\mu^N\} \subset \mathcal{P}([0,1] \times \R)$ is tight; hence, the second term vanishes by letting $R \to \infty$.
Hence,
\begin{multline*}
\limsup_{N \to \infty} \|\mu^N - \mu\|_{E^{-1}}=\lim_{N \to \infty} \int_{[0,1] \times \R} \phi^N ( \mu^N - \mu ) \; \rd \xi \rd x \\
= \lim_{N \to \infty} \int_{[0,1] \times \R} \phi ( \mu^N - \mu ) \; \rd \xi \rd x =0,
\end{multline*}
from (i). (As a remark, (i) itself only implies the last equality for $\phi \in C_c([0,1] \times \R)$. A density argument involving the tightness of $\mu^N$ is needed to prove (i) for general $\phi \in L^1([0,1]; C(\R))$.)
\end{proof}

\subsection{Proof of Proposition~\ref{prop:passing_limit_initial} (Convergence of the initial data)}\label{app:initial}

\begin{proof} [Proof of Proposition~\ref{prop:passing_limit_initial}]

To prove the first part of the proposition, notice that 
\begin{equation*}
\begin{aligned}
\int_{[0,1] \times \R} \bigg( \sum_{i = 1}^N \mu^{i;N}_0(x) \mathbbm{1}_{E^{i;N}}(\xi) \bigg) x^2 \;\rd \xi \rd x \leq \frac{1}{N} \sum_{i = 1}^N \E\big[|\bm{X}^{i;N}_0|^2\big].
\end{aligned}
\end{equation*}
Hence, the sequence of measures is uniformly tight. By Prokhorov's theorem,
up to the extraction of a subsequence, there exists $\mu_0 \in L^\infty([0,1]; \mathcal{M}(\R))$ such that
\begin{equation*}
\begin{aligned}
\bigg( \sum_{i = 1}^N \mu^{i;N}_0 \mathbbm{1}_{E^{i;N}} \bigg) \wsto \mu_ 0, \ \text{ as } \ N \to \infty.
\end{aligned}
\end{equation*}
It is then straightforward to check that $\mu_0(\xi,\cdot) \in \mathcal{P}(\R)$ for a.e. $\xi \in [0,1]$.
By Proposition~\ref{prop:equiv_topology_strong_plus}, the weak-* convergence is equivalent to the $\Phi_w^{-1}$-convergence,
\begin{equation*}
\begin{aligned}
\bigg\| \bigg( \sum_{i = 1}^N \mu^{i;N}_0 \mathbbm{1}_{E^{i;N}} \bigg) - \mu_ 0 \bigg\|_{\Phi_w^{-1}} \to 0, \ \text{ as } \ N \to \infty,
\end{aligned}
\end{equation*}
or the convergence in the negative Sobolev norm $H^{-1}([0,1])\otimes H^{-1}(\R)$,
\begin{equation*}
\begin{aligned}
\bigg\| \bigg( \sum_{i = 1}^N \mu^{i;N}_0 \mathbbm{1}_{E^{i;N}} \bigg) - \mu_ 0 \bigg\|_{H^{-1}([0,1])\otimes H^{-1}(\R)} \to 0, \ \text{ as } \ N \to \infty.
\end{aligned}
\end{equation*}
This concludes the proof of the first part.

Next, we prove the convergence of the extended empirical measure in the $H^{-1}([0,1])\otimes H^{-1}(\R)$ norm (which we denote $H^{-1}$ to lighten the notation). Our strategy starts with applying the triangle inequality
\begin{equation*}
\begin{aligned}
\big\| \bm{\mu}^N_0 - \mu_0 \big\|_{H^{-1}} \leq \bigg\|\bm{\mu}^N_0 - \bigg( \sum_{i = 1}^N \mu^{i;N}_0 \mathbbm{1}_{E^{i;N}} \bigg) \bigg\|_{H^{-1}} + \bigg\|\bigg( \sum_{i = 1}^N \mu^{i;N}_0 \mathbbm{1}_{E^{i;N}} \bigg) - \mu_0 \bigg\|_{H^{-1}}.
\end{aligned}
\end{equation*}
The second term converges to $0$ by our previous discussion. For the first term, we have that,
\begin{equation*}
\begin{aligned}
& \bigg\|\bm{\mu}^N_0 - \bigg( \sum_{i = 1}^N \mu^{i;N}_0 \mathbbm{1}_{E^{i;N}} \bigg) \bigg\|_{H^{-1}}^2
\\
&= \sum_{i = 1}^N \sum_{j = 1}^N \Big\langle \big( \delta_{\bm{X}^{i;N}_0} - \mu^{i;N}_0 \big) \mathbbm{1}_{E^{i;N}}, \big( \delta_{\bm{X}^{j;N}_0} - \mu^{j;N}_0 \big) \mathbbm{1}_{E^{j;N}} \Big\rangle_{H^{-1}}
\\
&= \sum_{i = 1}^N \sum_{j = 1}^N \int_{[0,1] \times \R} \big( \delta_{\bm{X}^{i;N}_0} - \mu^{i;N}_0 \big)(x) \ \mathbbm{1}_{E^{i;N}}(\xi) \ \Lambda \star \big( \delta_{\bm{X}^{j;N}_0} - \mu^{j;N}_0 \big)(x) \ \tilde \Lambda \star \mathbbm{1}_{E^{j;N}}(\xi) \ \rd \xi \rd x.
\end{aligned}
\end{equation*}
Since the $\{\bm{X}^{i;N}_0\}_{i=1}^N$ are independent, the expectation can be computed:
\begin{equation*}
\begin{aligned}
& \E \bigg[ \bigg\|\bm{\mu}^N_0 - \bigg( \sum_{i = 1}^N \mu^{i;N}_0 \mathbbm{1}_{E^{i;N}} \bigg) \bigg\|_{H^{-1}}^2 \bigg] 
\\
&= \E\sum_{i,j} \int_{[0,1] \times \R} \big( \delta_{\bm{X}^{i;N}_0} - \mu^{i;N}_0 \big)(x) \ \mathbbm{1}_{E^{i;N}}(\xi) \ \Lambda \star \big( \delta_{\bm{X}^{j;N}_0} - \mu^{j;N}_0 \big)(x) \ \tilde \Lambda \star \mathbbm{1}_{E^{j;N}}(\xi) \ \rd \xi \rd x 
\\
&= \E\sum_{i = 1}^N  \int_{[0,1] \times \R} \big( \delta_{\bm{X}^{i;N}_0} - \mu^{i;N}_0 \big)(x) \ \mathbbm{1}_{E^{i;N}}(\xi) \ \Lambda \star \big( \delta_{\bm{X}^{i;N}_0} - \mu^{i;N}_0 \big)(x) \ \tilde \Lambda \star \mathbbm{1}_{E^{i;N}}(\xi) \ \rd \xi \rd x ,
\end{aligned}
\end{equation*}
where, for the second equality, we used the fact that the terms with $i \neq j$ have expectation $0$. It is easy to verify that
\begin{equation*}
\begin{aligned}
& \|\Lambda \star \big( \delta_{\bm{X}^{i;N}_0} - \mu^{i;N}_0 \big)\|_{L^\infty} \leq \|\Lambda\|_{L^\infty} \|\delta_{\bm{X}^{i;N}_0} - \mu^{i;N}_0\|_{\mathrm{TV}} \leq 1, && \text{ a.s.},
\\
& \|\tilde \Lambda \star \mathbbm{1}_{E^{i;N}}\|_{L^\infty} \leq \|\tilde \Lambda\|_{L^\infty} \|\mathbbm{1}_{E^{i;N}}\|_{\mathrm{TV}} \leq \frac{2}{N},
\end{aligned}
\end{equation*}
where $\|\cdot\|_{\mathrm{TV}}$ denotes the total variation norm. Therefore,
\begin{equation*}
\begin{aligned}
& \E \bigg[ \bigg\|\bm{\mu}^N_0 - \bigg( \sum_{i = 1}^N \mu^{i;N}_0 \mathbbm{1}_{E^{i;N}} \bigg) \bigg\|_{H^{-1}}^2 \bigg] 
\\
& \quad \leq \E\sum_{i = 1}^N \big\| \delta_{\bm{X}^{i;N}_0} - \mu^{i;N}_0 \big\|_{\mathrm{TV}} \ \|\mathbbm{1}_{E^{i;N}}\|_{\mathrm{TV}} \|\Lambda \star \big( \delta_{\bm{X}^{i;N}_0} - \mu^{i;N}_0 \big)\|_{L^\infty} \|\tilde \Lambda \star \mathbbm{1}_{E^{i;N}}\|_{L^\infty} 
\\
& \quad = \frac{4}{N}.
\end{aligned}
\end{equation*}
Thus, 
\begin{equation*}
\begin{aligned}
\E \Big[ \big\| \bm{\mu}^N_0 - \mu_0 \big\|_{H^{-1}}\Big] \leq \frac{2}{\sqrt{N}} + \bigg\|\bigg( \sum_{i = 1}^N \mu^{i;N}_0 \mathbbm{1}_{E^{i;N}} \bigg) - \mu_0 \bigg\|_{H^{-1}} \to 0, \ \text{ as } \ N \to \infty,
\end{aligned}
\end{equation*}
which implies that
\begin{equation*}
\begin{aligned}
\big\| \bm{\mu}^N_0 - \mu_0 \big\|_{H^{-1}} \to 0, \ \text{ as } \ N \to \infty, \quad a.s.
\end{aligned}
\end{equation*}
By Proposition~\ref{prop:equiv_topology_strong_plus}, this convergence is equivalent to
\begin{equation*}
\begin{aligned}
\big\| \bm{\mu}^N_0 - \mu_0 \big\|_{\Phi_w^{-1}} \to 0, \ \text{ as } \ N \to \infty, \quad a.s.
\end{aligned}
\end{equation*}
Since the metric $\Phi_w^{-1}$ is bounded a.s. by definition, we conclude that
\begin{equation*}
\begin{aligned}
\E \Big[ \big\| \bm{\mu}^N_0 - \mu_0 \big\|_{\Phi_w^{-1}}\Big] \to 0, \ \text{ as } \ N \to \infty.
\end{aligned}
\end{equation*}
\end{proof}
\end{appendix}

\section*{Acknowledgments}

P.-E. J. is also affiliated with the Huck Institutes of the Life Sciences at Pennsylvania State University.  
D. Z. is also affiliated with the Laboratoire Jacques-Louis Lions at Sorbonne Université.

\begin{funding}
P.-E. J. and D. Z. were partially supported by NSF DMS Grants 2205694, 2219297.
D. Z. was supported by European Union’s Horizon 2020 research and innovation programme under the Marie Skłodowska-Curie grant agreement No. 101034255.
V. S. was supported by the Swiss National Science Foundation (no 200020\_207426) and a Royal Society Newton International Fellowship (NIF\textbackslash R1\textbackslash 231927).
\end{funding}

\bibliography{mybib_BW}{} 
\bibliographystyle{siam} 

\end{document}